\documentclass[11pt,oneside]{amsart}

\usepackage{tikz-cd}

\usepackage{newtxtext}
\usepackage{newtxmath}
\usepackage{fix-cm}

\DeclareMathAlphabet{\mathcal}{OMS}{cmsy}{m}{n}

\usepackage[arrow,curve,matrix,tips,2cell]{xy}

\usepackage[hmargin=2cm,vmargin=1.5cm]{geometry}

\usepackage{bm}

\usepackage{fancyhdr}

\usepackage[obeyspaces]{url}

\usepackage{tensor}

\usepackage{mathtools}

\usepackage{accents}

\usepackage{amsmath}
\usepackage{amscd}
\usepackage{amssymb}
\usepackage{latexsym}
\usepackage{url}

\usepackage{graphicx}

\newtheorem{theorem}{Theorem}[section]
\newtheorem*{theorem*}{Theorem}
\newtheorem{lemma}[theorem]{Lemma}
\newtheorem*{lemma*}{Lemma}
\newtheorem{corollary}[theorem]{Corollary}
\newtheorem{proposition}[theorem]{Proposition}

\newtheorem{remark}[theorem]{Remark}
\newtheorem{definition}[theorem]{Definition}
\newtheorem*{definition*}{Definition}

\newtheorem{question}[theorem]{Question}
\newtheorem*{question*}{Question}

\newtheorem{example}[theorem]{Example}
\newtheorem{examples}[theorem]{Examples}

\makeatletter
\def\revddots{\mathinner{\mkern1mu\raise\p@
\vbox{\kern7\p@\hbox{.}}\mkern2mu
\raise4\p@\hbox{.}\mkern2mu\raise7\p@\hbox{.}\mkern1mu}}
\makeatother 
\newcommand{\bgl}{\begin{equation}} 
\newcommand{\egl}{\end{equation}}
\newcommand{\bgloz}{\begin{equation*}} 
\newcommand{\egloz}{\end{equation*}}
\newcommand{\bgln}{\begin{eqnarray}} 
\newcommand{\egln}{\end{eqnarray}}
\newcommand{\bglnoz}{\begin{eqnarray*}} 
\newcommand{\eglnoz}{\end{eqnarray*}}
\newcommand{\btheo}{\begin{theorem}}
\newcommand{\etheo}{\end{theorem}}
\newcommand{\btheooz}{\begin{theorem*}}
\newcommand{\etheooz}{\end{theorem*}}

\newcommand{\blemma}{\begin{lemma}}
\newcommand{\elemma}{\end{lemma}}
\newcommand{\blemmaoz}{\begin{lemma*}}
\newcommand{\elemmaoz}{\end{lemma*}}
\newcommand{\bproof}{\begin{proof}}
\newcommand{\eproof}{\end{proof}}
\newcommand{\bbew}{\begin{beweis}}
\newcommand{\ebew}{\end{beweis}}
\newcommand{\bremark}{\begin{remark}\em}
\newcommand{\eremark}{\end{remark}}
\newcommand{\bdefin}{\begin{definition}}
\newcommand{\edefin}{\end{definition}}
\newcommand{\bdefinoz}{\begin{definition*}}
\newcommand{\edefinoz}{\end{definition*}}
\newcommand{\bex}{\begin{example}\em}
\newcommand{\eex}{\end{example}}
\newcommand{\bexs}{\begin{examples}}
\newcommand{\eexs}{\end{examples}}

\newcommand{\bprop}{\begin{proposition}}
\newcommand{\eprop}{\end{proposition}}
\newcommand{\bcor}{\begin{corollary}}
\newcommand{\ecor}{\end{corollary}}
\newcommand{\bfa}{\begin{cases}} 
\newcommand{\efa}{\end{cases}}

\newcommand{\bquestion}{\begin{question}}
\newcommand{\equestion}{\end{question}}
\newcommand{\bquestionoz}{\begin{question*}}
\newcommand{\equestionoz}{\end{question*}}
\newcommand{\cE}{\mathcal E}
\newcommand{\cF}{\mathcal F}

\newcommand{\cM}{\mathcal M}

\newcommand{\cO}{\mathcal O}

\newcommand{\cW}{\mathcal W}
\newcommand{\cX}{\mathcal X}
\newcommand{\cY}{\mathcal Y}
\newcommand{\cZ}{\mathcal Z}

\def\Cz{\mathbb{C}}

\def\Zz{\mathbb{Z}}

\def\1z{\mathbb{1}}
\newcommand{\fC}{\mathfrak C}
\newcommand{\fD}{\mathfrak D}

\newcommand{\fI}{\mathfrak I}

\newcommand{\fX}{\mathfrak X}
\newcommand{\fY}{\mathfrak Y}
\newcommand{\mfd}{\mathfrak d}
\newcommand{\mfe}{\mathfrak e}
\newcommand{\mff}{\mathfrak f}
\newcommand{\mfm}{\mathfrak m}
\newcommand{\mfn}{\mathfrak n}
\newcommand{\mfp}{\mathfrak p}
\newcommand{\mfq}{\mathfrak q}
\newcommand{\mfr}{\mathfrak r}
\newcommand{\mfs}{\mathfrak s}
\newcommand{\mft}{\mathfrak t}
\newcommand{\an}[1]{``#1''} 
\newcommand{\ti}{\tilde}

\newcommand{\lori}{\longrightarrow}

\newcommand{\ma}{\mapsto} 
\newcommand{\onto}{\twoheadrightarrow} 
\newcommand{\into}{\hookrightarrow} 
\newcommand{\tailarr}{\rightarrowtail} 
\newcommand{\isom}{\xrightarrow{\raisebox{-1ex}[0ex][0ex]{$\sim$}}} 
\newcommand{\Rarr}{\Rightarrow} 
\newcommand{\Larr}{\Leftarrow} 

\newcommand{\ve}{\varepsilon}

\def\SEMI{\mbox{$\times\kern-2pt\vrule height5pt width.6pt \kern3pt $}}

\newcommand{\halb}{\tfrac{1}{2}}

\newcommand{\Spec}{{\rm Spec\,}} 

\newcommand{\id}{{\rm id}}

\newcommand{\Ad}{{\rm Ad\,}}

\renewcommand{\dim}{{\rm dim}\,}
\newcommand{\rk}{{\rm rk}\,}
\newcommand{\img}{{\rm im\,}}
\renewcommand{\ker}{{\rm ker}\,}

\newcommand{\reg}{^\times} 
\newcommand{\ev}{\operatorname{ev}} 
\newcommand{\defeq}{\mathrel{:=}} 

\newcommand{\dop}{\text{: }} 
\newcommand{\ilim}{\varinjlim} 
\newcommand{\plim}{\varprojlim} 
\newcommand{\lge}{\left\{} 
\newcommand{\rge}{\right\}} 
\newcommand{\lru}{\left(} 
\newcommand{\rru}{\right)} 
\newcommand{\lsp}{\left\langle} 
\newcommand{\rsp}{\right\rangle} 
\newcommand{\rukl}[1]{\lru #1 \rru} 
\newcommand{\gekl}[1]{\lge #1 \rge} 
\newcommand{\spkl}[1]{\lsp #1 \rsp} 
\newcommand{\menge}[2]{\gekl{ #1 \dop #2 }} 
\def\bf1{\mathbf{1}}
\newcommand{\half}{\tfrac{1}{2}}
\newcommand{\dom}{{\rm dom\,}}
\newcommand{\bfdot}[1]{\accentset{\bullet}{#1}}
\newcommand{\accir}[1]{\accentset{\circ}{#1}}
\newcommand{\pc}{\sim_{\rm conn}}

\begin{document}

\title{Constructing Menger manifold C*-diagonals in classifiable C*-algebras}

\thispagestyle{fancy}

\author{Xin Li}

\address{Xin Li, School of Mathematics and Statistics, University of Glasgow, University Place, Glasgow G12 8QQ, United Kingdom}
\email{Xin.Li@glasgow.ac.uk}


\thanks{This project has received funding from the European Research Council (ERC) under the European Union's Horizon 2020 research
and innovation programme (grant agreement No. 817597).}

\begin{abstract}
We construct C*-diagonals with connected spectra in all classifiable stably finite C*-algebras which are unital or stably projectionless with continuous scale. For classifiable stably finite C*-algebras with torsion-free $K_0$ and trivial $K_1$, we further determine the spectra of the C*-diagonals up to homeomorphism. In the unital case, the underlying space turns out to be the Menger curve. In the stably projectionless case, the space is obtained by removing a non-locally-separating copy of the Cantor space from the Menger curve. We show that each of our classifiable C*-algebras has continuum many pairwise non-conjugate such Menger manifold C*-diagonals. Along the way, we also obtain a complete classification of C*-diagonals in all one-dimensional non-commutative CW complexes.
\end{abstract}

\maketitle


\setlength{\parindent}{0cm} \setlength{\parskip}{0.5cm}

\section{Introduction}

Classification of C*-algebras is a research programme initiated by the work of Glimm, Dixmier, Bratteli and Elliott. After some recent breakthroughs, the combination of work of many many mathematicians over several decades has culminated in the complete classification of unital separable simple nuclear $\cZ$-stable C*-algebras satisfying the UCT (see \cite{KP, Phi, GLN, GLNa, GLNb, EGLN, TWW} and the references therein). Further classification results are expected to cover the stably projectionless case as well (see for instance \cite{EN,EGLN17a,EGLN17b,GLI,GLII}). All in all, the final result classifies all separable simple nuclear $\cZ$-stable C*-algebras satisfying the UCT (which we refer to as \an{classifiable C*-algebras} in this paper) by their Elliott invariants.

Recently, it was shown in \cite{Li18} that every classifiable C*-algebra has a Cartan subalgebra. The interest here stems from the observation in \cite{Kum,Ren} that once a Cartan subalgebra has been found, it automatically produces an underlying topological groupoid such that the ambient C*-algebra can be written as the corresponding groupoid C*-algebra. Therefore, the results in \cite{Li18} build a strong connection between classification of C*-algebras and generalized topological dynamics (in the form of topological groupoids and their induced orbit structures). This connection has already proven to be very fruitful, for instance in the classification of Cantor minimal systems up to orbit equivalence \cite{GPS, GMPS08, GMPS10} or in the context of approximation properties \cite{Ker,KS}. Generally speaking, the notion of Cartan subalgebras in C*-algebras has attracted attention recently due to links to topological dynamics \cite{Li16,Li17,Li_DQH} and the UCT question \cite{BL16,BL17}.

More precisely, the construction in \cite{Li18} produces Cartan subalgebras in all the C*-algebra models from \cite{Ell, EV, Tho, GLN, GLNa} which exhaust all possible Elliott invariants of classifiable stably finite C*-algebras. Actually, we obtain C*-diagonals in this case (i.e., the underlying topological groupoid has no non-trivial stabilizers). Together with groupoid models (and hence Cartan subalgebras) which have already been constructed in the purely infinite case (see \cite{Spi} and also \cite[\S~5]{LR}), this produces Cartan subalgebras in all classifiable C*-algebras. An alternative approach to constructing groupoid models, based on topological dynamics, has been developed in \cite{DPS18,Put,DPS19a,DPS19b} and covers large classes of classifiable C*-algebras. In special cases, groupoid models have also been constructed in \cite{AM}.

The goal of this paper is to start a more detailed analysis of the C*-diagonals and the corresponding groupoids constructed in \cite{Li18}. A motivating question is whether the construction in \cite{Li18} produces a one-dimensional C*-diagonal of the Jiang-Su algebra $\cZ$ which is distinguished in some sense (here one-dimensional C*-diagonal means C*-diagonal whose spectrum has covering dimension one), or, put differently (see \cite[Problem~3]{BS_MFO}):
\bquestion
\label{q:UniqueCartanZ}
Does the Jiang-Su algebra $\cZ$ have any distinguished (one-dimensional) Cartan subalgebras?
\equestion 
Note that such uniqueness questions cannot have an affirmative answer without restrictions such as bounds on the dimension because every classifiable C*-algebra is $\cZ$-stable, so that taking tensor products produces Cartan subalgebras whose spectra have arbitrarily large covering dimension (see \cite[Proposition~5.1]{LR}). Instead of fixing the covering dimension, an even stronger restriction would be to fix the homeomorphism type of the spectrum and to look for a unique or distinguished Cartan subalgebra whose spectrum coincides with a given topological space. This leads to the question what we can say about the spectra of Cartan subalgebras of classifiable C*-algebras. In general, not much is known. Before the work in \cite{Li18}, it was for instance not even known whether $\cZ$ has any Cartan subalgebra with one-dimensional spectrum. Another example is the following question (see \cite[Problem~11]{BS_MFO}):
\bquestion
\label{q:UHF_connCartan}
Does the CAR algebra have a Cartan subalgebra with connected spectrum?
\equestion
\setlength{\parindent}{0cm} \setlength{\parskip}{0cm}

This question is motivated by a construction, due to Phillips and Wassermann \cite{PW_MFO}, of uncountably many pairwise non-conjugate MASAs (which are not Cartan subalgebras) in the CAR algebra whose spectra are all homeomorphic to the unit interval. In this context, we would like to mention that Kumjian \cite{Kum88} had constructed a Cartan subalgebra in an AF algebra with spectrum homeomorphic to the unit circle.
\setlength{\parindent}{0cm} \setlength{\parskip}{0.5cm}

The following are the main results of this paper, which shed some light on the above-mentioned questions.
\btheo
\label{thm:main1}
Every classifiable stably finite C*-algebra which is unital or stably projectionless with continuous scale (in the sense of \cite{Lin91,Lin04,GLI,GLII}) has a C*-diagonal with connected spectrum.
\etheo

\btheo
\label{thm:main2_unital}
Every classifiable stably finite unital C*-algebra with torsion-free $K_0$ and trivial $K_1$ has continuum many pairwise non-conjugate C*-diagonals whose spectra are all homeomorphic to the Menger curve. 
\etheo
\setlength{\parindent}{0cm} \setlength{\parskip}{0cm}

The Menger curve is also known as Menger universal curve, Menger cube, Menger sponge, Sierpinski cube or Sierpinski sponge. It was constructed by Menger \cite{Men} as a universal one-dimensional space, in the sense that every separable metrizable space of dimension at most one embeds into it. Anderson \cite{And58_1,And58_2} characterized the Menger curve by abstract topological properties. The reader may consult \cite{MOT} for more information about the Menger curve, including a concrete construction. 
\setlength{\parindent}{0cm} \setlength{\parskip}{0.5cm}

In order to obtain a version of Theorem~\ref{thm:main2_unital} in the stably projectionless setting, we need to replace the Menger curve $\bm{M}$ by another Menger manifold (a topological space locally homeomorphic to $\bm{M}$) of the form $\bm{M} \setminus \iota(C)$, where $\iota$ is an embedding of the Cantor space $C$ into $\bm{M}$ such that $\iota(C)$ is a non-locally-separating subset of $\bm{M}$, in the sense that for every connected open subset $U$ of $\bm{M}$, $U \setminus \iota(C)$ is still connected. Up to homoemorphism, the space $\bm{M} \setminus \iota(C)$ does not depend on the choice of $\iota$ (see \cite{MOT}), and we denote this space by $\bm{M}_{\setminus C} \defeq \bm{M} \setminus \iota(C)$. 
\btheo
\label{thm:main2_spl}
Every classifiable stably projectionless C*-algebra with continuous scale, torsion-free $K_0$ and trivial $K_1$ has continuum many pairwise non-conjugate C*-diagonals whose spectra are all homeomorphic to $\bm{M}_{\setminus C}$. 
\etheo
Theorem~\ref{thm:main1} answers Question~\ref{q:UHF_connCartan}. Note that in the stably projectionless case, the absence of projections only guarantees the absence of compact open subsets in the spectrum, but it does not automatically lead to a single connected component (see \cite[\S~8]{Li18}). Theorems~\ref{thm:main2_unital} and \ref{thm:main2_spl} show that the uniqueness question for Cartan subalgebras in classifiable C*-algebras has a negative answer unless we impose further conditions. Hence the problem asking for classification of Cartan subalgebras in classifiable C*-algebras (see \cite[Problem~3]{BS_MFO}) seems to be challenging (which maybe makes it interesting). It is interesting to point out that the situation for classifiable C*-algebras is very different from the corresponding one for von Neumann algebras. Theorems~\ref{thm:main2_unital} and \ref{thm:main2_spl} also tell us that in general, it seems that there is not much we can say about the induced map on K-theory of the natural inclusion of a Cartan subalgebra (see Remark~\ref{rem:KMenger}, which sheds some light on \cite[Problem~8]{BS_MFO}). 
\setlength{\parindent}{0.5cm} \setlength{\parskip}{0cm}

In particular, Theorem~\ref{thm:main2_unital} applies to all infinite-dimensional unital separable simple AF algebras, for instance all UHF algebras, and to $\cZ$. Theorem~\ref{thm:main2_spl} applies in particular to the Razak-Jacelon algebra $\cW$ and the stably projectionless version $\cZ_0$ of the Jiang-Su algebra of \cite[Definition~7.1]{GLII}. Even restricted to these special cases, Theorems~\ref{thm:main2_unital} and \ref{thm:main2_spl} yield new results (and we do not need the full strength of the classification theorem for all classifiable C*-algebras; the results in for instance \cite{Rob} suffice).
\setlength{\parindent}{0cm} \setlength{\parskip}{0.5cm}

The constructions we develop in order to prove our main theorems work in general, but only produce C*-diagonals with the desired properties under the conditions we impose in our main theorems. There are several reasons: In \cite{Li18}, C*-diagonals are constructed in all classifiable C*-algebras using the method of cutting down by suitable elements. This procedure, however, might not preserve connectedness. This is why Theorem~\ref{thm:main1} only covers unital C*-algebras and stably projectionless C*-algebras with continuous scale. Note that, however, this class of C*-algebras covers all classifiable C*-algebras up to stable isomorphism. The reason we further restrict to the case of torsion-free $K_0$ and trivial $K_1$ in Theorems~\ref{thm:main2_unital} and \ref{thm:main2_spl} is twofold: It is shown in \cite{Li18} that the spectra of the C*-diagonals constructed in \cite{Li18} will have dimension at least two as soon as torsion appears in K-theory. This rules out $\bm{M}$ or $\bm{M}_{\setminus C}$ as the spectrum in general. Even more serious is the obstruction that the path-lifting property established in Proposition~\ref{prop:path} for the connecting maps at the groupoid level, which plays a crucial role in establishing Theorems~\ref{thm:main2_unital} and \ref{thm:main2_spl}, does not hold anymore in the case where $K_0$ contains torsion or $K_1$ is non-trivial.

In order to prove our main results, the strategy is to adjust the constructions of C*-algebra models in \cite{Ell, EV, Tho, GLN, GLNa}, which arise as inductive limits of simpler building blocks and which exhaust all possible Elliott invariants, in such a way that the new, modified constructions produce C*-algebra models with C*-diagonals having various desired properties. The reader may find the corresponding versions of our main results in Theorems~\ref{thm:conn_Ell}, \ref{thm:ManyMenger_GPD_Ell} and \ref{thm:ManyMenger_Diag_Ell}, which do not depend on general classification results for all classifiable C*-algebras. These versions in combination with general classification results then yield our main theorems as stated above. To construct Cartan subalgebras in inductive limit C*-algebras, an important tool has been developed in \cite{Li18}. However, in \cite{Li18}, we were merely interested in existence results for C*-diagonals, whereas the present work requires several further modifications as well as a finer analysis of the construction of C*-diagonals in \cite{Li18} in order to ensure topological properties of the spectrum such as connectedness as well as abstract topological properties characterizing $\bm{M}$ or further properties characterizing $\bm{M}_{\setminus C}$. At the technical level, a crucial role is played by a new path-lifting property (see Proposition~\ref{prop:path}) of the connecting maps at the groupoid level. This is particularly powerful in combination with inverse limit descriptions of the spectra of the C*-diagonals we construct. Further fine-tuning of the construction is required to produce C*-diagonals for which we can completely determine the spectra up to homeomorphism. In order to show that the construction yields continuum many pairwise non-conjugate C*-diagonals, the key idea is to exploit connectedness not only of the spectra but of (parts of) the groupoid models themselves. This aspect of the construction seems to be interesting on its own right, because many important groupoid models which have been previously studied (for instance for AF algebras, Kirchberg algebras or coming from Cantor minimal systems) have totally disconnected unit spaces.

Important building blocks leading to the C*-algebra models in \cite{Ell, EV, Tho, GLN, GLNa} are given by one-dimensional non-commutative CW complexes and their generalizations. Therefore, as a starting point, we develop a complete classification of C*-diagonals in one-dimensional non-commutative CW complexes. Roughly speaking, the conjugacy class of C*-diagonals in these building blocks encode a particular set of data which can be used to construct the ambient non-commutative CW complex and which we can view as a one-dimensional CW complex in the classical sense (i.e., a graph). We refer to Theorem~\ref{thm:ClassAB} for more details. Our classification theorem generalizes the corresponding results for C*-diagonals in dimension drop algebras in \cite{BR}. It also puts into context the observation in \cite{BR} that in special cases, these C*-diagonals are classified up to conjugacy by the homeomorphism type of their spectra (see Theorem~\ref{thm:AB-B} for a generalization and Remark~\ref{rem:appBR}, Example~\ref{ex:appBR} for a clarification). 
\setlength{\parindent}{0.5cm} \setlength{\parskip}{0cm}

Several of the ideas and techniques leading to our main theorems already feature in the discussion of C*-diagonals in one-dimensional noncommutative CW complexes. However, even though a good understanding of these C*-diagonals played an important role in developing our main results, the actual classification results for this class of C*-diagonals are not needed in the proofs of Theorems~\ref{thm:main1}, \ref{thm:main2_unital} and \ref{thm:main2_spl}.
\setlength{\parindent}{0cm} \setlength{\parskip}{0.5cm}

\section{Classification of C*-diagonals in 1-dimensional NCCW complexes}
\label{s:nccw}

We set out to classify C*-diagonals in 1-dimensional non-commutative CW (NCCW) complexes up to conjugacy. The reader may find more about NCCW complexes in \cite{ET, Ell, ELP, Dea, Rob}. Let us start by introducing notations and some standing assumptions. Throughout this section, $\beta_0, \, \beta_1: \: F \to E$ denote *-homomorphisms between finite-dimensional C*-algebras $F$ and $E$. Let $F = \bigoplus_{i \in I} F^i$ and $E = \bigoplus_{p \in P} E^p$ denote the decompositions of $F$ and $E$ into matrix algebras and $DF^i$, $DE^p$ the canonical C*-diagonals of diagonal matrices. The 1-dimensional NCCW complex $A = A(E,F,\beta_0,\beta_1)$ is given by $A = \menge{(f,a) \in C([0,1],E) \oplus F}{f(\mfr) = \beta_\mfr(a) \ \rm{for} \ \mfr = 0,1}$. For $\mfr = 0,1$, we write $\beta_\mfr^p$ for the composition $F \overset{\beta_\mfr}{\lori} E \onto E^p$ where the second map is the canonical projection. We also write $\beta_\mfr^{p,i} \defeq \beta_\mfr^p \vert_{F^i}$ for the restriction of $\beta_\mfr^p$ to $F^i \subseteq F$. Throughout this section, we make the following assumptions:
\begin{enumerate}
\item[(A1)] For all $i$, $p$ and $\mfr = 0,1$, $\beta_\mfr^{p,i}$ is given by the composition
\begin{equation}
\label{e:beta=}
  F^i \overset{1 \otimes \id_{F^i}}{\lori} 1_{m_\mfr(p,i)} \otimes F^i \subseteq M_{m_\mfr(p,i)} \otimes F^i \tailarr E^p.
\end{equation}
\item[(A2)] $(\beta_0, \beta_1): \: F \to E \oplus E$ is injective.
\end{enumerate}
In \eqref{e:beta=}, an arrow $\tailarr$ denotes a *-homomorphism of multiplicity $1$, i.e., which preserves ranks of projections, and which sends diagonal matrices to diagonal matrices (in our case $DM_{m_\mfr(p,i)} \otimes DF^i$ to $DE^p$). Note that (A1) implies that $\beta_\mfr$ sends $DF$ to $DE$.
\setlength{\parindent}{0.5cm} \setlength{\parskip}{0cm}

There is no loss of generality assuming (A1) and (A2): Up to unitary equivalence, every *-homomorphism $F \to E$ is of the form as in \eqref{e:beta=}, so that we can always replace $\beta_\mfr^{p,i}$ by a map of the form \eqref{e:beta=} without changing the isomorphism class of $A$. And if (A2) does not hold, then $A$ decomposes as $A = A' \oplus F'$ where $A'$ is a 1-dimensional NCCW complex for which (A2) holds and $F' = \ker(\beta_0,\beta_1)$. Then the study of C*-diagonals in $A' \oplus F'$ reduces to the study of C*-diagonals in $A'$ and $F'$, and C*-diagonals in $F'$ are well-understood.

(A2) allows us to identify $A$ with the sub-C*-algebra $\menge{f \in C([0,1],E)}{(f(0), f(1)) \in \img(\beta_0,\beta_1)}$ of $C([0,1],E)$. We will do so frequently without explicitly mentioning it.
\setlength{\parindent}{0cm} \setlength{\parskip}{0.5cm}

Before we start to develop our classification results, we give an overview. If we let $\cX^i \defeq \Spec DF^i$, $\cX \defeq \Spec DF$ and $\cY^p \defeq \Spec DE^p$, $\cY \defeq \Spec DE$, then for $\mfr = 0, 1$, $\beta_{\mfr}$ corresponds to a collection $(\bm{b}_\mfr^p)_p$ of maps $\bm{b}_\mfr^p: \: \cY_\mfr^p \to \cX$ for some $\cY_\mfr^p \subseteq \cY^p$. Viewing $\cY^p$ as edges, $\cX$ as vertices and $\bm{b}_0^p$, $\bm{b}_1^p$ as source and target maps, this data gives rise to a collection of directed graphs $\Gamma^p$, or 1-dimensional CW complexes in the classical sense. (Strictly speaking, this is only correct when $A$ is unital; in the non-unital case, we obtain non-compact 1-dimensional CW complexes obtained by removing finitely many points from compact 1-dimensional CW complexes.)
Moreover, given a permutation $\bm{\sigma} = \coprod \bm{\sigma}^p$ of $\cY = \coprod \cY^p$, we obtain twisted graphs $\Gamma_{\bm{\sigma}}^p$ with the same edge set $\cY^p$, the same vertex set $\cX$, the same source map $\bm{b}_0^p$ and twisted target map $\bm{b}_1^p \circ \bm{\sigma}^p$.
Now it turns out that every C*-diagonal of a 1-dimensional NCCW complex corresponds to a permutation $\bm{\sigma}$ as above, and for two such permutations $\bm{\sigma}$ and $\bm{\tau}$, the corresponding C*-diagonals are conjugate if and only if the collections of oriented graphs $(\Gamma_{\bm{\sigma}}^p)_p$ and $(\Gamma_{\bm{\tau}}^p)_p$ are isomorphic in the sense that there exist isomorphisms of the individual graphs which are either orientation-preserving or orientation-reversing for each $p$. We refer to Theorem~\ref{thm:ClassAB} for more details.

As a first step, we provide models for C*-diagonals in $A$ up to conjugacy. Given a permutation matrix $\sigma$ in $E$, set 
$$
  A_{\sigma} 
  \defeq A(E,F,\beta_0,\Ad(\sigma) \circ \beta_1) 
  = \menge{(f,a) \in C([0,1],E) \oplus F}{f(0) = \beta_0(a), \, f(1) = \sigma \beta_1(a) \sigma^*}.
$$  
Moreover, define
$$
  B_{\sigma} \defeq \menge{(f,a) \in A_{\sigma}}{f(t) \in DE \ \forall \, t \in [0,1]}.
$$
Note that given $(f,a) \in A_{\sigma}$, the condition $f(t) \in DE$ for all $t \in [0,1]$ implies $a \in DF$ by (A1) and (A2). The following observation is a straightforward generalization of \cite[Proposition~5.1]{BR}.
\setlength{\parindent}{0cm} \setlength{\parskip}{0cm}

\blemma
$B_{\sigma}$ is a C*-diagonal of $A_{\sigma}$.
\elemma
\setlength{\parindent}{0cm} \setlength{\parskip}{0.5cm}

Conversely, it turns out that up to conjugacy, every C*-diagonal of $A$ is of this form.
\bprop
\label{prop:AB=AsBs}
For every C*-diagonal $B$ of $A$, there exists a permutation matrix $\sigma \in E$ such that $(A,B) \cong (A_{\sigma},B_{\sigma})$, i.e., there exists an isomorphism $A \isom A_{\sigma}$ sending $B$ onto $B_{\sigma}$.
\eprop
\setlength{\parindent}{0cm} \setlength{\parskip}{0cm}

\bproof
For a subset $S \subseteq [0,1]$, let $A_S \defeq \menge{f \vert_S}{f \in A} \subseteq C(S,E)$ and $B_S \defeq \menge{f \vert_S}{f \in B} \subseteq A_S$. It is easy to see (compare \cite[Proposition~4.1]{BR}) that for every $t \in (0,1)$, $B_{\gekl{t}}$ is a C*-diagonal of $A_{\gekl{t}} = E$, and that $B_{\gekl{0,1}}$ is a C*-diagonal of $A_{\gekl{0,1}}$. By (A2), $(\beta_0,\beta_1)$ defines an isomorphism $F \isom A_{\gekl{0,1}}$. Hence $(\beta_0,\beta_1)^{-1}(B_{\gekl{0,1}})$ is a C*-diagonal of $F$. Thus there is a unitary $u_F \in U(F)$ such that $u_F (\beta_0,\beta_1)^{-1}(B_{\gekl{0,1}}) u_F^* = DF$. Applying $(\beta_0,\beta_1)$ on both sides, we get 
$$
  (\beta_0(u_F),\beta_1(u_F)) (B_{\gekl{0,1}}) (\beta_0(u_F),\beta_1(u_F))^* = (\beta_0,\beta_1)(DF) \subseteq DE \oplus DE.
$$
Here we used that (A1) implies $\beta_\mfr(DF) \subseteq DE$ for $\mfr = 0,1$. Therefore, for $\mfr = 0,1$, $u_\mfr \defeq \beta_\mfr(u_F) + (1_E - \beta_\mfr(1_F))$ is a unitary in $E$ such that $u_\mfr B_{\gekl{\mfr}} u_\mfr^* = \beta_\mfr(u_F) B_{\gekl{\mfr}} \beta_\mfr(u_F)^* \subseteq DE$.
\setlength{\parindent}{0cm} \setlength{\parskip}{0.5cm}

Using \cite[Corollary~2.5 and Lemma~3.4]{BR}, it is straightforward to find $u: \: [0,\half] \to U(E)$ with $u(0) = u_0$ and $u \vert_{(0,1/2]} \in C((0,\half],U(E))$ such that $\Ad(u)$ induces an isomorphism $A_{[0,1/2]} \isom A_{[0,1/2]}$ sending $B_{[0,1/2]}$ to $\menge{f \in A_{[0,1/2]}}{f(t) \in DE \ \forall \, t \in [0,\half]}$. Similarly, find $v: \: [\half,1] \to U(E)$ satisfying $v(1) = u_1$ and $v \vert_{[1/2,1)} \in C([\half,1),U(E))$ such that $\Ad(v)$ induces $A_{[1/2,1]} \isom A_{[1/2,1]}$ sending $B_{[1/2,1]}$ to $\menge{f \in A_{[1/2,1]}}{f(t) \in DE \ \forall \, t \in [\half,1]}$.

Now consider $\sigma = u(\half) v(\half)^* \in U(E)$. We have $\sigma DE \sigma^* = u(\half) v(\half)^* DE v(\half) u(\half)^* = u(\half) B_{\gekl{1/2}} u(\half)^* = DE$. Thus $\sigma$ normalizes $DE$. This implies that $\sigma$ is the product of a unitary in $DE$ and a permutation matrix in $E$. By multiplying $u$ by a suitable unitary in $C([0,\half],U(DE))$, we can arrange that $\sigma$ is given by a permutation matrix in $E$. Define $w: \: [0,1] \to U(E)$ by $w(t) \defeq u(t)$ for $t \in [0,\half]$ and $w(t) \defeq \sigma v(t)$ for $t \in [\half,1]$. Then $w(t) B_{\gekl{t}} w(t)^* \subseteq DE$ for all $t \in [0,1]$. Hence $\Ad(w)$ induces an isomorphism $A \isom A_{\sigma}, \, (f,a) \ma (wfw^*, u_F a u_F^*)$ sending $B$ to $B_{\sigma} = \menge{f \in A_{\sigma}}{f(t) \in DE \ \forall \, t \in [0,1]}$.
\eproof
\setlength{\parindent}{0cm} \setlength{\parskip}{0.5cm}

By Proposition~\ref{prop:AB=AsBs}, the classification problem for C*-diagonals in $A$ reduces to the classification problem for Cartan pairs of the form $(A_{\sigma}, B_{\sigma})$. Our next goal is to further reduce to the situation where no index in $P$ is redundant. Let $A = A(E,F,\beta_0,\beta_1)$ be a 1-dimensional NCCW complex and $B = \menge{f \in A}{f(t) \in DE \ \forall \, t \in [0,1]}$. 
\bdefin
An index $q \in P$ is called redundant if there exists $\bar{q} \in P$ with $\bar{q} \neq q$ and $j \in I$, $\mfr, \mfs \in \gekl{0,1}$ such that $\beta_\mfr^{\bar{q},j}$ and $\beta_\mfs^{q,j}$ are isomorphisms and $\beta_\bullet^{p,j} = 0$ for all $p \notin \gekl{q,\bar{q}}$ and $\bullet = 0,1$.
\edefin
\setlength{\parindent}{0cm} \setlength{\parskip}{0cm}

Note that we must have $\beta_\mfr^{\bar{q},i} = 0$ and $\beta_\mfs^{q,i} = 0$ for all $i \neq j$.
\setlength{\parindent}{0cm} \setlength{\parskip}{0.5cm}

Given a redundant index $q$ as above, assume first that $\mfr = \mfs$, say $\mfr = \mfs = 0$ (the case $\mfr = \mfs = 1$ is treated analogously). Set $\check{\beta}_\bullet^p \defeq \beta_\bullet^p$ for all $p \neq q, \bar{q}$ and $\bullet = 0,1$, $\check{\beta}_0^{\bar{q}} \defeq \beta_1^{\bar{q}}$, write $\gamma = \beta_0^{\bar{q},j} (\beta_0^{q,j})^{-1}$ and set $\check{\beta}_1^{\bar{q}} \defeq \gamma \beta_1^q$. Set $\check{E} \defeq \bigoplus_{p \in P \setminus \gekl{q}} E^p$ and let $\check{\beta}_\bullet: \: F \to \check{E}$ be given by $\check{\beta}_\bullet = (\check{\beta}_\bullet^p)_{p \in P \setminus \gekl{q}}$ for $\bullet = 0,1$. Let $\check{A} \defeq A(\check{E},F,\check{\beta}_0,\check{\beta}_1)$ and $\check{B} \defeq \menge{f \in \check{A}}{f(t) \in D\check{E} \ \forall \, t \in [0,1]}$. The following is straightforward to check.
\blemma
\label{lem:A=vA_00}
We have an isomorphism $A \isom \check{A}, \, (f^p)_p \ma (\check{f}^p)_p$, where for $f^p \in C([0,1],E^p)$, $\check{f}^p = f^p$ if $p \neq q, \bar{q}$, $\check{f}^{\bar{q}}(t) = f^{\bar{q}}(1-2t)$ for $t \in [0,\half]$ and $\check{f}^{\bar{q}}(t) = \gamma(f^q(2t-1))$ for $t \in (\half,1]$. This isomorphism sends $B$ to $\check{B}$.
\elemma
\setlength{\parindent}{0cm} \setlength{\parskip}{0cm}

If $\mfr \neq \mfs$, say $\mfr = 0$ and $\mfs = 1$ (the other case is analogous), define $\check{\beta}_\bullet^p \defeq \beta_\bullet^p$ for all $p \neq q, \bar{q}$ and $\bullet = 0,1$, $\check{\beta}_0^{\bar{q}} \defeq \beta_0^q$, and $\check{\beta}_1^{\bar{q}} \defeq \gamma \beta_1^{\bar{q}}$, where $\gamma \defeq \beta_1^{q,j} (\beta_0^{\bar{q},j})^{-1}$, set $\check{E} \defeq \bigoplus_{p \in P \setminus \gekl{q}} E^p$, $\check{\beta}_\bullet \defeq (\check{\beta}_\bullet^p)_{p \in P \setminus \gekl{q}}$ for $\bullet = 0,1$, $\check{A} \defeq A(\check{E},F,\check{\beta}_0,\check{\beta}_1)$ and $\check{B} \defeq \menge{f \in \check{A}}{f(t) \in D\check{E} \ \forall \, t \in [0,1]}$. Then the following analogue of Lemma~\ref{lem:A=vA_00} is straightforward:
\blemma
\label{lem:A=vA_01}
We have an isomorphism $A \isom \check{A}, \, (f^p)_p \ma (\check{f}^p)_p$, where for $f^p \in C([0,1],E^p)$, $\check{f}^p \defeq f^p$ if $p \neq q, \bar{q}$, $\check{f}^{\bar{q}}(t) = f^q(2t)$ for $t \in [0,\half]$ and $\check{f}^{\bar{q}}(t) = \gamma(f^{\bar{q}}(2t-1))$ for $t \in (\half,1]$. This isomorphism sends $B$ to $\check{B}$.
\elemma
\setlength{\parindent}{0cm} \setlength{\parskip}{0.5cm}

\bdefin
We say that $A$ is in reduced form if no index in $P$ is redundant.
\edefin
Lemmas~\ref{lem:A=vA_00} and \ref{lem:A=vA_01} allow us to assume that $A$ is in reduced form from now on.

In the following, let us develop direct sum decompositions so that we can reduce our discussion to individual summands, i.e., to the case where $A$ is indecomposable. Let $\sim_P$ be the equivalence relation on $P$ generated by $q \sim_P \bar{q}$ if there are $i \in I$, $\mfr, \mfs \in \gekl{0,1}$ such that $\beta_\mfr^{q,i} \neq 0$ and $\beta_\mfs^{\bar{q},i} \neq 0$. Let $P = \coprod_{l \in L} P_l$ be the decomposition of $P$ into equivalence classes with respect to $\sim_P$. For each $l \in L$, let $E_l \defeq \bigoplus_{p \in P_l} E^p$, $I_l \defeq \{ i \in I: \: \beta_\bullet^{p,i} \neq 0 \ \text{for some} \ \bullet = 0,1 \ \text{and} \ p \in P_l \}$ and $F_l \defeq \bigoplus_{i \in I_l} F^i$. Define $\beta_{\bullet; l} \defeq (\beta_\bullet^{p,i})_{p \in P_l, \, i \in I_l}: \: \bigoplus_{i \in I_l} F^i \to \bigoplus_{p \in P_l} E^p$ for $\bullet = 0, 1$. Set $A_l \defeq A(E_l,F_l,\beta_{0; l},\beta_{1; l})$. The following is straightforward.
\blemma
\label{lem:dirsum}
We have $A = \bigoplus_{l \in L} A_l$, and for each $l \in L$, $A_l$ cannot be further decomposed into (non-trivial) direct summands. Moreover, the decomposition $A = \bigoplus_{l \in L} A_l$ is the unique direct sum decomposition of $A$ into indecomposable direct summands.
\elemma

\bremark
\label{rem:dirsumAB}
The direct sum decomposition in Lemma~\ref{lem:dirsum} is compatible with C*-diagonals in the sense that if $B = \menge{f \in A}{f(t) \in DE \ \forall \, t \in [0,1]}$, then under the direct sum decomposition $A = \bigoplus_{l \in L} A_l$ from Lemma~\ref{lem:dirsum}, we have $B = \bigoplus_{l \in L} B_l$, where $B_l = \menge{f \in A_l}{f(t) \in DE_l \ \forall \, t \in [0,1]}$.
\eremark

\bcor
\label{cor:dirsum}
Every isomorphism $A_{\sigma} \isom A_{\tau}$ restricts to isomorphisms $A_{\sigma; l} \isom A_{\tau; \lambda(l)}$ for all $l \in L$, where $A_{\sigma; l}$ and $A_{\tau; \lambda(l)}$ are the direct summands of $A_{\sigma}$ and $A_{\tau}$ provided by Lemma~\ref{lem:dirsum}, and $\lambda: \: L \isom L$ is a permutation of $L$.
\setlength{\parindent}{0.5cm} \setlength{\parskip}{0cm}

If the isomorphism $A_{\sigma} \isom A_{\tau}$ above sends $B_{\sigma}$ onto $B_{\tau}$, then for all $l \in L$, the isomorphism $A_{\sigma; l} \isom A_{\tau; \lambda(l)}$ above must send $B_{\sigma; l}$ onto $B_{\sigma; \lambda(l)}$, where $B_{\sigma; l}$ and $B_{\sigma; \lambda(l)}$ are as in Remark~\ref{rem:dirsumAB}.
\ecor
\setlength{\parindent}{0cm} \setlength{\parskip}{0cm}

Here, we are implicitly using that the equivalence relations $\sim_P$ does not depend on $\sigma$, $\tau$, i.e., they coincide for $A$, $A_{\sigma}$ and $A_{\tau}$. This is because $\sigma$ and $\tau$ decompose as $\sigma = (\sigma^p)$, $\tau = (\tau^p)$ for permutation matrices $\sigma^p$, $\tau^p$ in $E^p$.
\setlength{\parindent}{0cm} \setlength{\parskip}{0.5cm}

Lemma~\ref{lem:dirsum} and Corollary~\ref{cor:dirsum} allow us to reduce our discussion to the case where $A$ is indecomposable. So let us assume that we have $p_1 \sim_P p_2$ for all $p_1, p_2 \in P$.

Let us now describe the centre $Z(A)$ and its spectrum $\Spec Z(A)$. Let $\sim_Z$ be the equivalence relation on $[0,1] \times P$ generated by $(\mfr,q) \sim_Z (\mfs,\bar{q})$ if $\mfr, \mfs \in \gekl{0,1}$ and there exists $i \in I$ with $\beta_\mfr^{q,i} \neq 0$ and $\beta_\mfs^{\bar{q},i} \neq 0$. Note that on $(0,1) \times P$, $\sim_Z$ is trivial. We write $[\cdot]_Z$ for the canonical projection map $[0,1] \times P \onto ([0,1] \times P) / {}_{\sim_Z}$. Let $[0,1] \times_\bullet P \defeq \menge{(t,p) \in [0,1] \times P}{\beta_\mfr^q \ \text{is unital for all} \ (\mfr,q) \in [t,p]_Z \ \text{if} \ t \in \gekl{0,1}}$.
\blemma
\label{lem:Z}
The centre of $A$ is given by
\begin{align}
\label{e:Z}
  Z(A) = \big\{ &(f^p) = (g^p \cdot 1_{E^p}) \in C([0,1],Z(E)) = \bigoplus_p C([0,1],Z(E^p)):\\
  &g^p \in C[0,1], \, g^q(\mfr) = g^{\bar{q}}(\mfs) \ {\rm if} \ (\mfr,q) \sim_Z (\mfs,\bar{q}), \, g^q(\mfr) = 0 \ \text{if} \ (\mfr,q) \notin [0,1] \times_\bullet P \big\}. \nonumber
\end{align}

We have a homeomorphism $([0,1] \times_\bullet P) / {}_{\sim_Z} \isom \Spec Z(A)$ sending $[t,q]$ to the character $Z(A) \to \Cz, \, (f^p) = (g^p \cdot 1_{E^p}) \ma g^q(t)$. Here $[0,1]$ is given the usual topology and $P$ the discrete topology.
\elemma
\setlength{\parindent}{0cm} \setlength{\parskip}{0cm}

\bproof
If $f = (f^p)$ lies in $Z(A)$, then $f^p$ lies in $Z(C([0,1],E^p)) = C([0,1],Z(E^p))$ for all $p$, hence $f^p = g^p \cdot 1_{E^p}$ for some $g^p \in C[0,1]$. Moreover, if $a \in F$ satisfies $(f(0),f(1)) = \beta_0(a),\beta_1(a))$, then $a \in Z(F)$, i.e., $a = (\alpha^i \cdot 1_{F^i})$ with $\alpha^i \in \Cz$. Now $g^q(\mfr) \cdot 1_{E^q} = f^q(\mfr) = \beta_\mfr^q(a)$ and $\beta_\mfr^{q,i}(\alpha^i \cdot 1_{F^i}) = \alpha^i \beta_\mfr^{q,i}(1_{F^i})$ implies that $g^q(\mfr) = \alpha^i$ if $\beta_\mfr^{q,i} \neq 0$. Hence $g^q(\mfr) = \alpha^i = g^{\bar{q}}(\mfs)$ if both $\beta_\mfr^{q,i} \neq 0$ and $\beta_\mfs^{\bar{q},i} \neq 0$. In addition, we see that $g^q(\mfr) = 0$ and $(\alpha^i) = 0$ if $\beta_\mfr^q$ is not unital. This shows \an{$\subseteq$} in \eqref{e:Z}. For \an{$\supseteq$}, let $f = (g^p \cdot 1_{E^p})$ satisfy $g^p \in C[0,1]$, $g^q(\mfr) = g^{\bar{q}}(\mfs)$ if $(\mfr,q) \sim_Z (\mfs,\bar{q})$ and $g^q(\mfr) = 0$ if $(\mfr,q) \notin [0,1] \times_\bullet P$. For $i \in I$ take any $(\mfr,q) \in \gekl{0,1} \times P$ with $\beta_\mfr^{q,i} \neq 0$ and set $\alpha^i \defeq g^q(\mfr)$. This is well-defined by our assumption on $(g^p)$. Let $a \defeq (\alpha^i \cdot 1_{F^i}) \in F$. Then it is straightforward to see that $(f(0),f(1)) = (\beta_0(a),\beta_1(a))$. Hence $f \in A$, and thus $f \in Z(A)$.
\setlength{\parindent}{0.5cm} \setlength{\parskip}{0cm}

The second part describing $\Spec Z(A)$ is an immediate consequence.
\eproof
\setlength{\parindent}{0cm} \setlength{\parskip}{0cm}

In the following, we will always identify $\Spec Z(A)$ with $([0,1] \times_\bullet P) / {}_{\sim_Z}$ using the explicit homeomorphism from Lemma~\ref{lem:Z}.
\setlength{\parindent}{0cm} \setlength{\parskip}{0.5cm}

Let us show that the points in $\partial \defeq \menge{[\mfr,p]_Z \in \Spec Z(A)}{(\mfr,p) \in \gekl{0,1} \times P}$ are special. Suppose that $A$ is in reduced form, i.e., no index in $P$ is redundant. Further assume that $A$ is indecomposable, so that for all $p_1, p_2 \in P$, we have $p_1 \sim_P p_2$. Let $\sigma$ and $\tau$ be permutation matrices in $E$. Let $\phi: \: A_{\sigma} \isom A_{\tau}$ be an isomorphism. We denote its restriction to $Z(A_{\sigma})$ also by $\phi$, and let $\phi_Z^*$ be the induced homeomorphism $\Spec Z(A_{\tau}) \isom \Spec Z(A_{\sigma})$. Let $\partial_{\sigma} \defeq \menge{[\mfr,p]_Z \in \Spec Z(A_{\sigma})}{(\mfr,p) \in \gekl{0,1} \times P}$ and $\partial_{\tau} \defeq \menge{[\mfr,p]_Z \in \Spec Z(A_{\tau})}{(\mfr,p) \in \gekl{0,1} \times P}$.
\blemma
\label{lem:01special}
We have $\phi_Z^*(\partial_{\tau}) = \partial_{\sigma}$ unless $\# P = 1 = \# I$ and $\beta_0$, $\beta_1$ are isomorphisms.
\elemma
\setlength{\parindent}{0cm} \setlength{\parskip}{0cm}

\bproof
Assume that $\phi^*_Z [\mfr,p]_Z = [t,\bar{q}]_Z$ for some $(\mfr,p) \in \gekl{0,1} \times P$, $(t,\bar{q}) \in (0,1) \times P$. Let 
$$
  I_{[\mfr,p]_Z} \defeq \big\{ i \in I: \: \beta_\mfs^{q,i} \neq 0 \ \text{for some} \ (\mfs,q) \sim_Z (\mfr,p) \big\}.
$$
$\phi$ induces the following commutative diagram with exact rows:
$$
  \xymatrix{
  0 \ar[r] & \spkl{\ker ([t,\bar{q}]_Z)} \ar[d]_{\phi} \ar[r] & A_{\sigma} \ar[d]_{\phi} \ar[r] & E^{\bar{q}} \ar[d]^{\cong} \ar[r] & 0\\
  0 \ar[r] & \spkl{\ker ([\mfr,p]_Z)} \ar[r] & A_{\tau} \ar[r] & \bigoplus_{i \in I_{[\mfr,p]_Z}} F^i \ar[r] & 0
 }
$$
Here the map $ A_{\tau} \to \bigoplus_{i \in I_{[\mfr,p]_Z}} F^i$ sends $f \in A_{\tau}$ to the uniquely determined $a \in \bigoplus_{i \in I_{[\mfr,p]_Z}} F^i$ with $f^q(\mfs) = \beta_\mfs^q(a)$ for all $(\mfs,q) \in [\mfr,p]_Z$. $E^{\bar{q}} \cong \bigoplus_{i \in I_{[\mfr,p]_Z}} F^i$ implies that $\# I_{[\mfr,p]_Z} = 1$, say $I_{[\mfr,p]_Z} = \gekl{i}$.
\setlength{\parindent}{0cm} \setlength{\parskip}{0.5cm}

Moreover, for every sufficiently small open neighbourhood $U$ of $[t,\bar{q}]_Z$ in $\Spec Z(A_{\sigma})$, $U \setminus \gekl{[t,\bar{q}]_Z}$ is homeomorphic to $(0,1) \amalg (0,1)$, while for every sufficiently small neighbourhood $V$ of $[\mfr,p]_Z$ in $\Spec Z(A_{\tau})$, $V \setminus \gekl{[\mfr,p]_Z}$ is homeomorphic to $\coprod_{(\mfs,q) \in [\mfr,p]_Z} (0,1)$. Hence we must have $\# [\mfr,p]_Z = 2$.

Furthermore, if $U$ and $V$ are as above, then for all $\bm{u} \in U$, $A_{\sigma} / \spkl{\ker(\bm{u})}$ has the same dimension as $A_{\sigma} / \spkl{\ker ([t,\bar{q}]_Z)}$, whereas $A_{\tau} / \spkl{\ker(\bm{v})}$ has the same dimension as $A_{\tau} / \spkl{\ker ([\mfr,p]_Z)}$ for all $\bm{v} \in V$ only if $\beta_\mfs^{q,i}$ is an isomorphism $F^i \isom E^q$ for all $(\mfs,q) \in [\mfr,p]_Z$. Now if there exists $(\mfs,q) \in [\mfr,p]_Z$ with $q \neq p$, then $q$ (and equivalently $p$) would be a redundant index in $P$, which is impossible because $A$ is in reduced form. Hence we must have $[\mfr,p]_Z = \gekl{(0,p), (1,p)}$. But this implies that $\gekl{p}$ is an equivalence class with respect to $\sim_P$. Since $A$ is indecomposable, we must have $P = \gekl{p}$, and thus $I = I_{[\mfr,p]_Z}$. Thus, indeed, $\#P = 1 = \# I$, and $\beta_0$, $\beta_1$ are isomorphisms.
\eproof
\setlength{\parindent}{0cm} \setlength{\parskip}{0.5cm}

It is straightforward to deal with the remaining case where $\#P = 1 = \# I$ and $\beta_0$, $\beta_1$ are isomorphisms:
\blemma
\label{lem:PI=11}
If $\#P = 1 = \# I$ and $\beta_0$, $\beta_1$ are isomorphisms, then for all $\dot{t} \in (0,1)$, 
$A_{\tau} \to A_{\tau}, \, f \ma \ti{f}$, with $\ti{f}(t) \defeq \beta_0 \beta_1^{-1} f(t + (1-\dot{t}))$ for $t \in [0,\dot{t}]$ and $\ti{f}(t) \defeq f(t - \dot{t})$ for $t \in [\dot{t},1]$, is an isomorphism sending $B_{\tau}$ onto $B_{\tau}$ such that the induced map $\Spec Z(A_{\tau}) \isom \Spec Z(A_{\tau})$ sends $[\dot{t}]_Z$ to $[0]_Z = [1]_Z$.
\elemma
\setlength{\parindent}{0cm} \setlength{\parskip}{0cm}

Here we identify $[0,1] \times P$ with $[0,1]$, so that there is no need to carry around the $P$-coordinate.
\setlength{\parindent}{0cm} \setlength{\parskip}{0.5cm}

\bcor
\label{cor:01special}
If $(A_{\sigma},B_{\sigma}) \cong (A_{\tau},B_{\tau})$, then there exists an isomorphism $A_{\sigma} \isom A_{\tau}$ sending $B_{\sigma}$ onto $B_{\tau}$ such that the induced map $\Spec Z(A_{\tau}) \isom \Spec Z(A_{\sigma})$ sends $\partial_{\tau}$ to $\partial_{\sigma}$.
\ecor

Let $A = A(E,F,\beta_0,\beta_1)$ and $B = \menge{f \in A}{f(t) \in DE \ \forall \, t \in [0,1]}$. To describe $\Spec B$, let $\cY \defeq \Spec DE$, $\cY^p \defeq \Spec DE^p$, $\cX \defeq \Spec DF$, $\cX^i = \Spec DF^i$, and for $\mfr = 0,1$, let $\cY_\mfr = \Spec(DE \cdot \beta_\mfr(1_F) \cdot DE) = \menge{y \in \cY}{y(\beta_\mfr(1_F)) = 1}$. Let $\bm{b}_\mfr$ be the map $\cY_\mfr \to \cX$ dual to $\beta_\mfr \vert_{DF}: \: DF \to DE$, i.e., $\bm{b}_\mfr (y) = y \circ \beta_\mfr$. We have $\cY = \coprod_p \cY^p$, $\cX = \coprod_i \cX^i$, and with $\cY_\mfr^p \defeq \cY^p \cap \cY_\mfr$, the restriction $\bm{b}_\mfr^p \defeq \bm{b}_\mfr \vert_{\cY_\mfr^p}$ is dual to $\beta_\mfr^p \vert_{DF}: \: DF \to DE^p$. 

Define an equivalence relation $\sim_B$ on $[0,1] \times \cY$ by setting $(\mfr,y) \sim_B (\mfs,\bar{y})$ if $\mfr, \mfs \in \gekl{0,1}$, $y \in \cY_\mfr$, $\bar{y} \in \cY_\mfs$ and $\bm{b}_\mfr(y) = \bm{b}_\mfs(\bar{y})$. Note that on $(0,1) \times \cY$, $\sim_B$ is trivial. We write $[\cdot]_B$ for the canonical projection map $[0,1] \times \cY \onto ([0,1] \times \cY) / {}_{\sim_B}$. Set $[0,1] \times_\bullet \cY \defeq \menge{(t,y) \in [0,1] \times \cY}{y \in \cY_t \ \text{if} \ t \in \gekl{0,1}}$. Let $\bar{\Pi}: \: ([0,1] \times \cY) / {}_{\sim_B} \onto ([0,1] \times P) / {}_{\sim_Z}, \, [t,y] \ma [t,p]$ for $y \in \cY^p$ be the canonical projection. The following is straightforward:
\blemma
\label{lem:SpecB}
We have a homeomorphism $([0,1] \times_\bullet \cY) / {}_{\sim_B} \isom \Spec B$ sending $[t,y]$ to the character 
$$
  B \to \Cz, \, f \ma 
\bfa
y(f(t)) & \rm{if} \ t \in (0,1);\\
\bm{b}_t(y) \big( (\beta_0,\beta_1)^{-1}(f(0),f(1)) \big) & \rm{if} \ t \in \gekl{0,1}.
\efa
$$
Here $[0,1]$ is given the usual topology and $\cY$ the discrete topology.

Moreover, with respect to this description of $\Spec B$ and the description of $\Spec Z(A)$ from Lemma~\ref{lem:Z}, the map $\Pi: \: \Spec B \to \Spec Z(A)$ induced by the canonical inclusion $Z(A) \into B$ is given by the restriction of $\bar{\Pi}$ to $\dom \Pi \defeq \bar{\Pi}^{-1} (\Spec Z(A))$.
\elemma

We are now ready for our main classification theorem. Let $A = A(E,F,\beta_0,\beta_1)$ be in reduced form. Let $\sigma = (\sigma^p)$ and $\tau = (\tau^p)$ be permutation matrices in $E$. Write $\tensor*[_\sigma]{\beta}{_1} \defeq \Ad(\sigma) \circ \beta_1$ and $\tensor*[_{\tau}]{\beta}{_1} \defeq \Ad(\tau) \circ \beta_1$, and let $\tensor*[_\sigma]{\bm{b}}{_1^p}: \: \tensor*[_\sigma]{\cY}{_1^p} \to \cX$, $\tensor*[_\tau]{\bm{b}}{_1^p}: \: \tensor*[_\tau]{\cY}{_1^p} \to \cX$ be the maps dual to $\tensor*[_\sigma]{\beta}{_1^p}, \, \tensor*[_\tau]{\beta}{_1^p}: \: DF \to DE^p$.
\btheo
\label{thm:ClassAB}
We have $(A_\sigma,B_\sigma) \cong (A_\tau,B_\tau)$ if and only if there exist
\setlength{\parindent}{0cm} \setlength{\parskip}{0cm}

\begin{itemize}
\item a permutation $\rho$ of $P$ and for each $p \in P$ a bijection $\Theta^p: \: \cY^p \isom \cY^{\rho(p)}$,
\item a permutation $\kappa$ of $I$ and for each $i \in I$ a bijection $\Xi^i: \: \cX^i \isom \cX^{\kappa(i)}$ giving rise to the bijection $\Xi = \coprod_i \Xi^i: \: \cX = \coprod_i \cX^i \isom \coprod_i \cX^{\kappa(i)} = \cX$, 
\item a map $o: \: P \to \gekl{\pm 1}$
\end{itemize}
such that for every $p \in P$, we have commutative diagrams
\begin{equation}
\label{e:CDYYXX+}
\begin{tikzcd}
  \cY_0^p \arrow[d, "\bm{b}_0^p"'] \arrow[r, "\sim"', "\Theta^p"] & \cY_0^{\rho(p)} \arrow[d, "\bm{b}_0^{\rho(p)}"] 
  & & 
  \tensor*[_\tau]{\cY}{_1^p} 
  \arrow[d, "{\tensor*[_\tau]{\bm{b}}{_1^p}}"']
  \arrow[r, "\sim"', "\Theta^p"] & \tensor*[_\sigma]{\cY}{_1^{\rho(p)}}
  \arrow[d, "{\tensor*[_\sigma]{\bm{b}}{_1^{\rho(p)}}}"]
  \\
  \cX \arrow[r, "\sim"', "\Xi"] & \cX 
  & &
  \cX \arrow[r, "\sim"', "\Xi"] & \cX
\end{tikzcd}
\end{equation}
if $o(p) = +1$,
\begin{equation}
\label{e:CDYYXX-}
\begin{tikzcd}
  \cY_0^p \arrow[d, "\bm{b}_0^p"'] \arrow[r, "\sim"', "\Theta^p"] & \tensor*[_\sigma]{\cY}{_1^{\rho(p)}} \arrow[d, "{\tensor*[_\sigma]{\bm{b}}{_1^{\rho(p)}}}"] 
  & & 
  \tensor*[_\tau]{\cY}{_1^p} 
  \arrow[d, "{\tensor*[_\tau]{\bm{b}}{_1^p}}"']
  \arrow[r, "\sim"', "\Theta^p"] & \cY_0^{\rho(p)}
  \arrow[d, "\bm{b}_0^{\rho(p)}"]
  \\
  \cX \arrow[r, "\sim"', "\Xi"] & \cX 
  & &
  \cX \arrow[r, "\sim"', "\Xi"] & \cX
\end{tikzcd}
\end{equation}
if $o(p) = -1$.
\etheo
\setlength{\parindent}{0cm} \setlength{\parskip}{0cm}

\bproof
\an{$\Larr$}: The commutative diagrams \eqref{e:CDYYXX+} and \eqref{e:CDYYXX-} induce commutative diagrams
\begin{equation*}
\begin{tikzcd}
  \cY_0^p \times \cY_0^p 
  \arrow[r, "\sim"', "\Theta^p \times \Theta^p"] &
  \cY_0^{\rho(p)} \times \cY_0^{\rho(p)} 
  \\
  (\bm{b}_0^p \times \bm{b}_0^p)^{-1}(\coprod_i \cX^i \times \cX^i)
  \arrow[u, hook] 
  \arrow[d, "\bm{b}_0^p \times \bm{b}_0^p"'] 
  \arrow[r, "\sim"', "\Theta^p \times \Theta^p"] & 
  (\bm{b}_0^{\rho(p)} \times \bm{b}_0^{\rho(p)})^{-1}(\coprod_i \cX^{\kappa(i)} \times \cX^{\kappa(i)}) 
  \arrow[u, hook] 
  \arrow[d, "\bm{b}_0^{\rho(p)} \times \bm{b}_0^{\rho(p)}"] 
  \\
  \coprod_i \cX^i \times \cX^i 
  \arrow[r, "\sim"', "\Xi \times \Xi"] &
  \coprod_i \cX^{\kappa(i)} \times \cX^{\kappa(i)}
\end{tikzcd}
\end{equation*}
if $o(p) = +1$,
\begin{equation*}
\begin{tikzcd}
  \cY_0^p \times \cY_0^p 
  \arrow[r, "\sim"', "\Theta^p \times \Theta^p"] &
  \cY_0^{\rho(p)} \times \cY_0^{\rho(p)} 
  \\
  (\bm{b}_0^p \times \bm{b}_0^p)^{-1}(\coprod_i \cX^i \times \cX^i)
  \arrow[u, hook] 
  \arrow[d, "\bm{b}_0^p \times \bm{b}_0^p"'] 
  \arrow[r, "\sim"', "\Theta^p \times \Theta^p"] & 
  (\tensor*[_\sigma]{\bm{b}}{_1^{\rho(p)}} \times \tensor*[_\sigma]{\bm{b}}{_1^{\rho(p)}})^{-1}(\coprod_i \cX^{\kappa(i)} \times \cX^{\kappa(i)}) 
  \arrow[u, hook] 
  \arrow[d, "{\tensor*[_\sigma]{\bm{b}}{_1^{\rho(p)}} \times \tensor*[_\sigma]{\bm{b}}{_1^{\rho(p)}}}"] 
  \\
  \coprod_i \cX^i \times \cX^i 
  \arrow[r, "\sim"', "\Xi \times \Xi"] &
  \coprod_i \cX^{\kappa(i)} \times \cX^{\kappa(i)}
\end{tikzcd}
\end{equation*}
if $o(p) = -1$.
\setlength{\parindent}{0cm} \setlength{\parskip}{0.5cm}

Applying the groupoid C*-algebra construction, and using \cite[Proposition~5.4]{Li18} (see also \cite[Lemmas~3.2 and 3.4]{BL17}), we obtain the commutative diagram
\begin{equation}
\label{e:CDbeta0}
\begin{tikzcd}
  E^p & \arrow[l, "\sim", "\theta^p"'] E^{\rho(p)}
  \\
  F \arrow[u, "\beta_0^p"] & \arrow[l, "\sim", "\xi"'] F \arrow[u]
\end{tikzcd}
\end{equation}
where $\theta^p = (\Theta^p \times \Theta^p)^*$ is the map induced by $\Theta^p \times \Theta^p$, $\xi = (\Xi \times \Xi)^*$ is the map induced by $\Xi \times \Xi$, and the right vertical map is given by $\beta_0^{\rho(p)}$ if $o(p) = +1$ and $\tensor*[_\sigma]{\beta}{_1^{\rho(p)}}$ if $o(p) = -1$.

Similarly, we obtain a commutative diagram
\begin{equation}
\label{e:CDbeta1}
\begin{tikzcd}
  E^p & \arrow[l, "\sim", "\theta^p"'] E^{\rho(p)}
  \\
  F \arrow[u, "{\tensor*[_\tau]{\beta}{_1^p}}"] & \arrow[l, "\sim", "\xi"'] F \arrow[u]
\end{tikzcd}
\end{equation}
where the right vertical map is given by $\tensor*[_\sigma]{\beta}{_1^{\rho(p)}}$ if $o(p) = +1$ and $\beta_0^{\rho(p)}$ if $o(p) = -1$.

Now denote by $\theta$ the isomorphism $E \isom E$ given by $\bigoplus_p \theta^p: \: E = \bigoplus_p E^{\rho(p)} \isom \bigoplus_p E^p = E$. For $f = (f^p) \in C([0,1],E)$, $f^p \in C([0,1],E^p)$, define $\ti{f} \in C([0,1],E$ by $\ti{f} \defeq (\ti{f}^p)$, $\ti{f}^p \defeq f^p$ if $o(p) = +1$ and $\ti{f}^p \defeq f^p \circ (1 - \id)$ if $o(p) = -1$. We claim that $A_{\sigma} \to A_{\tau}, \, (f,a) \ma (\theta(\ti{f}),\xi(a))$ is an isomorphism sending $B_\sigma$ to $B_\tau$. All we have to show is that this map is well-defined, because we can construct an inverse by replacing $\theta$ by $\theta^{-1}$ and $\xi$ by $\xi^{-1}$, and the map clearly sends $B_\sigma$ to $B_\tau$. To see that it is well-defined, we compute
\begin{align*}
  (\theta(\ti{f})(0))^p &= \theta^p (\ti{f}^{\rho(p)}(0)) = \theta^p (f^{\rho(p)}(0)) = \theta^p (\beta_0^{\rho(p)}(a)) \overset{\eqref{e:CDbeta0}}{=} \beta_0^p(\xi(a)) && {\rm if} \ o(p) = +1,\\
  (\theta(\ti{f})(0))^p &= \theta^p (\ti{f}^{\rho(p)}(0)) = \theta^p (f^{\rho(p)}(1)) = \theta^p ( \tensor*[_\sigma]{\beta}{_1^{\rho(p)}} (a)) \overset{\eqref{e:CDbeta0}}{=} \beta_0^p(\xi(a)) && {\rm if} \ o(p) = -1.
\end{align*}
Similarly, $\theta(\ti{f})(1) = \tensor*[_\tau]{\beta}{_1} (\xi(a))$. This shows that $(\ti{f},\xi(a)) \in A_\tau$, as desired.

\an{$\Rarr$}: By Lemmas~\ref{lem:A=vA_00} and \ref{lem:A=vA_01}, we may assume that $A$ is in reduced form, i.e., no index in $P$ is redundant. By Corollary~\ref{cor:dirsum}, we may assume that $A$ is indecomposable, i.e., we have $p_1 \sim_P p_2$ for all $p \in P$. Let $\phi: \: A_\sigma \isom A_\tau$ be an isomorphism with $\phi(B_\sigma) = B_\tau$. Let $\phi^*_B$ be the induced homeomorphism $\Spec B_\tau \isom \Spec B_\sigma$ and $\phi^*_Z$ the induced homeomorphism $\Spec Z(A_\tau) \isom \Spec Z(A_\sigma)$. By Corollary~\ref{cor:01special}, we may assume that $\phi^*_Z(\partial_\tau) = \partial_\sigma$. We have a commutative diagram
\begin{equation*}
\begin{tikzcd}
  \dom \Pi_\tau \arrow[d, "\Pi_\tau"'] \arrow[r, "\sim"', "\phi^*_B"] & \dom \Pi_\sigma \arrow[d, "\Pi_\sigma"]
  \\
  \Spec Z(A_\tau) \arrow[r, "\sim"', "\phi^*_Z"] & \Spec Z(A_\sigma)
\end{tikzcd}
\end{equation*}
where the maps $\Pi_\tau$ and $\Pi_\sigma$ are the ones from Lemma~\ref{lem:SpecB}. $\phi^*_Z$ restricts to a homeomorphism $\Spec Z(A_\tau) \setminus \partial_\tau \isom \Spec Z(A_\sigma) \setminus \partial_\sigma$. As $\Spec Z(A_\tau) \setminus \partial_\tau \cong (0,1) \times P$ and $\Spec Z(A_\sigma) \setminus \partial_\sigma \cong (0,1) \times P$, there must exist a permutation $\rho$ of $P$ and for each $p \in P$ a homeomorphism $\lambda^p$ of $(0,1)$ such that $\phi^*_Z([t,p]_Z) = [\lambda^p(t),\rho(p)]_Z$. Set $o(p) \defeq +1$ if $\lambda^p$ is orientation-preserving and $o(p) \defeq -1$ if $\lambda^p$ is orientation-reversing. For fixed $p$, $\Pi_\tau^{-1} ((0,1) \times \gekl{p}) = (0,1) \times \cY^p$ and $\Pi_\sigma^{-1} ((0,1) \times \gekl{\rho(p)}) = (0,1) \times \cY^{\rho(p)}$, so that we obtain the commutative diagram
\begin{equation*}
\begin{tikzcd}
  (0,1) \times \cY^p \arrow[d, "\Pi_\tau"'] \arrow[r, "\sim"', "\phi^*_B"] & (0,1) \times \cY^{\rho(p)} \arrow[d, "\Pi_\sigma"]
  \\
  (0,1) \times \gekl{p} \arrow[r, "\sim"', "\phi^*_Z"] & (0,1) \times \gekl{\rho(p)}
\end{tikzcd}
\end{equation*}
It follows that there exists a bijection $\Theta^p: \: \cY^p \isom \cY^{\rho(p)}$ such that $\phi^*_B([t,y]_B) = [\lambda^p(t),\Theta^p(y)]_B$ for all $y \in \cY^p$. 

Now consider $\partial \bm{B}_\tau \defeq \menge{[\mfr,y]_B}{\mfr \in \gekl{0,1}, \, y \in \cY_{\mfr}}$, and define $\partial \bm{B}_\sigma$ analogously. $\phi^*_B$ restricts to a bijection $\partial \bm{B}_\tau \isom \partial \bm{B}_\sigma$ because $\partial \bm{B}_\tau = \Spec B_\tau \setminus \Pi_\tau^{-1}( \Spec Z(A_\tau) \setminus \partial_\tau )$ and similarly for $\partial \bm{B}_\sigma$. In addition, we have a bijection $\partial \bm{B}_\tau \isom \cX$ sending $[0,y]_B$ to $\bm{b}_0(y)$ and $[1,y]_B$ to $\tensor*[_\tau]{\bm{b}}{_1}(y)$, and an analogous bijection $\partial \bm{B}_\sigma \isom \cX$. Thus we obtain a bijection $\cX \isom \cX$ which fits into the commutative diagram
\begin{equation}
\label{e:CDBdBdXX}
\begin{tikzcd}
  \partial \bm{B}_\tau \arrow[d, "\cong"'] \arrow[r, "\sim"', "\phi^*_B"] & \partial \bm{B}_\sigma \arrow[d, "\cong"]
  \\
  \cX \arrow[r, "\sim"] & \cX
\end{tikzcd}
\end{equation}
As this bijection $\cX \isom \cX$ corresponds to an isomorphism $F \isom F$ which fits into the commutative diagram
\begin{equation*}
\begin{tikzcd}
  A_\sigma 
  \arrow[d, "{(\beta_0, \tensor*[_\sigma]{\beta}{_1})^{-1} \circ (\ev_0, \ev_1)}"'] 
\arrow[r, "\sim"', "\phi"] 
  & A_\tau \arrow[d, "{(\beta_0, \tensor*[_\tau]{\beta}{_1})^{-1} \circ (\ev_0, \ev_1)}"]
  \\
  F \arrow[r, "\sim"] & F
\end{tikzcd}
\end{equation*}
there must exist a permutation $\kappa$ of $I$ and bijections $\Xi^i: \: \cX^i \isom \cX^{\kappa(i)}$ such that the bijection $\cX \isom \cX$ in \eqref{e:CDBdBdXX} is given by $\Xi \defeq \coprod_i \Xi^i: \: \coprod_i \cX^i \isom \coprod_i \cX^{\kappa(i)}$. 

Now take $p \in P$ with $o(p) = +1$. For $y \in \cY^p$, $[0,\Theta^p(y)]_B$ is mapped under the right vertical map in \eqref{e:CDBdBdXX} to $\bm{b}_0^{\rho(p)}(\Theta^p(y))$. At the same time $[0,\Theta^p(y)]_B = \lim_{t \, {\scriptscriptstyle \searrow} \, 0} \, [\lambda^p(t), \Theta^p(y)]_B = \lim_{t \, {\scriptscriptstyle \searrow} \, 0} \, \phi^*_B([t,y]_B) = \phi^*_B([0,y]_B)$. By commutativity of \eqref{e:CDBdBdXX}, the right vertical map in \eqref{e:CDBdBdXX} sends $\phi^*_B([0,y]_B)$ to $\Xi(\bm{b}_0^p(y))$. Hence $\Xi \circ \bm{b}_0^p = \bm{b}_0^{\rho(p)} \circ \Theta^p$. Similarly, $\Xi \circ \tensor*[_\tau]{\bm{b}}{_1^p} = \tensor*[_\sigma]{\bm{b}}{_1^{\rho(p)}} \circ \Theta^p$. If $o(p) = -1$, then an analogous argument shows that $\Xi \circ \bm{b}_0^p = \tensor*[_\sigma]{\bm{b}}{_1^{\rho(p)}} \circ \Theta^p$ and $\Xi \circ \tensor*[_\tau]{\bm{b}}{_1^p} = \bm{b}_0^{\rho(p)} \circ \Theta^p$. 
\eproof
\setlength{\parindent}{0cm} \setlength{\parskip}{0.5cm}

In \cite{BR}, the authors identify particular 1-dimensional NCCW complexes $A$ (certain dimension drop algebras) with the property that given any two C*-diagonals $B_1$ and $B_2$ of $A$, we have $(A,B_1) \cong (A,B_2)$ if and only if $\Spec B_1 \cong \Spec B_2$ (i.e., $B_1 \cong B_2$). As the latter is obviously a necessary condition, this can be viewed as a rigidity result. Moving towards a rigidity result in our general setting, let us first prove a weaker statement.
\btheo
\label{thm:AB-BZ}
Suppose that $A$ is a 1-dimensional NCCW complex such that for all $(\mfr,p) \in \gekl{0,1} \times P$, $\# \{i \in I: \: \beta_\mfr^{p,i} \neq 0\} \leq 1$. Given two C*-diagonals $B_1$ and $B_2$ of $A$, we have $(A,B_1) \cong (A,B_2)$ if and only if there exists an isomorphism $B_1 \isom B_2$ sending $Z(A)$ onto $Z(A)$.
\etheo
\setlength{\parindent}{0cm} \setlength{\parskip}{0cm}

\bproof
\an{$\Rarr$} is clear. Let us prove \an{$\Larr$}. By Proposition~\ref{prop:AB=AsBs}, it suffices to show that given permutation matrices $\sigma$ and $\tau$ in $E$, $(B_\sigma,Z(A_\sigma)) \cong (B_\tau,Z(A_\tau))$ implies that $(A_\sigma,B_\sigma) \cong (A_\tau,B_\tau)$. As in the proof of Theorem~\ref{thm:ClassAB}, we may assume that $A$ is in reduced form and that $A$ is indecomposable. Suppose that we have an isomorphism $\phi: \: B_\sigma \isom B_\tau$ with $\phi(Z(A_\sigma)) = Z(A_\tau)$. Let $\phi^*: \: \Spec B_\tau \isom \Spec B_\sigma$ be the homeomorphism induced by $\phi$. Define $\partial \bm{B}_\tau$ and $\partial \bm{B}_\sigma$ as in the proof of Theorem~\ref{thm:ClassAB}. Using Lemma~\ref{lem:PI=11} as for Theorem~\ref{thm:ClassAB}, we can without loss of generality assume that $\phi^*(\partial \bm{B}_\tau) = \partial \bm{B}_\sigma$.
\setlength{\parindent}{0.5cm} \setlength{\parskip}{0cm}

As in the proof of Theorem~\ref{thm:ClassAB} we get a bijection $\Xi: \: \cX \isom \cX$ which fits into a commutative diagram
\begin{equation*}
\begin{tikzcd}
  \partial \bm{B}_\tau \arrow[d, "\cong"'] \arrow[r, "\sim"', "\phi^*_B"] & \partial \bm{B}_\sigma \arrow[d, "\cong"]
  \\
  \cX \arrow[r, "\sim"] & \cX
\end{tikzcd}
\end{equation*}
It remains to show that there exists a permutation $\kappa$ of $I$ and bijections $\Xi^i: \: \cX^i \isom \cX^{\kappa(i)}$ such that $\Xi = \coprod_i \Xi^i$. This follows from the observation --- which is a consequence of our assumption --- that $[\mfr,y]_B$, $[\mfs,\ti{y}]_B$ in $\partial \bm{B}_\tau$ are mapped to elements in $\cX^i$ for the same index $i \in I$ if and only if we have for all open neighbourhoods $U$ and $V$ of $[\mfr,y]_B$ and $[\mfs,\ti{y}]_B$ that $\Pi_\tau( U \cap \dom \Pi_\tau ) \cap \Pi_\tau( V \cap \dom \Pi_\tau ) \neq \emptyset$.

Now the rest of the proof proceeds in exactly the same way as the proof of Theorem~\ref{thm:ClassAB}.
\eproof
\setlength{\parindent}{0cm} \setlength{\parskip}{0.5cm}

Now let us present a strong rigidity result in our general context.
\btheo
\label{thm:AB-B}
Suppose that $A$ is a 1-dimensional NCCW complex such that 
\setlength{\parindent}{0cm} \setlength{\parskip}{0cm}

\begin{itemize}
\item $\beta_0$ and $\beta_1$ are unital,
\item for all $i \in I$, there exists exactly one $(\mfr,p) \in \gekl{0,1} \times P$ such that $\beta_\mfr^{p,i} \neq 0$,
\item for all $(\mfr,p) \in \gekl{0,1} \times P$, there exists exactly one $i \in I$ such that $\beta_\mfr^{p,i} \neq 0$,
\item for all these triples $(\mfr, p, i)$, we have $m_\mfr(p,i) \neq 2$,
\item the map defined on such triples sending $(\mfr, p, i)$ to $m_\mfr(p,i)$ must be injective.
\end{itemize}
Then given any two C*-diagonals $B_1$ and $B_2$ of $A$, we have $(A,B_1) \cong (A,B_2)$ if and only if $\Spec B_1 \cong \Spec B_2$.
\etheo
\setlength{\parindent}{0cm} \setlength{\parskip}{0cm}

\bproof
\an{$\Rarr$} is clear. To prove \an{$\Larr$}, by Proposition~\ref{prop:AB=AsBs}, it suffices to show that for any two permutation matrices $\sigma$ and $\tau$ in $E$, $B_\sigma \cong B_\tau$ implies that $(A_\sigma,B_\sigma) \cong (A_\tau,B_\tau)$. Let $\phi: \: B_\sigma \isom B_\tau$ be an isomorphism and $\phi^*: \: \Spec B_\tau \isom \Spec B_\sigma$ its dual map. Our condition $m_\mfr(p,i) = 1$ or $m_\mfr(p,i) \geq 3$ implies that $\phi^*(\partial \bm{B}_\tau) = \partial \bm{B}_\sigma$, where $\partial \bm{B}_\tau = \menge{[\mfr,y]_B}{(\mfr,y) \in \gekl{0,1} \times \cY}$ and similarly for $\partial \bm{B}_\sigma$. Thus $\phi^*$ induces a permutation $\Xi: \: \cX \isom \cX$, which must be of the desired form as in the proof of Theorem~\ref{thm:ClassAB} by our assumptions on $\beta_\mfr^{p,i}$. Similarly, $\phi^*$ induces homeomorphisms $(0,1) \times \cY = \Spec B_\tau \setminus \partial \bm{B}_\tau \isom \Spec B_\sigma \setminus \partial \bm{B}_\sigma = (0,1) \times \cY$ as in the proof of Theorem~\ref{thm:ClassAB}, again using our assumptions on $\beta_\mfr^{p,i}$. Now proceed in exactly the same way as the proof of Theorem~\ref{thm:ClassAB}.
\eproof
\setlength{\parindent}{0cm} \setlength{\parskip}{0.5cm}

\bremark
\label{rem:appBR}
The hypotheses of Theorem~\ref{thm:AB-B} are for instance satisfied in the case of stabilized dimension drop algebras in the sense of \cite{BR}, i.e., where $P = \gekl{p}$, $I = \gekl{i_0,i_1}$, $E = E^p = M_m \otimes M_n \otimes M_o$, $F^{i_0} = M_m \otimes M_o$, $F^{i_1} = M_n \otimes M_o$, $\beta_0 = \beta_0^{p,i_0}: \: M_m \otimes M_o \to M_m \otimes 1_n \otimes M_o \subseteq E^p$ is given by $\id \otimes 1_n \otimes \id$, $\beta_1 = \beta_1^{p,i_1}: \: M_n \otimes M_o \to 1_m \otimes M_n \otimes M_o \subseteq E^p$ is given by $1_m \otimes \id \otimes \id$, and $m, n \geq 3$, $m \neq n$. 
\setlength{\parindent}{0.5cm} \setlength{\parskip}{0cm}

The conclusion of Theorem~\ref{thm:AB-B} is also shown to be true in \cite{BR} using ad-hoc methods in the case where exactly one of $m$ or $n$ is equal to $2$ or when $(m,n,o) = (2,2,1)$. However, contrary to what is claimed in \cite[Theorem~7.8]{BR}, the conclusion of Theorem~\ref{thm:AB-B} is not true in case $m=n$ and $m, n \geq 3$. The problem is that \cite[Remark~7.7]{BR} is not true in this case, as the following example shows.
\eremark
\setlength{\parindent}{0cm} \setlength{\parskip}{0.5cm}

\bex
\label{ex:appBR}
Let $\nu \geq 6$ be an integer and suppose that $\nu$ is not prime, so that $\nu$ has a divisor $\delta \in \gekl{3, \dotsc, \nu - 3}$. Let $M$ be a $\nu \times \nu$-matrix with two identical rows and pairwise distinct columns such that each row and each column has exactly $\delta$ ones, and zeros everywhere else. Such a matrix has been constructed for example in \cite{NJLSB}. Now consider the matrices
$$
M_\sigma \defeq
\rukl{
\begin{smallmatrix}
\frac{2 \nu}{\delta} \cdot M & 0 \\
0 & \frac{2 \nu}{\delta} \cdot M
\end{smallmatrix}
},
\qquad
M_\tau \defeq
\rukl{
\begin{smallmatrix}
\frac{2 \nu}{\delta} \cdot M & 0 \\
0 & (\frac{2 \nu}{\delta} \cdot M)^t
\end{smallmatrix}
}
$$
Then each row and each column of $M_\sigma$ and $M_\tau$ has sum equal to $2\nu$. The columns of $M_\sigma$ are pairwise distinct, whereas $M_\tau$ has two identical rows and two identical columns. Hence we cannot find permutation matrices $P$ and $Q$ such that $M_\sigma = P M_\tau Q$ or 
$M_\sigma = P M_\tau^t Q$. Thus $M_\sigma$ and $M_\tau$ are not congruent in the language of \cite{BR}. However, it is straightforward to see that the bipartite graphs $\Gamma_\sigma$ and $\Gamma_\tau$ attached to $M_\sigma$ and $M_\tau$ are isomorphic (though not in a way which either consistently preserves orientation or consistently reverses orientation), where the bipartite graphs $\Gamma_\bullet = (V, E_\bullet, E_\bullet \to V \times V)$, for $\bullet = \sigma,\tau$, are defined as follows: Let $V \defeq \gekl{1, \dotsc, 2\nu} \times \gekl{0,1}$, $E_\bullet \defeq \menge{(v_0,v_1,\mu)}{v_0, v_1 \in \gekl{1, \dotsc, 2\nu}, \, \mu \in \gekl{1, \dotsc, (M_\bullet)_{v_0,v_1}}}$ for $\bullet = \sigma, \tau$, and define $E_\bullet \to V \times V, \, (v_0,v_1,\mu) \ma ((v_0,0), (v_1,1))$. This shows that \cite[Remark~7.7]{BR} is not true.
\setlength{\parindent}{0.5cm} \setlength{\parskip}{0cm}

This leads to an example of a 1-dimensional NCCW complex in the same form as in Remark~\ref{rem:appBR} with $m = n = 2\nu$, $o = 1$, and permutation matrices $\sigma, \tau \in E^p$ such that $B_\sigma \cong B_\tau$ but $(A_\sigma, B_\sigma) \not\cong (A_\tau, B_\tau)$: For $\cY = \Spec DE$, $\cX = \Spec DF$, we have canonical identifications $\cY \cong \gekl{1, \dotsc, 2\nu} \times \gekl{1, \dotsc, 2\nu}$, $\cX = \cX^{i_0} \amalg \cX^{i_1}$ with $\cX^{i_0} \cong \gekl{1, \dotsc, 2\nu}$ and $\cX^{i_1} \cong \gekl{1, \dotsc, 2\nu}$ such that the maps $\bm{b}_\bullet$ dual to $\beta_\bullet$ are given by $\bm{b}_0: \: \cY \to \cX^{i_0} \subseteq \cX, \, (y_0,y_1) \ma y_0$ and $\bm{b}_1: \: \cY \to \cX^{i_1} \subseteq \cX, \, (y_0,y_1) \ma y_1$. Let $\bm{\sigma}$ be the permutation of $\cY$ such that for all $x_0, x_1 \in \cX$, we have $\# \menge{y \in \cY}{(\bm{b}_0(y), \bm{b}_1(\bm{\sigma}(y))) = (x_0,x_1)} = (M_{\sigma})_{x_0,x_1}$, and let $\bm{\tau}$ be a permutation of $\cY$ with the analogous property for $M_\tau$ instead of $M_\sigma$. The proof of \cite[Proposition~6.9]{BR} gives a precise recipe to find such $\bm{\sigma}$ and $\bm{\tau}$. Now let $\sigma$ be the permutation matrix in $E$ given by $\sigma_{\bar{y},y} = 1$ if and only if $\bar{y} = \bm{\sigma}^{-1}(y)$, and define $\tau$ similarly. Then $\Gamma_\sigma \cong \Gamma_\tau$ implies that $\Spec B_\sigma \cong \Spec B_\tau$, hence $B_\sigma \cong B_\tau$, while we cannot have $(A_\sigma, B_\sigma) \cong (A_\tau, B_\tau)$ since otherwise Theorem~\ref{thm:ClassAB} would imply that $M_\sigma$ and $M_\tau$ would have to be congruent in the sense of \cite{BR}.
\eex
\setlength{\parindent}{0cm} \setlength{\parskip}{0.5cm}

\section{Construction of C*-diagonals with connected spectra}

We set out to construct C*-diagonals with connected spectra in classifiable stably finite C*-algebras.

\subsection{Construction of C*-diagonals in classifiable stably finite C*-algebras}
\label{ss:CCCC}

We recall the construction in \cite[\S~4]{Li18} which is a modified version of the constructions in \cite{Ell, EV, GLN} (see \cite[\S~2]{Li18}). The construction provides a model for every classifiable stably finite C*-algebra which is unital or stably projectionless with continuous scale, with prescribed Elliott invariant $\cE = (G_0, G_0^+, u, T, r, G_1)$ as in \cite[Theorem~1.2]{Li18} or $\cE = (G_0, \gekl{0}, T, \rho, G_1)$ as in \cite[Theorem~1.3]{Li18}, in the form of an inductive limit $\ilim_n \gekl{A_n, \varphi_n}$. In addition, the crucial point in \cite{Li18} is to identify C*-diagonals $B_n$ of $A_n$ which are preserved under the connecting maps and which satisfy the hypothesis of \cite[Theorem~1.10]{Li18} so that $\ilim_n \gekl{B_n, \varphi_n}$ becomes a C*-diagonal of $\ilim_n \gekl{A_n, \varphi_n}$. Here $A_n = \menge{(f,a) \in C([0,1],E_n) \oplus F_n}{f(\mfr) = \beta_{n, \mfr}(a) \ {\rm for} \ \mfr = 0,1}$ where $E_n = \bigoplus_p E_n^p$, $E_n^p = M_{\gekl{n,p}}$, $F_n = \bigoplus_i F_n^i$, $F_n^{\bm{i}} = P_n^{\bm{i}} M_{\infty}(C(Z_n)) P_n^{\bm{i}}$ for a distinguished index $\bm{i}$, $Z_n$ is a path-connected space with base point $\theta^{\bm{i}}_n$, $P^{\bm{i}}_n$ is a projection corresponding to a vector bundle over $Z_n$ of dimension $[n,\bm{i}]$, $P^{\bm{i}}_n(\theta^{\bm{i}}_n)$ is up to conjugation by a permutation matrix given by $1_{[n,\bm{i}]}$, $F_n^i = M_{[n,i]}$ for $i \neq \bm{i}$, $\hat{F}_n = \bigoplus \hat{F}_n^i$, $\hat{F}_n^{\bm{i}} = M_{[n,\bm{i}]}$, $\hat{F}_n^i = F_n^i$ if $i \neq \bm{i}$, $\pi_n: \: F_n \onto \hat{F}_n$ is given by $\pi_n = \ev_{\theta_n^{\bm{i}}} \oplus \bigoplus_{i \neq \bm{i}} \id_{F_n^i}$, and $\beta_{n, \bullet} = \hat{\beta}_{n, \bullet} \circ \pi_n$, where $\hat{\beta}_{n, \bullet}: \: \hat{F}_n \to E_n$ is of the same form as in \eqref{e:beta=}. In the stably projectionless case, we can (and will) always arrange that for all $n$, there exists exactly one index $\grave{p}$ such that $\beta_{n,0}^{\grave{p}}$ is unital and $\beta_{n,1}^{\grave{p}}$ is non-unital, while $\beta_{n,\bullet}^p$ is unital for all other $p \neq \grave{p}$. 

The connecting maps $\varphi \defeq \varphi_n: \: A_n \to A_{n+1}$ are determined by $\varphi_C: \: A_n \overset{\varphi}{\lori} A_{n+1} \to C([0,1],E_{n+1})$ and $\varphi_F: \: A_n \overset{\varphi}{\lori} A_{n+1} \to F_{n+1}$. $\varphi_C(f,a)$ is of block diagonal form, with block diagonal entries of the form 
\begin{equation}
\label{e:fplp}
f^p \circ \lambda,
\end{equation}
for a continuous map $\lambda: \: [0,1] \to [0,1]$ with $\lambda^{-1}(\gekl{0,1}) \subseteq \gekl{0,1}$, where $f^p$ is the image of $f$ under the canonical projection $C([0,1],E_n) \onto C([0,1],E_n^p)$ (see \cite[Equation~(16)]{Li18}), or of the form
\begin{equation}
\label{e:tauax}
[0,1] \to E_{n+1}^q, \, t \ma \tau(t) a(x(t)), 
\end{equation}
where $x: \: [0,1] \to Z_n$ is continuous and $\tau(t): \: P_n(x(t)) M_{\infty} P_n(x(t)) \cong P_n(\theta_n^i) M_{\infty}P_n(\theta_n^i)$ is an isomorphism depending continuously on $t$, with $\theta_n^i$ in the same connected component as $x(t)$, and $\tau(t) = \id$ if $x(t) = \theta_n^i$ (see \cite[Equation~(17)]{Li18}). Moreover, we can always arrange the following conditions:
\begin{eqnarray}
\label{e:phiCfp}
&&\forall \ p, \, q \ \exists \ \text{a block diagonal entry in} \ C([0,1],E_{n+1}^q) \ \text{of the form} \ f^p \circ \lambda^p \ \text{as in} \ \eqref{e:fplp} \ \text{with} \ \lambda^p(0) = 0, \, \lambda^p(1) = 1. \\
\label{e:phiCl}
&&\forall \ \lambda \ \text{as in} \ \eqref{e:fplp} \ \text{and} \ \mfr \in \gekl{0,1}, \, \lambda(\mfr) \notin \gekl{0,1} \ \Rarr \ \lambda(\mfr^*) \in \gekl{0,1}, \ \text{where} \ \mfr^* = 1 - \mfr.\\
\label{e:phiCx}
&&\forall \ x \ \text{as in} \ \eqref{e:tauax}, \ \text{if} \ \img(x) \subseteq Z_n^{\bm{i}}, \ \text{then} \ x(0) = \theta_n^{\bm{i}} \ \text{or} \ x(1) = \theta_n^{\bm{i}}.
\end{eqnarray}
Note that a crucial (though basic) modification of the constructions in \cite{Ell, EV, GLN} is to push unitary conjugation all into $\beta_{n+1,\bullet}$, so that $\varphi_C$ can be arranged to be always of this block diagonal form (see \cite[Remark~4.1]{Li18} for details).

$\varphi_F(f,a)$ is up to permutation given by
\begin{equation}
\label{e:phiF}
\varphi_F(f,a) = 
\rukl{
\begin{smallmatrix}
\varphi_{F,C}(f) & 0 \\
0 & \varphi_{F,F}(a)
\end{smallmatrix}
},
\quad
\text{where} \ 
\varphi_{F,F}(a) = (\varphi^{j,i}(a^i))_j
\ \text{for} \ F_n = \bigoplus_i F_n^i, \, F_{n+1} = \bigoplus_j F_{n+1}^j.
\end{equation}
Moreover, with $\pi \defeq \pi_{n+1}$, $\pi \circ \varphi^{j,i}$ is given by the composition
\begin{equation}
\label{e:piphiji}
F_n^i = \hat{F}_n^i \overset{1 \otimes \id_{\hat{F}_n^i}}{\lori} 1_{m(j,i)} \otimes \hat{F}_n^i \subseteq M_{m(j,i)} \otimes \hat{F}_n^i \tailarr \hat{F}_{n+1}^j \quad \text{if} \ i \neq \bm{i},
\end{equation}
$\pi \circ \varphi^{j,\bm{i}}$ is given by
\begin{equation}
\label{e:piphijbi}
\rukl{
\begin{smallmatrix}
\pi \circ \varphi_\theta^{j,\bm{i}} & 0 \\
0 & \pi \circ \varphi_{\bm{Z}}^{j,\bm{i}}
\end{smallmatrix}
}
\end{equation}
where $\pi \circ \varphi_\theta^{j,\bm{i}}$ is given by the composition
\begin{equation}
\label{e:piphijbitheta}
F_n^{\bm{i}} \overset{\ev_{\theta_n^{\bm{i}}}}{\lori} \hat{F}_n^{\bm{i}} \overset{1 \otimes \id_{\hat{F}_n^{\bm{i}}}}{\lori} 1_{m(j,\bm{i})} \otimes F_n^{\bm{i}} \subseteq M_{m(j,\bm{i})} \otimes F_n^{\bm{i}} \tailarr \hat{F}_{n+1}^j,
\end{equation}
and $\pi \circ \varphi_{\bm{Z}}^{j,\bm{i}}$ consists of block diagonals of a similar form as $\pi \circ \varphi_\theta^{j,\bm{i}}$, but starting with $\ev_{\bm{z}}$ instead of $\ev_{\theta_n^{\bm{i}}}$, for $\bm{z} \in \bm{Z} \subseteq Z_n \setminus \gekl{\theta_n^{\bm{i}}}$. As in \eqref{e:beta=}, the arrow $\tailarr$ denotes a *-homomorphism of multiplicity $1$ sending diagonal matrices to diagonal matrices. It is convenient to collect $\pi \circ \varphi^{j,i}$ for all $j$ into a single map $\pi \circ \varphi^{-,i}: \: F_n^i \to \hat{F}_{n+1}$ given by
\begin{equation}
\label{e:piphi-i}
F_n^i = \hat{F}_n^i \to (1_{m(j,i)} \otimes \hat{F}_n^i)_j \subseteq \Big( \bigoplus_j M_{m(j,i)} \Big) \otimes \hat{F}_n^i \tailarr \hat{F}_{n+1} \quad \text{if} \ i \neq \bm{i},
\end{equation}
and in a similar way, we obtain $\pi \circ \varphi_{\theta}^{-,\bm{i}}: \: F_n^{\bm{i}} \to \hat{F}_{n+1}$ given by
\begin{equation}
\label{e:piphi-bi}
F_n^{\bm{i}} \overset{\ev_{\theta_n^{\bm{i}}}}{\lori} \hat{F}_n^{\bm{i}} \to (1_{m(j,{\bm{i}})} \otimes \hat{F}_n^{\bm{i}})_j \subseteq \Big( \bigoplus_j M_{m(j,{\bm{i}})} \Big) \otimes \hat{F}_n^{\bm{i}} \tailarr \hat{F}_{n+1}.
\end{equation}

\subsection{Modification (conn)}
\label{ss:(conn)}

We modify the construction described in \S~\ref{ss:CCCC} to obtain C*-diagonals with connected spectra. We start with the inductive limit decomposition as in \cite[\S~2]{Li18} and construct C*-algebras $F_n$ as in \cite[\S~3,4]{Li18}. Now the original construction recalled in \S~\ref{ss:CCCC} produces $\bfdot{A}_n$ and $\bfdot{\varphi}_{n-1}$ inductively on $n$. Suppose that the original construction starts with the C*-algebra $\bfdot{A}_1$ of the form as in \S~\ref{ss:CCCC}. Let us explain how to modify it. We will use the same notation as in \S~\ref{s:nccw}. Let $[1,I] \defeq \sum_i [1,i]$. Choose an index $\mfp$ and define $E_1^{\mfp} \defeq M_{\gekl{1,\mfp} + [1,I]}$. View $\bfdot{E}_1^{\mfp}$ and $\hat{F}_1$ as embedded into $E_1^{\mfp}$ via $\bfdot{E}_1^{\mfp} \oplus \hat{F}_1 = M_{\gekl{1,\mfp}} \oplus (\bigoplus_i M_{[1,i]}) \subseteq M_{\gekl{1,\mfp}} \oplus M_{[1,I]} \subseteq E_1^{\mfp}$. Define $E_1^p \defeq \bfdot{E}_1^p$ for all $p \neq \mfp$ and $E_1 \defeq \bigoplus_p E_1^p$. Let $d_l$, $1 \leq l \leq [1,I]$, be the rank-one projections in $DM_{[1,I]}$ and $w \in M_{[1,I]}$ the permutation matrix such that $w d_l w^* = d_{l+1}$ if $1 \leq l \leq [1,I] - 1$ and $w d_{[1,I]} w^* = d_1$. Define $\beta_{1,\mfr}^p \defeq \bfdot{\beta}_{1,\mfr}^p$ for $\mfr = 0,1$, $p \neq \mfp$, $\beta_{1,0}^{\mfp} \defeq (\bfdot{\beta}_{1,0}^{\mfp}, \pi)$ and $\beta_{1,1}^{\mfp} \defeq \Ad(1_{\bfdot{E}_1^{\mfp}}, w) \circ (\beta_{1,1}^{\mfp}, \pi_1)$ as maps $F_1 \to \bfdot{E}_1^{\mfp} \oplus \hat{F}_1 \subseteq E_1^{\mfp}$. Now define $A_1 \defeq \menge{(f,a) \in C([0,1],E_1) \oplus F_1}{f(\mfr) = \beta_{1,\mfr}(a) \ {\rm for} \ \mfr = 0,1}$. 

Now suppose that our new construction produced
$$
A_1 \overset{\varphi_1}{\lori} A_2 \overset{\varphi_2}{\lori} \dotso \overset{\varphi_{n-1}}{\lori} A_n.
$$
Let $\bfdot{A}_{n+1}$ and $\bfdot{\varphi}_n: \: A_n \to \bfdot{A}_{n+1}$ be given by the original construction as recalled in \S~\ref{ss:CCCC}. In order to modify $\bfdot{A}_{n+1}$ and $\bfdot{\varphi}_n$, we use the same notation for $\bfdot{\varphi} \defeq \bfdot{\varphi}_n$ as in \S~\ref{ss:CCCC}. Let $[n+1,J] \defeq \sum_j [n+1,j]$. Choose an index $\mfq$. Define $E_{n+1}^{\mfq} \defeq M_{\gekl{n+1,\mfq} + [n+1,J]}$. View $\bfdot{E}^{\mfq}_{n+1}$ and $\hat{F}_{n+1}$ as embedded into $E^{\mfq}_{n+1}$ via $\bfdot{E}^{\mfq}_{n+1} \oplus \hat{F}_{n+1} = M_{\gekl{n+1,\mfq}} \oplus (\bigoplus_j M_{[n+1,j]}) \subseteq M_{\gekl{n+1,\mfq}} \oplus M_{[n+1,J]} \subseteq E^{\mfq}_{n+1}$. Define $E^q_{n+1} \defeq \bfdot{E}^q_{n+1}$ for all $q \neq \mfq$ and $E_{n+1} \defeq \bigoplus_q E^q_{n+1}$. Set $\beta_\mfr^q \defeq \bfdot{\beta}_\mfr^q$ for $\mfr = 0,1$, $q \neq \mfq$ and $\beta_0^{\mfq} \defeq (\bfdot{\beta}_0^{\mfq}, \pi)$ as a map $F_{n+1} \to \bfdot{E}^{\mfq}_{n+1} \oplus \hat{F}_{n+1} \subseteq E^{\mfq}_{n+1}$. Let us now define $\beta_1$. Consider the descriptions of $\pi \circ \bfdot{\varphi}^{j,i}$ for $i \neq \bm{i}$ in \eqref{e:piphiji} and \eqref{e:piphi-i} and of $\pi \circ \bfdot{\varphi}^{j,\bm{i}}_\theta$ in \eqref{e:piphijbitheta} and \eqref{e:piphi-bi}. Let $d_l^i$, $1 \leq l \leq \sum_j m(j,i)$, be the rank-one projections in $\bigoplus_j DM_{m(j,i)}$ and $w^i, w^{\bm{i}}_\theta \in M_{[n+1,J]}$ permutation matrices such that, if we identify $d_l^i \otimes \mff$ with its image in $E_{n+1}^\mfq$ under the compositions of the embeddings $\big( \bigoplus_j M_{m(j,i)} \big) \otimes \hat{F}_n^i \tailarr \hat{F}_{n+1}$ from \eqref{e:piphi-i}, \eqref{e:piphi-bi} and $\hat{F}_{n+1} = \bigoplus_j M_{[n+1,j]} \subseteq M_{[n+1,J]}$ from above, we have $w^i (d_l^i \otimes \mff) (w^i)^* = d_{l+1}^i \otimes \mff$ if $1 \leq l \leq \sum_j m(j,i) - 1$, $w^i (d_l^i \otimes \mff) (w^i)^* = d_1^i \otimes \mff$ if $l = \sum_j m(j,i)$, $w^{\bm{i}}_\theta (d_l^{\bm{i}} \otimes \mff) (w^{\bm{i}}_\theta)^* = d_{l+1}^{\bm{i}} \otimes \mff$ if $1 \leq l \leq \sum_j m(j,\bm{i}) - 1$ and $w^{\bm{i}}_\theta (d_l^{\bm{i}} \otimes \mff) (w^{\bm{i}}_\theta)^* = d_1^{\bm{i}} \otimes \mff$ if $l = \sum_j m(j,\bm{i})$, for all $\mff \in D\hat{F}_n^{\bm{i}}$. Let $\mfe^i$ be the unit of $\big( \bigoplus_j M_{m(j,i)} \big) \otimes \hat{F}_n^i$, viewed as a projection in $M_{[n+1,J]}$ via the above embedding into $M_{[n+1,J]}$, and let $\mfe^{\bm{i}}_\theta$ be the unit of $\big( \bigoplus_j M_{m(j,\bm{i})} \big) \otimes \hat{F}_n^{\bm{i}}$, viewed as a projection in $M_{[n+1,J]}$ via the above embedding into $M_{[n+1,J]}$. Let $\mfe_{F,C}$ and $\mfe_{F,F}$ be the projections in $\hat{F}_{n+1}$ corresponding to the decomposition of $\bfdot{\varphi}_F$ (or rather $\pi \circ \bfdot{\varphi}_F$) in \eqref{e:phiF} so that $1_{\hat{F}_{n+1}} = \mfe_{F,C} + \mfe_{F,F}$, and set $\mfe_{\bm{Z}} \defeq \mfe_{F,F} - (\sum_{i \neq \bm{i}} \mfe^i) - \mfe_\theta^{\bm{i}}$. Now define $w_{F,F} \defeq (\sum_{i \neq \bm{i}} \mfe^i w^i \mfe^i) + \mfe_\theta^{\bm{i}} w_\theta^{\bm{i}} \mfe_\theta^{\bm{i}} + \mfe_{\bm{Z}}$ and 
\begin{equation}
\label{e:w=1w}
w \defeq 
\rukl{
\begin{smallmatrix}
\mfe_{F,C} & 0 \\
0 & w_{F,F}
\end{smallmatrix}
}
\end{equation}
with respect to the decomposition of $\bfdot{\varphi}_F$ in \eqref{e:phiF}. Set
$$
\beta_1^{\mfq} \defeq
\Ad
\rukl{
\begin{smallmatrix}
1_{\bfdot{E}_{n+1}^{\mfq}} & 0 \\
0 & w
\end{smallmatrix}
}
\circ
\rukl{
\begin{smallmatrix}
\bfdot{\beta}_1^{\mfq} & 0 \\
0 & \pi
\end{smallmatrix}
}
: \: F_{n+1} \to \bfdot{E}_{n+1}^{\mfq} \oplus \hat{F}_{n+1} \subseteq E_{n+1}^{\mfq}. 
$$
Finally, define $A_{n+1} \defeq \menge{(f,a) \in C([0,1],E_{n+1}) \oplus F_{n+1}}{f(\mfr) = \beta_\mfr(a) \ \text{for} \ \mfr = 0,1}$ and $\varphi = \varphi_n: \: A_n \to A_{n+1}$ by $\varphi_F \defeq \bfdot{\varphi}_F$, $\varphi_C^q \defeq \bfdot{\varphi}_C^q$ for $q \neq \mfq$, and 
\begin{equation}
\label{e:dotphibqC}
\varphi_C^{\mfq} \defeq
\rukl{
\begin{smallmatrix}
\bfdot{\varphi}_C^{\mfq} & 0 \\
0 & \pi \circ \bfdot{\varphi}_F
\end{smallmatrix}
}
: \: A_n \to C([0,1], \bfdot{E}_{n+1}^{\mfq} \oplus \hat{F}_{n+1}) \subseteq C([0,1], E_{n+1}^{\mfq}).
\end{equation}
By construction, $\varphi_n$ is well-defined, i.e., $\varphi_n(f,a)$ satisfies the defining boundary conditions for $A_{n+1}$ for all $(f,a) \in A_n$. Proceeding in this way, we obtain an inductive system $\gekl{A_n, \varphi_n}_n$.

\blemma
\label{lem:conn:Ell}
$A \defeq \ilim_n \gekl{A_n, \varphi_n}$ is a classifiable C*-algebra with ${\rm Ell}(A) \cong \cE$, where $\cE = (G_0, G_0^+, u, T, r, G_1)$ as in \cite[Theorem~1.2]{Li18} or $\cE = (G_0, \gekl{0}, T, \rho, G_1)$ as in \cite[Theorem~1.3]{Li18}. In the latter case, $A$ has continuous scale.
\setlength{\parindent}{0.5cm} \setlength{\parskip}{0cm}

If we set $B_n \defeq \menge{(f,a) \in A_n}{f(t) \in DE_n \ \forall \ t \in [0,1], \, a \in DF_n}$, then $B \defeq \ilim_n \gekl{B_n, \varphi_n}$ is a C*-diagonal of $A$.
\elemma
\setlength{\parindent}{0cm} \setlength{\parskip}{0cm}

Here $DF_n$ is the C*-diagonal of $F_n$ defined in \cite[\S~6.1]{Li18}.

\bproof
$A$ is classifiable and unital or stably projectionless with continuous scale for the same reasons why the original construction recalled in \S~\ref{ss:CCCC} yields C*-algebras with these properties (see \cite{Li18, Ell, EV, GLN} for details). We also have ${\rm Ell}(A) \cong \cE$ for the same reasons as for the original construction. This is straightforward for K-theory, as $A_{n+1}$ and $\bfdot{A}_{n+1}$ have the same K-theory and $\varphi_n$ induces the same map on K-theory as $\bfdot{\varphi}_n$. It is also straightforward to see that modification (conn) yields the desired trace simplex and pairing between $K_0$ and traces. Indeed, we can think of our modification taking place already at the first stage of the construction summarized in \cite[\S~2]{Li18}, where a non-simple C*-algebra with the prescribed Elliott invariant is constructed. And that this non-simple C*-algebra has the desired trace simplex and pairing is enforced in the construction summarized in \cite[\S~2]{Li18} by making sure that for the analogues of $\bfdot{\varphi}_C$ and $\bfdot{\varphi}_F$, the block diagonal entries of the form $t \ma \tau(t) a(x(t))$ as in \eqref{e:tauax} and $\bfdot{\varphi}_{F,F}$ as in \eqref{e:phiF} take up larger and larger portions of $C([0,1],\bfdot{E}_{n+1})$. But our modification only increases these portions.
\setlength{\parindent}{0.5cm} \setlength{\parskip}{0cm}

Finally, the connecting maps $\varphi_n$ are of the same form as in \cite[\S~4]{Li18}, and hence admit groupoid models as in \cite[\S~6]{Li18}. Hence $B$ is indeed a C*-diagonal of $A$ by the same argument as in \cite[\S~5 -- 7]{Li18}.  
\eproof
\setlength{\parindent}{0cm} \setlength{\parskip}{0.5cm}

\subsection{Building block C*-diagonals with path-connected spectra}
\label{ss:BBCDiagPathConn}

Let us now show that modification (conn) yields C*-diagonals with connected spectra. We need the following notations: Let $\cY_n \defeq \Spec DE_n$, $\bfdot{\cY}_n^p \defeq \Spec D\bfdot{E}_n^p$, $\cY_n^{\mfp} \defeq \Spec DE_n^{\mfp}$ so that $\cY_n = \cY_n^{\mfp} \amalg \coprod_{p \neq \mfp} \bfdot{\cY}_n^p$, $\cX_n \defeq \Spec D\hat{F}_n$, $\cX_n \defeq \Spec D\hat{F}_n$, $\cX_n^i \defeq \Spec D\hat{F}_n^i$ and $\cF_n^{(0)} \defeq \Spec DF_n$. We have $\cF_n^{(0)} \cong (Z_n \times \cX^{\bm{i}}_n) \amalg (\coprod_{i \neq \bm{i}} \gekl{\theta_n^i} \times \cX_n^i)$. $\pi: \: F_n \onto \hat{F}_n$ restricts to $DF_n \onto D\hat{F}_n$, which induces $\cX_n \into \cF_n^{(0)}$ given by $\cX_n^i \into \gekl{\theta_n^i} \times \cX_n^i, \, x \ma (\theta_n^i, x)$ with respect to the identification of $\cF_n^{(0)}$ we just mentioned. We identify $\cX_n$ with a subset of $\cF_n^{(0)}$ in this way. Let $\bm{b}_{n,\mfr}^p: \: \cY_{n,\mfr}^p \to \cX_n$, where $\cY_{n,\mfr}^p \defeq \dom \bm{b}_{n,\mfr}^p \subseteq \cY_n^p$, be the map inducing $\beta_{n,\mfr}^p$, define $\bm{b}_{n,\mfr}: \: \cY_{n,\mfr} \to \cX_n$ correspondingly, and let $\bfdot{\bm{b}}_{n,\mfr}^{\mfp}: \: \bfdot{\cY}_{n,\mfr}^{\mfp} \to \cX_n$, with $\bfdot{\cY}_{n,\mfr}^{\mfp} \defeq \dom \bfdot{\bm{b}}_{n,\mfr}^{\mfp} \subseteq \bfdot{\cY}_n^{\mfp}$, be the map inducing $\bfdot{\beta}_{n,\mfr}^{\mfp}$. Let $\sim$ be the equivalence relation on $([0,1] \times \cY_n) \amalg \cF_n^{(0)}$ generated by $(\mfr,y) \sim \bm{b}_{n,\mfr}(y) \in \cX_n \subseteq \cF_n^{(0)}$ for $\mfr \in \gekl{0,1}$ and $y \in \cY_{n,\mfr}$. We write $[\cdot]$ for the canonical projection map $([0,1] \times \cY_n) \amalg \cF_n^{(0)} \onto \big( ([0,1] \times \cY_n) \amalg \cF_n^{(0)} \big) / {}_{\sim}$ and identify $\cF_n^{(0)}$ with its image under $[\cdot]$. Set $[0,1] \times_\bullet \cY_n \defeq \menge{(t,y) \in [0,1] \times \cY_n}{y \in \cY_{n,t} \ \text{if} \ t \in \gekl{0,1}}$. 
The following generalization of Lemma~\ref{lem:SpecB} is straightforward:
\blemma
\label{lem:SpecB_gen}
We have a homeomorphism $\big( ([0,1] \times_\bullet \cY_n) \amalg \cF_n^{(0)} \big) / {}_{\sim} \isom \Spec B_n$ sending $[t,y]$ (for $(t,y) \in [0,1] \times \cY_n$) to the character $B_n \to \Cz, \, (f,a) \ma y(f(t))$ and $\bm{x} \in \cF_n^{(0)}$ to the character $B_n \to \Cz, \, (f,a) \ma \bm{x}(a)$.
\elemma

Moreover, we always have $\cY_n^{\mfp} = \bfdot{\cY}_n^{\mfp} \amalg \cX_n$, $\cY_{n,\mfr}^{\mfp} = \bfdot{\cY}_{n,\mfr}^{\mfp} \amalg \cX_n$ and
\begin{equation}
\label{e:bn0bp}
\bm{b}_{n,0}^{\mfp} = \bfdot{\bm{b}}_{n,0}^{\mfp} \amalg \id_{\cX_n}.
\end{equation}
For $n=1$, $\bm{b}_{1,1}$ is given by $\bm{b}_{1,1}^{\mfp} \vert_{\cY_{1,1}^{\mfp}} = \bfdot{\bm{b}}_{1,1}^{\mfp}$, and if $\gekl{x_l}_{1 \leq l \leq [1,I]} = \cX_1$ according to the enumeration of rank-one projections in $DM_{[1,I]}$ in \S~\ref{ss:(conn)}, we have
\begin{equation}
\label{e:b11bp}
\bm{b}_{1,1}^{\mfp}(x_l) = x_{l-1} \ \text{if} \ 2 \leq l \leq [1,I] \qquad \text{and} \ \bm{b}_{1,1}^{\mfp}(x_1) = x_{[1,I]}.
\end{equation}
Now we need to describe the groupoid model $\bm{p} \defeq \bm{p}_n$ for $\varphi_n$. Let us drop the index $n+1$ and write $\cY \defeq \cY_{n+1}$, $\cX \defeq \cX_{n+1}$ and so on. By construction of $E^{\mfq}$, we have a decomposition $\cY^{\mfq} = \bfdot{\cY}^{\mfq} \amalg \cX$. Moreover, according to the decomposition of $\pi \circ \bfdot{\varphi}_F$ in \S~\ref{ss:CCCC} (see \eqref{e:phiF} -- \eqref{e:piphijbitheta} in combination with the definition of $\mfe^i$, $\mfe_\theta^{\bm{i}}$, $\mfe_{\bm{Z}}$, $\mfe_{F,C}$ in \S~\ref{ss:CCCC}), we have a decomposition of $\cX \subseteq \cY^{\mfq}$ as $\cX = (\coprod_{i \neq \bm{i}} \cX[\mfe^i]) \amalg (\cX[\mfe_\theta^{\bm{i}}] \amalg \cX[\mfe_{\bm{Z}}]) \amalg \cX[\mfe_{F,C}]$, where $\cX[\mfe] = \menge{x \in \cX}{x(\mfe) = 1}$. Define $\cY_{\rm conn}^{\mfq} \defeq (\coprod_{i \neq \bm{i}} \cX[\mfe^i]) \amalg \cX[\mfe_\theta^{\bm{i}}]$ and $\cY_{\rm rest}^{\mfq} \defeq \cX[\mfe_{\bm{Z}}] \amalg \cX[\mfe_{F,C}]$. We have $\cY_{\rm conn}^{\mfq} \subseteq \cY_\mfr^{\mfq}$ for $\mfr = 0,1$. According to the construction of $\beta_0^{\mfq} \defeq \beta_{n+1,0}^{\mfq}$ in \S~\ref{ss:(conn)}, $\bm{b}_0 \defeq \bm{b}_{n+1,0}$ sends $x \in \cY_{\rm conn}^{\mfq}$ to $x \in \cX$. To describe $\bm{b}_1 \defeq \bm{b}_{n+1,1}$ on $\cY_{\rm conn}^{\mfq}$, note that we have an identification
\begin{equation}
\label{e:X1X1=MX}
\Big( \coprod_{i \neq \bm{i}} \cX[\mfe^i] \Big) \amalg \cX[\mfe_\theta^{\bm{i}}] \isom \Big( \coprod_{i \neq \bm{i}} \Big( \coprod_j \cM(j,i) \times \cX_n^i \Big) \Big) \amalg \Big( \coprod_j \cM(j,\bm{i}) \times \cX_n^{\bm{i}} \Big)
= \Big( \coprod_{i \neq \bm{i}} \cM^i \times \cX_n^i \Big) \amalg (\cM^{\bm{i}} \times \cX_n^{\bm{i}})
\end{equation}
corresponding to the decomposition of $\pi \circ \bfdot{\varphi}_F$ in \S~\ref{ss:CCCC} (see \eqref{e:phiF} -- \eqref{e:piphi-bi} in combination with the definition of $\mfe^i$, $\mfe_\theta^{\bm{i}}$, $\mfe_{\bm{Z}}$, $\mfe_{F,C}$ in \S~\ref{ss:CCCC}), where $\cM^i = \coprod_j \cM(j,i)$ and $\cM^{\bm{i}} = \coprod_j \cM(j,\bm{i})$. With respect to \eqref{e:X1X1=MX}, if $\cM^i = \big\{ \mu^i_1, \dotsc, \mu^i_{\sum_j m(j,i)} \big\}$ corresponding to the enumeration of rank-one projections in $\bigoplus_j DM_{m(j,i)}$ in \S~\ref{ss:(conn)}, we have
\begin{equation}
\label{e:b1mux}
\bm{b}_1(\mu^i_l,x) = (\mu^i_{l-1},x) \ {\rm if} \ 2 \leq l \leq \sum_j m(j,i) \qquad \text{and} \ \bm{b}_1(\mu^i_1,x) = (\mu^i_{\sum_i m(j,i)},x) 
\qquad \qquad \forall \ x \in \cX_n^i,
\end{equation}
according to the construction of $\beta_1^{\mfq} \defeq \beta_{n+1,1}^{\mfq}$ in \S~\ref{ss:(conn)}. We also have $\cY_{\rm rest}^{\mfq} \subseteq \cY_{\mfr}^{\mfq}$ for $\mfr = 0,1$, and $\bm{b}_\mfr$ sends $x \in \cY_{\rm rest}^{\mfq}$ to $x \in \cX$ for $\mfr = 0,1$ according to the construction of $\beta_\mfr^{\mfq}$ in \S~\ref{ss:(conn)}.

On $( \coprod_{i \neq \bm{i}} \cX[\mfe^i] ) \amalg \cX[\mfe_\theta^{\bm{i}}]$, using the identification \eqref{e:X1X1=MX}, we have
\begin{equation}
\label{e:pmux}
\bm{p}(\mu,x) = x \in \cX_n^i \qquad \forall \ \mu \in \cM^i, \, x \in \cX_n^i
\end{equation}
according to the descriptions of the components of $\pi \circ \bfdot{\varphi}_F$ in \eqref{e:piphiji}, \eqref{e:piphijbitheta}, \eqref{e:piphi-i} and \eqref{e:piphi-bi}. Furthermore, note that condition \eqref{e:phiCfp} implies that we have embeddings
\begin{equation}
\label{e:YnYn+1}
\cY_n = \cY_n^{\mfp} \amalg \coprod_{p \neq \mfp} \bfdot{\cY}_n^p \into \cY_\mfr^{\mfq}, \qquad \mfr = 0,1,
\end{equation}
sending $\cY_{n,\mfr}$ into $\bm{b}_\mfr^{-1}(( \coprod_{i \neq \bm{i}} \cX[\mfe^i] ) \amalg \cX[\mfe_\theta^{\bm{i}}])$ such that the following diagram commutes for $\mfr = 0,1$:
\begin{equation}
\label{e:CDbnb}
\begin{tikzcd}
  \cY_{n,\mfr} \arrow[d, "\bm{b}_{n,\mfr}"'] \arrow[r, hook] & \bm{b}_\mfr^{-1}((\coprod_{i \neq \bm{i}} \cX[\mfe^i] ) \amalg \cX[\mfe_\theta^{\bm{i}}]) \arrow[d, "\bm{b}_\mfr"]
  \\
  \cX_n & \arrow[l, "\bm{p}"'] ( \coprod_{i \neq \bm{i}} \cX[\mfe^i] ) \amalg \cX[\mfe_\theta^{\bm{i}}]
\end{tikzcd}
\end{equation}

\bprop
\label{prop:conn}
The C*-diagonals $B_n$ as in Lemma~\ref{lem:conn:Ell} have path-connected spectra for all $n = 1, 2, 3, \dotsc$.
\eprop
\setlength{\parindent}{0cm} \setlength{\parskip}{0cm}

\bproof
In the following, for two points $x_1$ and $x_2$, we write $x_1 \pc x_2$ if there exists a continuous path from $x_1$ to $x_2$, in a space which will be clear from the context or specified otherwise. We start with the observation that given $\bm{x} \in \cF_n^{(0)} \setminus \cX_n$, i.e., $\bm{x} \in (Z_n \setminus \gekl{\theta^{\bm{i}}_n}) \times \cX_n^{\bm{i}}$, since $Z_n$ is path-connected, we have $\bm{x} \pc x \in \gekl{\theta_n^{\bm{i}}} \times \cX_n^{\bm{i}} \subseteq \cX_n$. Hence, to show that $\Spec B_n$ is path-connected, it suffices to show that $[[0,1] \times_\bullet \cY_n]$ is path-connected. Let us prove inductively on $n$ that $[[0,1] \times_\bullet \cY_n] \subseteq \Spec B_n$ is path-connected. Note that we can always make the following reduction: For all $(t,y) \in [0,1] \times \cY_n$, we have $y \in \cY_{n,0}$ because $\beta_{n,0}$ is always unital, and $[t,y] \pc [0,y]$. Moreover, given $\mfr \in \gekl{0,1}$ and $y \in \cY_{n,\mfr}$, since $\bm{b}_{n,0}: \: \cX_n \subseteq \cY_{n,0}^{\mfp} \to \cX_n$ is surjective by \eqref{e:bn0bp}, there exists $\bar{y} \in \cX_n \subseteq \cY_{n,0}^{\mfp}$ such that $\bm{b}_{n,0}(\bar{y}) = \bm{b}_{n,\mfr}(y)$ and thus $[0,\bar{y}] = [\mfr,y]$. Hence it is enough to show for $y, \bar{y} \in \cX_n \subseteq \cY_{n,0}^{\mfp}$ that $[0,y] \pc [0,\bar{y}]$.
\setlength{\parindent}{0cm} \setlength{\parskip}{0.5cm}

For $n=1$, this follows from the observation that we have $[0,x_{l+1}] \pc [1,x_{l+1}] = [0,x_l]$ for all $1 \leq l \leq [1,I] - 1$ because of \eqref{e:b11bp}. Now let us assume that $[[0,1] \times_\bullet \cY_n]$ is path-connected, and let us show that $[[0,1] \times_\bullet \cY_{n+1}]$ is path-connected. We use the same notation as in the description of $\bm{p}$ above (we also drop the index $n+1$). It suffices to show that for all $y, \bar{y} \in \cX \subseteq \cY_0^{\mfq}$ that $[0,y] \pc [0,\bar{y}]$. We further reduce to $y, \bar{y} \in \cY_{\rm conn}^{\mfq}$: If $y \in \cX[\mfe_{\bm{Z}}]$, then there exists $\mfr \in \gekl{0,1}$ and $\ti{y} \in \cY_\mfr$ with $\bm{b}_\mfr(\ti{y}) = \bm{b}_0(y) \ (= y)$, and by \eqref{e:phiCx}, we must have $\bm{b}_{\mfr^*}(\ti{y}) \in \bm{p}^{-1}(\cX_n) \cap \cX = ( \coprod_{i \neq \bm{i}} \cX[\mfe^i] ) \amalg \cX[\mfe_\theta^{\bm{i}}] = \cY_{\rm conn}^{\mfq}$, where $\mfr^* = 1 - \mfr$. Hence $[0,y] = [\mfr,\ti{y}] \pc [\mfr^*,\ti{y}] = [0,y']$ for some $y' \in \cY_{\rm conn}^{\mfq}$. If $y \in \cX[\mfe_{F,C}]$, then there must exist $\mfr \in \gekl{0,1}$ and $\ti{y} \in \cY_\mfr$ with $\bm{b}_\mfr(\ti{y}) = \bm{b}_0(y) \ (= y)$, and by \eqref{e:phiCl}, we must have $\bm{b}_{\mfr^*}(\ti{y}) \in \bm{p}^{-1}(\cX_n) \cap \cX = ( \coprod_{i \neq \bm{i}} \cX[\mfe^i] ) \amalg \cX[\mfe_\theta^{\bm{i}}] = \cY_{\rm conn}^{\mfq}$ (here $\cX$ is viewed as a subset of $\Spec B_{n+1}$), where $\mfr^* = 1 - \mfr$. Hence $[0,y] = [\mfr,\ti{y}] \pc [\mfr^*,\ti{y}] = [0,y']$ for some $y' \in \cY_{\rm conn}^{\mfq}$.

Moreover, given $y \in \cY_{\rm conn}^{\mfq}$ (for which we have $\bm{b}_0(y) = y \in \cX$), there exists by \eqref{e:CDbnb} $y' \in \cY_{n,0} \subseteq \cY_0$ such that $\bm{p}(\bm{b}_0(y')) = \bm{p}(\bm{b}_0(y))$. Viewing $\bm{b}_0(y')$ as an element in $\cY_{\rm conn}^{\mfq}$, let us now show that
\begin{equation}
\label{e:0bypc0by'}
[0,\bm{b}_0(y')] = [0,\bm{b}_0(y)]:
\end{equation}
Under the bijection \eqref{e:X1X1=MX}, we have $y = (\mu,x)$ and $\bm{b}_0(y') = (\mu',x)$, where $x = \bm{p}(\bm{b}_0(y')) = \bm{p}(\bm{b}_0(y))$. Hence \eqref{e:0bypc0by'} follows from the following claim:
\begin{equation}
\label{e:0mxpc0m'x}
\text{Under the bijection} \ \eqref{e:X1X1=MX}, \ \text{we have} \ [0,(\mu,x)] \pc [0,(\mu',x)] \ \text{for all} \ \mu, \mu' \in \cM^i, \, x \in \cX_n^i.
\end{equation}
This in turn follows from the observation that for all $l \in \gekl{1, \dotsc, (\sum_j m(j,i)) - 1}$ and $x \in \cX_n^i$, we have $[0,(\mu_{l+1},x)] \pc [1,(\mu_{l+1},x)] = [0,(\mu_l,x)]$. The last equation follows from \eqref{e:b1mux}. So we have $[0,y] \pc [0,\bm{b}_0(y')] = [0,y']$. Hence it suffices to show $[0,y] \pc [0,\bar{y}]$ in $[[0,1] \times_\bullet \cY_{n+1}] \subseteq \Spec B_{n+1}$ for all $y, \bar{y} \in \cY_{n,0}$.

By induction hypothesis, we have $[0,y] \pc [0,\bar{y}]$ in $[[0,1] \times_\bullet \cY_n] \subseteq \Spec B_n$. Hence there exist $(\mfr_k,y_k) \in \gekl{0,1} \times \cY_n$, $0 \leq k \leq K$, such that $(\mfr_0,y_0) = (0,y)$, $(\mfr_K,y_K) = (0,\bar{y})$ and for all $0 \leq k \leq K - 1$, we have $[\mfr_k,y_k] = [\mfr_{k+1},y_{k+1}]$ in $\Spec B_n$ or $y_k = y_{k+1}$, $\mfr_{k+1} = \mfr_k^*$ (where $\mfr_k^* = 1 - \mfr_k$). Clearly, in the latter case, we have $[\mfr_k,y_k] \pc [\mfr_k^*,y_k] = [\mfr_{k+1},y_{k+1}]$ in $[[0,1] \times_\bullet \cY_{n+1}] \subseteq \Spec B_{n+1}$. To treat the former case, we need to show that $[\mfr_k,y_k] = [\mfr_{k+1},y_{k+1}]$ in $\Spec B_n$ (i.e., $\bm{b}_{n,\mfr_k}(y_k) = \bm{b}_{n,\mfr_{k+1}}(y_{k+1})$) implies $[\mfr_k,y_k] \pc [\mfr_{k+1},y_{k+1}]$ in $[[0,1] \times_\bullet \cY_{n+1}] \subseteq \Spec B_{n+1}$, where we view $y_k$ and $y_{k+1}$ as elements of $\cY_{\mfr_k}^{\mfq}$ and $\cY_{\mfr_{k+1}}^{\mfq}$ using \eqref{e:YnYn+1}. We have 
$$
\bm{p}(\bm{b}_{\mfr_k}(y_k)) \overset{\eqref{e:CDbnb}}{=} \bm{b}_{n,\mfr_k}(y_k) = \bm{b}_{n,\mfr_{k+1}}(y_{k+1}) = \bm{p} (\bm{b}_{\mfr_{k+1}}(y_{k+1})).
$$
Thus, viewing $\bm{b}_{\mfr_k}(y_k)$, $\bm{b}_{\mfr_{k+1}}(y_{k+1})$ as elements of $\cY_{\rm conn}^{\mfq}$, we have $\bm{b}_{\mfr_k}(y_k) = (\mu,x)$ and $\bm{b}_{\mfr_{k+1}}(y_{k+1}) = (\mu',x)$ for some $\mu, \mu' \in \cM^i$ and $x \in \cX_n^i$ with respect to \eqref{e:X1X1=MX}. Hence \eqref{e:0mxpc0m'x} implies that, in $[[0,1] \times_\bullet \cY_{n+1}] \subseteq \Spec B_{n+1}$, we have
\begin{equation*}
[\mfr_k,y_k] = [0, \bm{b}_{\mfr_k}(y_k)] = [0, (\mu,x)] \pc [0, (\mu',x)] = [0, \bm{b}_{\mfr_{k+1}}(y_{k+1})] = [\mfr_{k+1},y_{k+1}]. \qedhere
\end{equation*}
\eproof

\bremark
\label{rem:Yconnconn}
The proof of \eqref{e:0mxpc0m'x} yields that for all $y, \bar{y} \in \cY_{\rm conn}^{\mfq}$ with $\bm{p}[0,y] = \bm{p}[0,\bar{y}]$, there exists a continuous path $\xi$ in $[[0,1] \times_\bullet \cY_{n+1}]$ with $\xi(0) = [0,y]$, $\xi(1) = [0,\bar{y}]$ and $\bm{p} \circ \xi \equiv \bm{p}[0,y] = \bm{p}[0,\bar{y}]$.
\eremark
\setlength{\parindent}{0cm} \setlength{\parskip}{0.5cm}

\bcor
\label{cor:conn_unital}
In the unital case, modification (conn) yields C*-diagonals with connected spectra.
\ecor
\setlength{\parindent}{0cm} \setlength{\parskip}{0cm}

\bproof
The C*-diagonal is given by $B = \ilim_n \gekl{B_n, \varphi_n}$, so that its spectrum is $\Spec B \cong \plim_n \gekl{\Spec B_n, \bm{p}_n}$. In the unital case, $B_n$ is unital for all $n$, so that $\Spec B_n$ is compact for all $n$. By Proposition~\ref{prop:conn}, $\Spec B_n$ is path-connected, in particular connected. Now our claim follows from the general fact that inverse limits of compact connected spaces are again connected (see for instance \cite[Theorem~6.1.20]{Eng}). 
\eproof
\setlength{\parindent}{0cm} \setlength{\parskip}{0.5cm}

In the stably projectionless case, we cannot argue as for Corollary~\ref{cor:conn_unital} because it is no longer true in general that inverse limits of locally compact, non-compact, connected spaces are again connected. Instead, by conjugating $\beta_{n+1,\bullet}$ by suitable permutation matrices and adjusting $\varphi$ accordingly, we can always arrange that the $\lambda$s in \eqref{e:fplp} are monotonous and that, in addition to \eqref{e:phiCfp} -- \eqref{e:phiCx}, we have the following:
\begin{align}
\label{e:phiCl*}
& \forall \ \lambda, \ \text{corresponding block diagonal entry} \ f^p \circ \lambda \ \text{in} \ \varphi_C(f,a) \ \text{as in} \ \eqref{e:fplp}, \ \mfr \ \text{as in} \ \eqref{e:phiCl} \ \text{with} \ \lambda(\mfr) \notin \gekl{0,1} \tag{11$^*$} \\
& \exists \ \text{a block diagonal entry} \ f^p \circ \lambda^* \ \text{in} \ \varphi_C(f,a) \ \text{as in} \ \eqref{e:fplp} \ \text{with} \ \lambda^*(\mfr^*) = \lambda(\mfr), \, \lambda^*(\mfr) = \lambda(\mfr^*)^*, \nonumber \\ 
& \text{unless} \ p = \grave{p}, \ \text{in which case} \ \lambda(\mfr^*) = 0; \nonumber \\
\label{e:phiClx}
& \forall \ \lambda, \mfr \ \text{as in} \ \eqref{e:phiCl} \ \text{with} \ \bm{t} \defeq \lambda(\mfr) \notin \gekl{0,1} \ \text{and the corresponding block diagonal entry} \ f^p \circ \lambda \ \text{in} \ \eqref{e:fplp}, \tag{11$\reg$}\\
& \text{we have that} \ f^p(\bm{t}) \ \text{appears as exactly one block diagonal entry in} \ \varphi_{F,C}(f) \ \text{in} \ \eqref{e:phiF}. \nonumber
\end{align}
\bprop
\label{prop:conn_spl}
In the stably projectionless case, modification (conn) with the above-mentioned adjustments yields C*-diagonals with connected spectra.
\eprop
\setlength{\parindent}{0cm} \setlength{\parskip}{0cm}

\bproof
The C*-diagonal is given by $B = \ilim_n \gekl{B_n, \varphi_n}$, so that its spectrum is $\Spec B \cong \plim_n \gekl{\Spec B_n, \bm{p}_n}$. Let $\bm{p}_{n,\infty}: \: \Spec B \onto \Spec B_n$ be the canonical map from the inverse limit structure of $\Spec B$, and denote by $\bm{p}_{n,\bar{n}}: \: \Spec B_{\bar{n}+1} \onto \Spec B_n$ the composition $\bm{p}_{\bar{n}} \circ \dotso \circ \bm{p}_n$. Now define for each $N \geq 1$ the intervals $I_y \defeq [0,1]$ for $y \in \cY_{1,1}$, $I_y \defeq [0,1 - \tfrac{1}{N}]$ for $y \notin \cY_{1,1}$, and the subset $K_{N,1} \defeq [(\bigcup_{y \in \cY_1} I_y \times \gekl{y}) \amalg \cF_1^{(0)}] \subseteq \Spec B_1$. Now it is straightforward to check by induction on $n$ that $\bm{p}_{1,n}^{-1}(K_{N,1}) = (\bigcup_{y \in \fY} [\fI_y \times \gekl{y}]) \cup (\bigcup_{x \in \fX} [Z_n \times \gekl{x}])$ where $\fY$ is a subset of $\cY_n$, $\fI_y$ is of the form $[0,1]$, $[0,t]$ or $[t,1]$ for some $t \in [0,1]$, $\fX$ is a subset of $\cX_n^{\bm{i}}$, for all $\ti{y} \in \fY$ with $\fI_{\ti{y}} \neq [0,1]$ there exists $y \in \fY$ with $\fI_y = [0,1]$ and $[\fI_{\ti{y}} \times \gekl{\ti{y}}] \cap [\fI_y \times \gekl{y}] \neq \emptyset$, and for all $x \in \fX$ there exists $y \in \fY$ with $\fI_y = [0,1]$ and $[Z_n \times \gekl{x}] \cap [\fI_y \times \gekl{y}] \neq \emptyset$. Now we proceed inductively on $n$ to show that $\bm{p}_{1,n}^{-1}(K_{N,1})$ is path-connected for all $n$. The case $n=1$ is checked as in Proposition~\ref{prop:conn}. For the induction step, first reduce as in Proposition~\ref{prop:conn} to showing that $\bigcup_y [\fI_y \times \gekl{y}]$ is path-connected, where the union is taken over all $y \in \fY$ with $\fI_y = [0,1]$. Further reduce as in Proposition~\ref{prop:conn} to the statement that for all $y, \bar{y} \in \fY$ with $\fI_y = [0,1]$, $\fI_{\bar{y}} = [0,1]$ and $y, \bar{y} \in \cY_{\rm conn}^\mfq$ that $[0,y] \pc [0,\bar{y}]$ in $\bm{p}_{1,n}^{-1}(K_{N,1})$. Here the case $y \in \cX[\mfe_{\bm{Z}}]$ is treated as in Proposition~\ref{prop:conn}, while the case $y \in \cX[\mfe_{F,C}]$ uses \eqref{e:phiCl*} and \eqref{e:phiClx}. Now use the induction hypothesis as in Proposition~\ref{prop:conn} to show that we indeed have $[0,y] \pc [0,\bar{y}]$ in $\bm{p}_{1,n}^{-1}(K_{N,1})$ for all $y, \bar{y} \in \fY$ with $\fI_y = [0,1]$, $\fI_{\bar{y}} = [0,1]$ and $y, \bar{y} \in \cY_{\rm conn}^\mfq$. As $\bm{p}_{1,n}$ is proper, $\bm{p}_{1,n}^{-1}(K_{N,1})$ is compact. Hence it follows that $K_N \defeq \bm{p}_{1,\infty}^{-1}(K_{N,1}) \cong \plim_n \big\{ \bm{p}_{1,n}^{-1}(K_{N,1}), \bm{p}_n \big\}$ is connected (see for instance \cite[Theorem~6.1.20]{Eng}). Therefore $\Spec B = \bigcup_N K_N$ is connected as it is the increasing union of connected subsets. 
\eproof

All in all, we obtain the following
\btheo
\label{thm:conn_Ell}
For every prescribed Elliott invariant $(G_0, G_0^+, u, T, r, G_1)$ as in \cite[Theorem~1.2]{Li18}, our construction produces a twisted groupoid $(G,\Sigma)$ with the same properties as in \cite[Theorem~1.2]{Li18} (in particular, $C^*_r(G,\Sigma)$ is a classifiable unital C*-algebra with ${\rm Ell}(C^*_r(G,\Sigma)) \cong (G_0, G_0^+, u, T, r, G_1)$) such that $G$ has connected unit space.
\setlength{\parindent}{0.5cm} \setlength{\parskip}{0cm}

For every prescribed Elliott invariant $(G_0, T, \rho, G_1)$ as in \cite[Theorem~1.3]{Li18}, our construction produces a twisted groupoid $(G,\Sigma)$ with the same properties as in \cite[Theorem~1.3]{Li18} (in particular, $C^*_r(G,\Sigma)$ is classifiable stably projectionless with continuous scale, and ${\rm Ell}(C^*_r(G,\Sigma)) \cong (G_0, \gekl{0}, T, \rho, G_1)$) such that $G$ has connected unit space.
\etheo

This theorem, in combination with classification results for all classifiable C*-algebras, implies Theorem~\ref{thm:main1}.
\setlength{\parindent}{0cm} \setlength{\parskip}{0.5cm}

\section{Further modification of the construction leading to the path-lifting property}
\label{s:FurtherMod}

Let us now present a further modification of the construction recalled in \S~\ref{ss:CCCC} which will allow us to produce C*-diagonals with Menger manifold spectra. We focus on constructing classifiable C*-algebras (unital or stably projectionless with continuous scale) with torsion-free $K_0$ and trivial $K_1$. In that case, the construction recalled in \S~\ref{ss:CCCC} simplifies because $F_n = \hat{F}_n$ for all $n$, so that we can (and will) think of $A_n$ as a sub-C*-algebra of $C([0,1],E_n)$.

\subsection{Modification (path)}
\label{ss:(path)}

Suppose that we are given a tuple $\cE = (G_0, G_0^+, u, T, r, G_1)$ as in \cite[Theorem~1.2]{Li18} or $\cE = (G_0, \gekl{0}, T, \rho, G_1)$ as in \cite[Theorem~1.3]{Li18} which we want to realize as the Elliott invariant of a classifiable C*-algebra, with $G_0$ torsion-free and $G_1 = \gekl{0}$. As explained in \cite[\S~2]{Li18}, the construction recalled in \S~\ref{ss:CCCC} proceeds in two steps. First an inductive system $\{ \bfdot{A}_n, \accir{\varphi}_n \}$ is constructed so that $\ilim_n \{ \bfdot{A}_n, \accir{\varphi}_n \}$ has the desired Elliott invariant, but is not simple, and then a further modification yields an inductive system $\{ \bfdot{A}_n, \bfdot{\varphi}_n \}$ such that $\ilim_n \{ \bfdot{A}_n, \bfdot{\varphi}_n \}$ has the same Elliott invariant and in addition is simple.

The first step in our modification (path) is as in the previous modification (conn) (see \S~\ref{ss:(conn)}) and produces the first building block $A_1$. Now suppose that we have produced 
$$
A_1 \overset{\varphi_1}{\lori} A_2 \overset{\varphi_2}{\lori} \dotso \overset{\varphi_{n-1}}{\lori} A_n,
$$
and that the first step of the original construction as in \cite[\S~2]{Li18} yields $\accir{\varphi}_n: \: A_n \to \bfdot{A}_{n+1}$. We modify $\accir{\varphi}_n$ in two steps, first to $\bfdot{\varphi}_n: \: A_n \to \bfdot{A}_{n+1}$, then to $\varphi_n: \: A_n \to A_{n+1}$. Let us start with the first step. We use the same notation as in \S~\ref{ss:CCCC} and \S~\ref{ss:(conn)}. Recall the description of $\beta_{n,\mfr}^{p,i}$ in \eqref{e:beta=}; it is a composition of the form
$$
  F_n^i \overset{1 \otimes \id_{F_n^i}}{\lori} 1_{m_\mfr(p,i)} \otimes F_n^i \subseteq M_{m_\mfr(p,i)} \otimes F_n^i \tailarr E_n^p.
$$
Here and in the sequel, an arrow $\tailarr$ denotes a *-homomorphism of multiplicity $1$ sending diagonal matrices to diagonal matrices as before. Let $\psi_n: \: F_n \to F_{n+1}$ be as in \cite[\S~2]{Li18}. The map
$$
  \psi_n^{j,i}: \: F_n^i \into F_n \overset{\psi_n}{\lori} F_{n+1} \onto F_{n+1}^j
$$
is given by the following composition:
\begin{equation*}
F_n^i \overset{1 \otimes \id_{\hat{F}_n^i}}{\lori} 1_{m(j,i)} \otimes F_n^i \subseteq M_{m(j,i)} \otimes F_n^i \tailarr F_{n+1}^j.
\end{equation*}
By choosing $G'$ in \cite[\S~2]{Li18} suitably and because of \cite[Inequality~(2)]{Li18}, we can always arrange that there exist pairwise distinct indices $\gekl{j_0^p}_p \cup \gekl{j_1^p}_{p \neq \grave{p}}$ such that we have $m(j_\mfr^p,i) \geq m_\mfr(p,i)$ for all $p$, $i$, $\mfr = 0,1$ ($p \neq \grave{p}$  if $\mfr = 1$). Then for suitable embeddings $E_n^p \tailarr F_{n+1}^{j_\bullet^p}$ sending $DE_n^p$ into $DF_{n+1}^{j_\bullet^p}$, $\psi^{j_\bullet^p}$ is of the form $F_n \to E_n^p \oplus \bar{F}_{n+1}^{j_\bullet^p} \tailarr F_{n+1}^{j_\bullet^p}$, for some finite-dimensional algebra $\bar{F}_{n+1}^{j_\bullet^p}$, where the first map is given by
$
\rukl{
\begin{smallmatrix}
\beta_{n,\bullet}^p & 0 \\
0 & \bar{\psi}^{j_\bullet^p}
\end{smallmatrix}
}
$
for some map $\bar{\psi}^{j_\bullet^p}: \: F_n \to \bar{F}_{n+1}^{j_\bullet^p}$. Let 
$\ve_\beta^{j_\bullet^p} \defeq 
\rukl{
\begin{smallmatrix}
\beta_{n,\bullet}^p(1_{F_n}) & 0 \\
0 & 0
\end{smallmatrix}
}
$,
viewed as a projection in $F_{n+1}^{j_\bullet^p}$ via the second embedding $E_n^p \oplus \bar{F}_{n+1}^{j_\bullet^p} \tailarr F_{n+1}^{j_\bullet^p}$.

We start discussing the connecting map and will drop indices whenever convenient. $\accir{\varphi} \defeq \accir{\varphi}_n: \: A_n \to \bfdot{A}_{n+1}$ is given by $\accir{\varphi}_F: \: A_n \overset{\accir{\varphi}}{\lori} \bfdot{A}_{n+1} \to F_{n+1}$, $\accir{\varphi}_F(f,a) = \psi(a)$, and $\accir{\varphi}_C: \: A_n \overset{\accir{\varphi}}{\lori} \bfdot{A}_{n+1} \to C([0,1],\bfdot{E}_{n+1})$, 
$\accir{\varphi}_C(f,a) = 
\rukl{
\begin{smallmatrix}
\Phi(f) & 0 \\
0 & \Phi_F(a)
\end{smallmatrix}
}
$.
Let $\ve_\Phi$ be the smallest projection in $D\bfdot{E}$ such that $\Phi(f)(t) = \ve_\Phi \cdot \Phi(f)(t) \cdot \ve_\Phi$ for all $t \in [0,1]$, and let $\ve_{C,F} \in D\bfdot{E}$ be such that $\Phi_F(1_{F_n}) \equiv \ve_{C,F}$. We have a decomposition $\ve_\Phi = \sum_{q,p} \ve_\Phi^{q,p}$, $\ve_\Phi^{q,p} = \ve^+_{q,p} + \ve_+^{q,p} + \ve^-_{q,p} + \ve_-^{q,p}$ into pairwise orthogonal projections in $D\bfdot{E}$ such that, for all $q, p$,
\begin{align*}
\ve^+_{q,p} \cdot \Phi(f) \cdot \ve^+_{q,p} &= e^+_{q,p} \otimes f^p,
&\ve_+^{q,p} \cdot \Phi(f) \cdot \ve_+^{q,p} &= e_+^{q,p} \otimes f^p,\\
\ve^-_{q,p} \cdot \Phi(f) \cdot \ve^-_{q,p} &= e^-_{q,p} \otimes f^p \circ (1 - \id),
&\ve_-^{q,p} \cdot \Phi(f) \cdot \ve_-^{q,p} &= e_-^{q,p} \otimes f^p \circ (1 - \id),
\end{align*}
for some finite-rank projections $e^+_{q,p}, e_+^{q,p}, e^-_{q,p}, e_-^{q,p}$ encoding multiplicities of block diagonal entries in $\Phi$. In the unital case, we can always arrange
\begin{equation}
\label{e:e>1_unital}
\rk e^+_{q,p}, \rk e_+^{q,p}, \rk e^-_{q,p}, \rk e_-^{q,p} \geq 1 \quad \forall \ q, p. 
\end{equation}
In the stably projectionless case, we can always arrange that 
\begin{equation}
\label{e:e>1_spl1}
\rk e^+_{q,p}, \rk e_+^{q,p}, \rk e^-_{q,p}, \rk e_-^{q,p} \geq 1 \quad \forall \ q \neq \grave{q}, \, p \neq \grave{p}, \qquad \text{and} \ \rk e^+_{\grave{q},\grave{p}} \geq 1, 
\end{equation}
as well as $\rk e^+_{q,p}, \rk e_+^{q,p}, \rk e^-_{q,p}, \rk e_-^{q,p} = 0$ for all $q = \grave{q}, p \neq \grave{p}$ or $q \neq \grave{q}, p = \grave{p}$, and $\rk e_+^{\grave{q},\grave{p}}, \rk e^-_{\grave{q},\grave{p}}, \rk e_-^{\grave{q},\grave{p}} = 0$.

$\beta_\mfr^{q,j}$ is a composition as in \eqref{e:beta=} of the form $F^j \overset{1 \otimes \id_{F^j}}{\lori} 1_{m_\mfr(q,j)} \otimes F^j \subseteq M_{m_\mfr(q,j)} \otimes F^j \tailarr \bfdot{E}^q$. By replacing $\bfdot{E}^q$ by $M_{{n+1,q} + N \cdot [n+1,J]}$ containing $\bfdot{E}^q \oplus F^{\oplus N}$ in the canonical way, and by replacing $\beta_\mfr^q$ by $\beta_\mfr^q \oplus \id_{F^{\oplus N}}$ as in modification (conn), we can arrange that, for all $q,p$,
\begin{equation*}
m_0(q,j_0^p) \geq \rk e^+_{q,p}, \quad
m_1(q,j_1^p) \geq \rk e_+^{q,p}, \quad
m_1(q,j_0^p) \geq \rk e^-_{q,p}, \quad
m_0(q,j_1^p) \geq \rk e_-^{q,p}.
\end{equation*}
By further enlarging $\bfdot{E}^q$ as above, and by conjugating $\beta_\mfr^q$ by suitable permutation matrices if necessary, we can arrange that there exist a decomposition $\ve_{C,F} = (\sum_{q,p} \underline{\ve}^{q,p}) + (\sum_{q,p} \bar{\ve}_{q,p}) + \ve_{\rm const}$ into pairwise orthogonal projections in $D\bfdot{E}$ such that for all $q,p$ and $\mfr = 0,1$, 
\begin{align*}
\beta_\mfr^q \circ (\ve_\beta^{j_\mfs^p} \cdot \psi^{j_\mfs^p} \cdot \ve_\beta^{j_\mfs^p}) 
=& \ \ve^+_{q,p} \cdot ( \beta_\mfr^q \circ (\ve_\beta^{j_\mfs^p} \cdot \psi^{j_\mfs^p} \cdot \ve_\beta^{j_\mfs^p}) ) \cdot \ve^+_{q,p} 
+ \ve_+^{q,p} \cdot ( \beta_\mfr^q \circ (\ve_\beta^{j_\mfs^p} \cdot \psi^{j_\mfs^p} \cdot \ve_\beta^{j_\mfs^p}) ) \cdot \ve_+^{q,p}\\
&+ \ve^-_{q,p} \cdot ( \beta_\mfr^q \circ (\ve_\beta^{j_\mfs^p} \cdot \psi^{j_\mfs^p} \cdot \ve_\beta^{j_\mfs^p}) ) \cdot \ve^-_{q,p} 
+ \ve_-^{q,p} \cdot ( \beta_\mfr^q \circ (\ve_\beta^{j_\mfs^p} \cdot \psi^{j_\mfs^p} \cdot \ve_\beta^{j_\mfs^p}) ) \cdot \ve_-^{q,p}\\
&+ \underline{\ve}^{q,p} \cdot ( \beta_\mfr^q \circ (\ve_\beta^{j_\mfs^p} \cdot \psi^{j_\mfs^p} \cdot \ve_\beta^{j_\mfs^p}) ) \cdot \underline{\ve}^{q,p} + \bar{\ve}_{q,p} \cdot ( \beta_\mfr^q \circ (\ve_\beta^{j_\mfs^p} \cdot \psi^{j_\mfs^p} \cdot \ve_\beta^{j_\mfs^p}) ) \cdot \bar{\ve}_{q,p},
\end{align*}
and pairwise orthogonal finite-rank projections $\underline{e}^{q,p}, e_{\scriptscriptstyle (\diagup)}^{q,p}, e_{\scriptscriptstyle (\diagdown)}^{q,p}, \bar{e}_{q,p}, e^{\scriptscriptstyle (\diagup)}_{q,p}, e^{\scriptscriptstyle (\diagdown)}_{q,p}$ encoding multiplicities of block diagonal entries in $\Phi$, such that we have, for all $q,p$, 
\begin{align*}
&\ve^+_{q,p} \cdot ( \beta_0^q \circ (\ve_\beta^{j_0^p} \cdot \psi^{j_0^p} \cdot \ve_\beta^{j_0^p}) ) \cdot \ve^+_{q,p}
= e^+_{q,p} \otimes \beta_{n,0}^p, \ \
\underline{\ve}^{q,p} \cdot ( \beta_0^q \circ (\ve_\beta^{j_0^p} \cdot \psi^{j_0^p} \cdot \ve_\beta^{j_0^p}) ) \cdot \underline{\ve}^{q,p}
= \underline{e}^{q,p} \otimes \beta_{n,0}^p + e_{\scriptscriptstyle (\diagdown)}^{q,p} \otimes \beta_{n,0}^p,\\
&\ve^-_{q,p} \cdot ( \beta_1^q \circ (\ve_\beta^{j_0^p} \cdot \psi^{j_0^p} \cdot \ve_\beta^{j_0^p}) ) \cdot \ve^-_{q,p}
= e^-_{q,p} \otimes \beta_{n,0}^p, \ \
\underline{\ve}^{q,p} \cdot ( \beta_1^q \circ (\ve_\beta^{j_0^p} \cdot \psi^{j_0^p} \cdot \ve_\beta^{j_0^p}) ) \cdot \underline{\ve}^{q,p}
= \underline{e}^{q,p} \otimes \beta_{n,0}^p + e_{\scriptscriptstyle (\diagup)}^{q,p} \otimes \beta_{n,0}^p,\\
&\ve_-^{q,p} \cdot ( \beta_0^q \circ (\ve_\beta^{j_1^p} \cdot \psi^{j_1^p} \cdot \ve_\beta^{j_1^p}) ) \cdot \ve_-^{q,p}
= e_-^{q,p} \otimes \beta_{n,1}^p, \ \
\bar{\ve}_{q,p} \cdot ( \beta_0^q \circ (\ve_\beta^{j_1^p} \cdot \psi^{j_1^p} \cdot \ve_\beta^{j_1^p}) ) \cdot \bar{\ve}_{q,p}
= \bar{e}_{q,p} \otimes \beta_{n,1}^p + e^{\scriptscriptstyle (\diagup)}_{q,p} \otimes \beta_{n,1}^p,\\
&\ve_+^{q,p} \cdot ( \beta_1^q \circ (\ve_\beta^{j_1^p} \cdot \psi^{j_1^p} \cdot \ve_\beta^{j_1^p}) ) \cdot \ve_+^{q,p}
= e_+^{q,p} \otimes \beta_{n,1}^p, \ \
\bar{\ve}_{q,p} \cdot ( \beta_1^q \circ (\ve_\beta^{j_1^p} \cdot \psi^{j_1^p} \cdot \ve_\beta^{j_1^p}) ) \cdot \bar{\ve}_{q,p}
= \bar{e}_{q,p} \otimes \beta_{n,1}^p + e^{\scriptscriptstyle (\diagdown)}_{q,p} \otimes \beta_{n,1}^p,
\end{align*}
and $\ve \cdot ( \beta_\mfr^q \circ (\ve_\beta^{j_\mfs^p} \cdot \psi^{j_\mfs^p} \cdot \ve_\beta^{j_\mfs^p}) ) \cdot \ve = 0$ for all remaining choices of $\mfr, \mfs \in \gekl{0,1}$ and $\ve \in \gekl{\ve^+_{q,p}, \ve_+^{q,p}, \ve^-_{q,p}, \ve_-^{q,p}, \bar{\ve}_{q,p}, \underline{\ve}^{q,p}}$. In the stably projectionless case, we have $\bar{\ve}_{q,\grave{p}} = 0$ for all $q$ by arrangement. Moreover, we can always arrange that
\begin{equation}
\label{e:e>1_spl2}
\rk e_{\scriptscriptstyle (\diagdown)}^{\grave{q},\grave{p}} + \rk e_{\scriptscriptstyle (\diagup)}^{\grave{q},\grave{p}} \geq 1.
\end{equation}

Now define $\bfdot{\varphi} = \bfdot{\varphi}_n: \: A_n \to \bfdot{A}_{n+1}$ by setting $\bfdot{\varphi}_F^j: \: A_n \overset{\bfdot{\varphi}}{\lori} \bfdot{A}_{n+1} \to F_{n+1} \onto F_{n+1}^j$ and $\bfdot{\varphi}_C: \: A_n \overset{\bfdot{\varphi}}{\lori} \bfdot{A}_{n+1} \onto C([0,1],\bfdot{E}_{n+1})$ as follows:
\begin{eqnarray}
\label{e:phiF=1/2}
\bfdot{\varphi}_F^{j_\mfr^p}(f,a) &\defeq& 
\rukl{
\begin{smallmatrix}
f^p(\half) & 0 \\
0 & \bar{\psi}^{j_\mfr^p}(a)
\end{smallmatrix}
}, 
\quad {\rm and} \ \bfdot{\varphi}_F^j(f,a) \defeq \accir{\varphi}_F^j(f,a) \ {\rm for} \ j \notin \gekl{j_0^p, j_1^p};\\
\label{e:phiC_path}
\bfdot{\varphi}_C &=& \sum_{q,p} \Big( \ve^+_{q,p} \cdot \bfdot{\varphi}_C \cdot \ve^+_{q,p} + \ve_+^{q,p} \cdot \bfdot{\varphi}_C \cdot \ve_+^{q,p} + \ve^-_{q,p} \cdot \bfdot{\varphi}_C \cdot \ve^-_{q,p} + \ve_-^{q,p} \cdot \bfdot{\varphi}_C \cdot \ve_-^{q,p}\\
&&+ \ \underline{\ve}^{q,p} \cdot \bfdot{\varphi}_C \cdot \underline{\ve}^{q,p} + \bar{\ve}_{q,p} \cdot \bfdot{\varphi}_C \cdot \bar{\ve}_{q,p} \Big) + \ve_{\rm const} \cdot \bfdot{\varphi}_C \cdot \ve_{\rm const}; \nonumber \\
\ve^+_{q,p} \cdot \bfdot{\varphi}_C(f,a) \cdot \ve^+_{q,p} &\defeq& e^+_{q,p} \otimes f^p \circ (\half + \half \cdot \id), \quad 
\ve_+^{q,p} \cdot \bfdot{\varphi}_C(f,a) \cdot \ve_+^{q,p} \defeq e_+^{q,p} \otimes f^p \circ (\half \cdot \id), \nonumber \\
\ve^-_{q,p} \cdot \bfdot{\varphi}_C(f,a) \cdot \ve^-_{q,p} &\defeq& e^-_{q,p} \otimes f^p \circ (1 - \half \cdot \id), \quad
\ve_-^{q,p} \cdot \bfdot{\varphi}_C(f,a) \cdot \ve_-^{q,p} \defeq e_-^{q,p} \otimes f^p \circ (\half - \half \cdot \id); \nonumber \\
\underline{\ve}_{q,p} \cdot \bfdot{\varphi}_C(f,a) \cdot \underline{\ve}_{q,p}
&\defeq& \underline{e}^{q,p} \otimes f^p (\half) + e_{\scriptscriptstyle (\diagdown)}^{q,p} \otimes f^p \circ (\half - \half \cdot \id) + e_{\scriptscriptstyle (\diagup)}^{q,p} \otimes f^p \circ (\half \cdot \id), \nonumber \\
\bar{\ve}_{q,p} \cdot \bfdot{\varphi}_C(f,a) \cdot \bar{\ve}_{q,p}
&\defeq& \bar{e}_{q,p} \otimes f^p (\half) + e^{\scriptscriptstyle (\diagup)}_{q,p} \otimes f^p \circ (\half + \half \cdot \id) + e^{\scriptscriptstyle (\diagdown)}_{q,p} \otimes f^p \circ (1 - \half \cdot \id); \nonumber \\
\ve_{\rm const} \cdot \bfdot{\varphi}_C \cdot \ve_{\rm const} &\defeq& \ve_{\rm const} \cdot \accir{\varphi}_C \cdot \ve_{\rm const}. \nonumber 
\end{eqnarray}

Let us now continue with the second step and modify $\bfdot{\varphi}_n$ to $\varphi_n: \: A_n \to A_{n+1}$. This second step in our modification proceeds exactly in the same way as modification (conn), with the following difference: The embeddings $E_n^p \subseteq F_{n+1}^{j_\bullet^p}$, lead to the embedding $(\bigoplus_p E_n^p) \oplus (\bigoplus_{p \neq \grave{p}} E_n^p) \subseteq (\bigoplus_p F_{n+1}^{j_0^p}) \oplus (\bigoplus_{p \neq \grave{p}} F_{n+1}^{j_1^p}) \subseteq F_{n+1} \subseteq E_{n+1}^\mfq$, where $E_{n+1}^\mfq$ and the embedding $F_{n+1} \subseteq E_{n+1}^\mfq$ are constructed as in modification (conn). Let $\mfe_{EE} \in E_{n+1}^\mfq$ be the image of the unit of $(\bigoplus_p E_n^p) \oplus (\bigoplus_{p \neq \grave{p}} E_n^p)$ under the above embedding. Let $w_{EE}$ be a permutation matrix in $E_{n+1}^\mfq$ inducing the flip automorphism on $E_n^p \oplus E_n^p$ (i.e., the automorphism $E_n^p \oplus E_n^p \isom E_n^p \oplus E_n^p, \, (e,e') \ma (e',e)$) for all $p \neq \grave{p}$. Using the same notation as in modification (conn), note that $\mfe_{EE} \leq \mfe_{F,C}$, and define $w_{F,C} \defeq \mfe_{EE} \cdot w_{EE} \cdot \mfe_{EE} + (\mfe_{F,C} - \mfe_{EE})$. Define $w_{F,F}$ as in modification (conn). Now replace $w$ defined by \eqref{e:w=1w} in modification (conn) by 
$w \defeq 
\rukl{
\begin{smallmatrix}
w_{F,C} & 0 \\
0 & w_{F,F}
\end{smallmatrix}
}
$.
Furthermore, define $\beta_{n+1}^\mfq$, $A_{n+1}$ and $\varphi_n$ in the same way as in modification (conn).

Now it is straightforward to check that $\varphi_n$ is well-defined, i.e., $\varphi_n(f,a)$ satisfies the defining boundary conditions for $A_{n+1}$ for all $(f,a) \in A_n$. Proceeding recursively in this way, we obtain an inductive system $\gekl{A_n, \varphi_n}_n$.

\blemma
\label{lem:path:Ell}
$A = \ilim_n \gekl{A_n, \varphi_n}$ is a classifiable C*-algebra with ${\rm Ell}(A) \cong \cE$. In the stably projectionless case $A$ has continuous scale. If we define $B_n \defeq \menge{(f,a) \in A_n}{f(t) \in DE_n \ \forall \ t \in [0,1], \, a \in DF_n}$, then $B \defeq \ilim_n \gekl{B_n, \varphi_n}$ is a C*-diagonal of $A$.
\elemma
\setlength{\parindent}{0cm} \setlength{\parskip}{0cm}

\bproof
$A$ is classifiable and unital or stably projectionless with continuous scale for the same reasons why the original construction recalled in \S~\ref{ss:CCCC} yields classifiable C*-algebras with these properties (see \cite{Li18, Ell, EV, GLN} for details). Indeed, to see for instance that $A$ is simple, note that with $\varphi_{N,n}$ denoting the composition 
$$
A_n \overset{\varphi_n}{\lori} A_{n+1} \to \dotso \to A_{N-1} \overset{\varphi_{N-1}}{\lori} A_N,
$$
we have for $f \in A_n \subseteq C([0,1],E_n)$ and $t \in [0,1]$ that $(\varphi_{N,n}(f))^q(t) = 0$ for some $q$ only if $f^p(\bar{t}) = 0$ for all $\bar{t} \in \menge{\tfrac{t+k}{2^{N-n}}}{0 \leq k \leq 2^{N-n}-1}$ for all $p$ in the unital case and for all $p \neq \grave{p}$ in the stably projectionless case, and similarly for $p = \grave{p}$ in the stably projectionless case. Hence we see that for all $p$, $\bar{t}$ runs through subsets of $[0,1]$ which become arbitrarily dense in $[0,1]$. This shows simplicity of $A$.
\setlength{\parindent}{0cm} \setlength{\parskip}{0.5cm}

It is clear that $A_n$ has the same K-theory as $\bfdot{A}_n$ and that $\varphi_n$ induces the same map on K-theory as $\bfdot{\varphi}_n$. To see that $\bfdot{\varphi}_n$ induces the same K-theoretic map as $\accir{\varphi}_n$, we construct a homotopy between $\accir{\varphi}_n$ and $\bfdot{\varphi}_n$ as follows: For $s \in [0,1]$, define $\bfdot{\varphi}_s: \: A_n \to \bfdot{A}_{n+1}$ by setting $\bfdot{\varphi}_{s,F}^j: \: A_n \to \bfdot{A}_{n+1} \to F_{n+1}^j$ and $\bfdot{\varphi}_{s,C}: \: A_n \to \bfdot{A}_{n+1} \to C([0,1],E_{n+1})$ as  
\begin{align*}
\bfdot{\varphi}_{s,F}^{j_0^p}(f,a) &\defeq 
\rukl{
\begin{smallmatrix}
f^p(s \cdot \half) & 0 \\
0 & \bar{\psi}^{j_0^p}(a)
\end{smallmatrix}
}, \
\bfdot{\varphi}_{s,F}^{j_1^p}(f,a) \defeq 
\rukl{
\begin{smallmatrix}
f^p(1 - s \cdot \half) & 0 \\
0 & \bar{\psi}^{j_1^p}(a)
\end{smallmatrix}
}, \
\bfdot{\varphi}_{s,F}^j(f,a) \defeq \accir{\varphi}_F^j(f,a) \ {\rm for} \ j \notin \gekl{j_0^p, j_1^p};\\
\bfdot{\varphi}_{s,C} &= \sum_{q,p} \Big( \ve^+_{q,p} \cdot \bfdot{\varphi}_{s,C} \cdot \ve^+_{q,p} + \ve_+^{q,p} \cdot \bfdot{\varphi}_{s,C} \cdot \ve_+^{q,p} + \ve^-_{q,p} \cdot \bfdot{\varphi}_{s,C} \cdot \ve^-_{q,p} + \ve_-^{q,p} \cdot \bfdot{\varphi}_{s,C} \cdot \ve_-^{q,p}\\
&+ \ \underline{\ve}^{q,p} \cdot \bfdot{\varphi}_{s,C} \cdot \underline{\ve}^{q,p} + \bar{\ve}_{q,p} \cdot \bfdot{\varphi}_{s,C} \cdot \bar{\ve}_{q,p} \Big) + \ve_{\rm const} \cdot \bfdot{\varphi}_{s,C} \cdot \ve_{\rm const};\\
\end{align*}
\vspace*{-1.25cm}
\begin{align*}
\ve^+_{q,p} \cdot \bfdot{\varphi}_{s,C}(f,a) \cdot \ve^+_{q,p} &\defeq e^+_{q,p} \otimes f^p \circ (s \cdot \half + (1 - s \cdot \half) \cdot \id),\\ 
\ve_+^{q,p} \cdot \bfdot{\varphi}_{s,C}(f,a) \cdot \ve_+^{q,p} &\defeq e_+^{q,p} \otimes f^p \circ ((1 - s \cdot \half) \cdot \id),\\
\ve^-_{q,p} \cdot \bfdot{\varphi}_{s,C}(f,a) \cdot \ve^-_{q,p} &\defeq e^-_{q,p} \otimes f^p \circ (1 - (1 - s \cdot \half) \cdot \id),\\
\ve_-^{q,p} \cdot \bfdot{\varphi}_{s,C}(f,a) \cdot \ve_-^{q,p} &\defeq e_-^{q,p} \otimes f^p \circ (1 - s \cdot \half - (1 - s \cdot \half) \cdot \id);\\
\underline{\ve}^{q,p} \cdot \bfdot{\varphi}_{s,C}(f,a) \cdot \underline{\ve}^{q,p}
&\defeq \underline{e}^{q,p} \otimes f^p (s \cdot \half) + e_{\scriptscriptstyle (\diagdown)}^{q,p} \otimes f^p \circ (s \cdot \half - s \cdot \half \cdot \id) + \ e_{\scriptscriptstyle (\diagup)}^{q,p} \otimes f^p \circ (s \cdot \half \cdot \id),\\
\bar{\ve}_{q,p} \cdot \bfdot{\varphi}_{s,C}(f,a) \cdot \bar{\ve}_{q,p}
&\defeq \bar{e}_{q,p} \otimes f^p (1 - s \cdot \half) + e^{\scriptscriptstyle (\diagup)}_{q,p} \otimes f^p \circ (1 - s \cdot \half + s \cdot \half \cdot \id) + e^{\scriptscriptstyle (\diagdown)}_{q,p} \otimes f^p \circ (1 - s \cdot \half \cdot \id);\\
\ve_{\rm const} \cdot \bfdot{\varphi}_{s,C} \cdot \ve_{\rm const} &\defeq \ve_{\rm const} \cdot \accir{\varphi}_C \cdot \ve_{\rm const}.
\end{align*}
Then $s \ma \bfdot{\varphi}_s$ is a continuous path connecting $\accir{\varphi}_n$ with $\bfdot{\varphi}_n$. Hence $\accir{\varphi}_n$ and $\bfdot{\varphi}_n$ induce the same map on K-theory.
\setlength{\parindent}{0.5cm} \setlength{\parskip}{0cm}

The same argument as for modification (conn) (see Lemma~\ref{lem:conn:Ell}) shows that our modification (path) yields a C*-algebra $A$ with the desired trace simplex and prescribed pairing between $K_0$ and traces. 
\setlength{\parindent}{0cm} \setlength{\parskip}{0.5cm}

Finally, the connecting maps $\varphi_n$ are of the same form as in \cite[\S~4]{Li18}, and hence admit groupoid models as in \cite[\S~6]{Li18}. Hence $B$ is indeed a C*-diagonal of $A$ by the same argument as in \cite[\S~5 -- 7]{Li18}.  
\eproof
\setlength{\parindent}{0cm} \setlength{\parskip}{0.5cm}

\subsection{Groupoid models for building blocks and connecting maps}
\label{ss:GPDModels}

Before we establish the path-lifting property for our connecting maps, let us first develop a groupoid model for them. Suppose that modification (path) gives us the inductive system 
$$
  A_1 \overset{\varphi_1}{\lori} A_2 \overset{\varphi_2}{\lori} A_3 \overset{\varphi_3}{\lori} \dotso,
$$
with $A_n = \menge{(f,a) \in C([0,1],E_n) \oplus F_n}{f(\mfr) = \beta_{n, \mfr}(a) \ {\rm for} \ \mfr = 0,1}$ for finite-dimensional algebras $E_n$ and $F_n$ as in \S~\ref{ss:(path)} (we use the same notation as in \S~\ref{ss:CCCC}). To describe the connecting map $\varphi \defeq \varphi_n$, we describe
$$
\varphi_C^q: \: A_n \overset{\varphi}{\lori} A_{n+1} \to C([0,1],E_{n+1}) \onto C([0,1],E_{n+1}^q) \quad \text{and} \quad 
\varphi_F^j: \: A_n \overset{\varphi}{\lori} A_{n+1} \to F_{n+1} \onto F_{n+1}^j.
$$
For $q \neq \mfq$, $\varphi_C^q$ is given by the following composition:
\begin{align}
\label{e:varphiC}
(f,a) \ma & \ \big(1_{m^+(q,p)} \otimes f^p \circ \lambda^+, 1_{m_+(q,p)} \otimes f^p \circ \lambda_+, 1_{m^-(q,p)} \otimes f^p \circ \lambda^-, 1_{m_-(q,p)} \otimes f^p \circ \lambda_- )_p,\\
& \ (1_{\underline{m}(q,p)} \otimes f^p(\half))_p, (1_{\overline{m}(q,p)} \otimes f^p(\half))_p , (1_{m^{q,i}} \otimes a^i)_i \big) \nonumber \\
\in & \ C \Big( [0,1], \big(\bigoplus_p (M_{m^+(q,p)} \oplus M_{m_+(q,p)} \oplus M_{m^-(q,p)} \oplus M_{m_-(q,p)}) \otimes E_n^p \big) \nonumber \\
& \ \oplus \big(\bigoplus_p (M_{\underline{m}(q,p)} \oplus M_{\overline{m}(q,p)}) \otimes E_n^p \big)
\oplus \big( \bigoplus_i M_{m^{q,i}} \otimes F_n^i \big) \Big) \nonumber \\
\tailarr & \ C([0,1], E_{n+1}^q). \nonumber
\end{align}
Here $\lambda^+ = \half + \half \cdot \id$, $\lambda_+ = \half \cdot \id$, $\lambda^- = 1 - \half \cdot \id$ and $\lambda_- = \half - \half \cdot \id$. The last arrow is induced by an embedding $\big(\bigoplus_p (M_{m^+(q,p)} \oplus M_{m_+(q,p)} \oplus M_{m^-(q,p)} \oplus M_{m_-(q,p)}) \otimes E_n^p \big) \oplus \big(\bigoplus_p (M_{\underline{m}(q,p)} \oplus M_{\overline{m}(q,p)}) \otimes E_n^p \big) \oplus \big( \bigoplus_i M_{m^{q,i}} \otimes F_n^i \big) \tailarr E_{n+1}^q$ of multiplicity $1$ sending diagonal matrices to diagonal matrices as in \eqref{e:beta=}. Note that $m^+(q,p) = \rk e^+_{q,p} + \rk e^{\scriptscriptstyle (\diagup)}_{q,p}$, $m_+(q,p) = \rk e_+^{q,p} + \rk e_{\scriptscriptstyle (\diagup)}^{q,p}$, $m^-(q,p) = \rk e^-_{q,p} + \rk e^{\scriptscriptstyle (\diagdown)}_{q,p}$ and $m_-(q,p) = \rk e_-^{q,p} + \rk e_{\scriptscriptstyle (\diagdown)}^{q,p}$. By \eqref{e:e>1_unital} -- \eqref{e:e>1_spl2}, we have 
\begin{align}
\label{e:m>1}
& m^+(q,p), m_+(q,p), m^-(q,p), m_-(q,p) \geq 1 \quad \forall q,p && \text{in the unital case};\\
& m^+(q,p), m_+(q,p), m^-(q,p), m_-(q,p) \geq 1 \quad \forall q \neq \grave{q}, p \neq \grave{p}, && \nonumber\\
& m^+(\grave{q},\grave{p}) \geq 1 \quad \text{and} \ m_+(\grave{q},\grave{p}) \ \text{or} \ m_-(\grave{q},\grave{p}) \geq 1 && \text{in the stably projectionless case}.\nonumber
\end{align}
$\varphi_C^\mfq$ is of a similar form, but has an additional component given by $\varphi_F(f,a)$ going into $C([0,1],F_{n+1}) \subseteq C([0,1],E_{n+1}^\mfq)$ (see the second step of modification (path)).

$\varphi_F^j$ is given by the following composition:
\begin{equation}
\label{e:varphiFj_NEW}
(f,a) \ma
\bfa
(f^p(\half), (1_{m(j,i)} \otimes a^i)_i) \in E_n^p \oplus \bigoplus_i M_{m(j,i)} \otimes F_n^i \tailarr F_{n+1}^j & \text{if} \ j = j_\bullet^p,\\
(1_{m(j,i)} \otimes a^i)_i \in \bigoplus_i M_{m(j,i)} \otimes F_n^i \tailarr F_{n+1}^j & \text{if} \ j \notin \gekl{j_0^p, j_1^p}.
\efa
\end{equation}

Recall that $\beta_{n,\mfr} = (\beta^{p,i}_{n,\mfr})_{p,i}$ and that $\beta^{p,i}_{n,\mfr}$ is a composition of the form 
\begin{equation}
\label{e:betanpi}
  F_n^i \overset{1 \otimes \id_{F_n^i}}{\lori} 1_{m_\mfr(p,i)} \otimes F_n^i \subseteq M_{m_\mfr(p,i)} \otimes F_n^i \tailarr E_n^p.
\end{equation}
The groupoid morphism $\bm{b}_{n,\mfr}$ inducing $\beta_{n,\mfr}$ is given on $\cE_{n,\mfr}^p$, the intersection of the domain $\cE_{n,\mfr}$ of $\bm{b}_{n,\mfr}$ with $\cE_n^p$, by
\begin{equation}
\label{e:bmbnp}
  \bm{b}_{n,\mfr}^p: \: \cE_{n,\mfr}^p \cong \coprod_i \cM_\mfr(p,i) \times \cF_n^i \to \coprod_i \cF_n^i = \cF_n,
\end{equation}
where $\cE_n$, $\cE_n^p$, $\cF_n$ and $\cF_n^i$ are groupoid models for $E_n$, $E_n^p$, $F_n$ and $F_n^i$. Now a groupoid model for $(A_n,B_n)$ is given by $G_n \defeq \big( ([0,1] \times_\bullet \cE_n) \amalg \cF_n \big) / {}_\sim$, where $[0,1] \times_\bullet \cE_n \defeq \menge{(t,\gamma) \in [0,1] \times \cE_n}{\gamma \in \cE_{n,t} \ \text{if} \ t = 0,1}$, and $\sim$ is the equivalence relation on $([0,1] \times_\bullet \cE_n) \amalg \cF_n$ generated by $(\mfr,\gamma) \sim \bm{b}_{n,\mfr}(\gamma)$ for all $\mfr = 0,1$, $\gamma \in \cE_{n,\mfr}$. For details, we refer to \cite[\S~6.1]{Li18}. Note that we have $G_n = [[0,1] \times_\bullet \cE_n]$ just as in \S~\ref{s:nccw}, i.e., the extra copy of $\cF_n$ is not needed; it is just convenient to describe the groupoid model $\bm{p}_n$ for $\varphi_n$.

Let us now describe a groupoid model $\bm{p} \defeq \bm{p}_n$ for the connecting map $\varphi_n$ (see \cite[\S~6.2]{Li18} for details). Let $H_n$ be the subgroupoid of $G_{n+1}$ given by $H_n \defeq \big( ([0,1] \times_\bullet \cE_{n+1}[\bm{p}]) \amalg \cF_{n+1}[\bm{p}] \big) / {}_{\sim}$, where, with $\lambda_\mu \defeq \lambda^+$ if $\mu \in \cM^+(q,p)$, $\lambda_\mu \defeq \lambda_+$ if $\mu \in \cM_+(q,p)$, $\lambda_\mu \defeq \lambda^-$ if $\mu \in \cM^-(q,p)$, $\lambda_\mu \defeq \lambda_-$ if $\mu \in \cM_-(q,p)$, $\lambda_\mu \equiv \half$ if $\mu \in \underline{\cM}(q,p) \amalg \overline{\cM}(q,p)$, $\cE_{n+1}[\bm{p}] = \coprod_q \cE_{n+1}^q[\bm{p}]$, and we have identifications
\begin{align}
\label{e:Eqp=}
\cE_{n+1}^q[\bm{p}] \cong \ & {} \big( \coprod_p (\cM^+(q,p) \amalg \cM_+(q,p) \amalg \cM^-(q,p) \amalg \cM_-(q,p)) \times \cE_n^p \big)\\
&\amalg \big( \coprod_p (\underline{\cM}(q,p) \amalg \overline{\cM}(q,p)) \times \cE_n^p \big) 
\amalg \big( \coprod_i \cM^{q,i} \times \cF_n^i \big) \qquad \text{if} \ q \neq \mfq, \nonumber
\end{align}
and similarly for $\cE_{n+1}^\mfq[\bm{p}]$, but with an additional copy of $\cF_{n+1}[\bm{p}]$, and for $\mfr = 0,1$, we have with respect to \eqref{e:Eqp=}:
\begin{eqnarray*}
\cE_{n+1,\mfr}^q[\bm{p}] &=& \big\{ (\mu,\gamma) \in \cE_{n+1}^q[\bm{p}]: \: \gamma \in \cE_{n,\lambda_\mu(\mfr)}^p \ \text{if} \ \mu \in \cM^+(q,p) \amalg \cM_+(q,p) \amalg \cM^-(q,p) \amalg \cM_-(q,p) \big\},\\
\cE_{n+1,\mfr}[\bm{p}] &=& \coprod_q \cE_{n+1,\mfr}^q[\bm{p}]; \
[0,1] \times_\bullet \cE_{n+1}[\bm{p}] 
\defeq \menge{(t,(\mu,\gamma)) \in [0,1] \times \cE_{n+1}[\bm{p}]}{(\mu,\gamma) \in \cE_{n+1,t}[\bm{p}] \ \text{if} \ t \in \gekl{0,1}}.
\end{eqnarray*}
Now $\bm{p}$ is given by $\bm{p}[t,(\mu,\gamma)] = [\lambda_\mu(t),\gamma]$ for $\gamma \in \cE_n^p$ and $\bm{p}(\mu,\gamma) = \gamma$ for $\gamma \in \cF_n^i$. Moreover there are identifications
\begin{equation}
\label{e:F=EMForMF}
\cF_{n+1}^j[\bm{p}] \cong \cE_n^p \amalg \big( \coprod_i \cM(j,i) \times \cF_n^i \big) \quad \text{if} \ j = j_\bullet^p; \qquad
\cF_{n+1}^j[\bm{p}] \cong \coprod_i \cM(j,i) \times \cF_n^i \quad \text{if} \ j \notin \gekl{j_0^p, j_1^p},
\end{equation}
such that $\bm{p}(\mu,\gamma) = \gamma$ for $(\mu,\gamma) \in \cM(j,i) \times \cF_n^i$, and
\begin{equation}
\label{e:bpgamma_FFinE}
\bm{p}(\gamma) = [\half,\gamma] \quad \forall \ \gamma \in \cE_n^p \subseteq \cF_{n+1}^{j_\bullet^p}[\bm{p}] \qquad \text{and} \ \bm{p}[t,\gamma] = \bm{p}(\gamma) \quad \forall \ \gamma \in \cF_{n+1}[\bm{p}] \subseteq \cE_{n+1}^\mfq[\bm{p}], \, t \in [0,1].
\end{equation}
We will often work with the identifications \eqref{e:Eqp=} and \eqref{e:F=EMForMF} without explicitly mentioning them.

That $\varphi_n(f,a)$ satisfies the defining boundary condition for $A_{n+1}$ for all $(f,a) \in A_n$ translates to the following compatibility conditions for $\bm{b}_\bullet$ and $\bm{p}$: We have a commutative diagram
\begin{equation}
\label{e:EEFF}
\begin{tikzcd}
  \cE_{n+1,\mfr}^q[\bm{p}] \arrow[d, "\bm{b}_\mfr"'] \arrow[r, "\subseteq"] & \cE_{n+1,\mfr}^q \arrow[d, "\bm{b}_\mfr"] 
  \\
  \cF_{n+1}[\bm{p}] \arrow[r, "\subseteq"] & \cF_{n+1}
\end{tikzcd}
\end{equation}
For every $\mu \in \cM^+(q,p) \amalg \cM_+(q,p) \amalg \cM^-(q,p) \amalg \cM_-(q,p)$, $\mfr, \mfs \in \gekl{0,1}$ with $\lambda_\mu(\mfr) = \mfs$, the restriction of \eqref{e:EEFF} to $\gekl{\mu} \times \cE_{n,\mfs}^p \subseteq \cE_{n+1,\mfr}^q[\bm{p}]$ fits into the following commutative diagram
\begin{equation}
\label{e:BiggestCD}
\begin{tikzcd}
\gekl{\mfs} \times \cE_{n,\mfs}^p \arrow[d, "\cong"'] & \arrow[l] \gekl{\mfr} \times \gekl{\mu} \times \cE_{n,\mfs}^p \arrow[d, "\cong"] \arrow[r, "\subseteq"] & \gekl{\mfr} \times \cE_{n+1,\mfr}^q \arrow[d, "\cong"]\\
\cE_{n,\mfs}^p \arrow[d, "\cong"'] \arrow[dd, "\bm{b}_{n,\mfs}"', bend right=90] & \arrow[l, "{\gamma \mapsfrom (\mu,\gamma)}"'] \gekl{\mu} \times \cE_{n,\mfs}^p \arrow[dd, "\bm{b}_\mfr"] \arrow[r, "\subseteq"] & \cE_{n+1,\mfr}^q \arrow[d, "\cong"] \arrow[dd, "\bm{b}_\mfr", bend left=90]\\
\coprod_i \cM_\mfs(p,i) \times \cF_n^i \arrow[d] & & \coprod_j \cM_\mfr(q,j) \times \cF_{n+1}^j \arrow[d]\\
\coprod_i \cF_n^i & \arrow[l, "\bm{p}"'] \coprod_j ( \coprod_i \cM(j,i) \times \cF_n^i ) \arrow[r, "\subseteq"] & \coprod_j \cF_{n+1}^j
\end{tikzcd}
\end{equation}
where we identify $\gekl{\mu} \times \cE_{n,\mfs}^p$ and $\coprod_j ( \coprod_i \cM(j,i) \times \cF_n^i )$ with subsets of $\cE_{n+1,\mfr}^q[\bm{p}]$ and $\coprod_j \cF_{n+1}^j[\bm{p}]$ via \eqref{e:Eqp=} and \eqref{e:F=EMForMF}, and the lower vertical arrows on the left and right are given by the canonical projection maps.

Moreover, for all $q, p$, we have $\bm{b}_\mfr(\mu,\gamma) = \gamma \in \cE_n^p \subseteq \cF_{n+1}^{j_\mfr^p}[\bm{p}]$ for all $\mu \in \cM^+(q,p) \amalg \cM_+(q,p)$ and $\bm{b}_\mfr(\mu,\gamma) = \gamma \in \cE_n^p \subseteq \cF_{n+1}^{j_{\mfr^*}^p}[\bm{p}]$ for all $\mu \in \cM^-(q,p) \amalg \cM_-(q,p)$, where $\mfr \in \gekl{0,1}$ satisfies $\lambda_\mu(\mfr) = \halb$, and $\mfr^* = 1 - \mfr$, $\bm{b}_\mfr(\mu,\gamma) = \gamma \in \cE_n^p \subseteq \cF_{n+1}^{j_0^p}[\bm{p}]$ for all $\mu \in \underline{\cM}(q,p)$ and $\mfr = 0,1$, and $\bm{b}_\mfr(\mu,\gamma) = \gamma \in \cE_n^p \subseteq \cF_{n+1}^{j_1^p}[\bm{p}]$ for all $\mu \in\overline{\cM}(q,p)$ and $\mfr = 0,1$. On $\cF_{n+1} \subseteq \cE_{n+1}^\mfq$, $\bm{b}_0$ is given by $\id$ and $\bm{b}_1$ is of a similar form as in \S~\ref{ss:BBCDiagPathConn} and in addition sends $\cE_n^p \subseteq \cF_{n+1}^{j_0^p}$ identically onto $\cE_n^p \subseteq \cF_{n+1}^{j_1^p}$ and $\cE_n^p \subseteq \cF_{n+1}^{j_1^p}$ identically onto $\cE_n^p \subseteq \cF_{n+1}^{j_0^p}$ for all $p \neq \grave{p}$.

Finally, the restriction of \eqref{e:EEFF} to $\coprod_i \cM^{q,i} \times \cF_n^i \subseteq \cE_{n+1,\mfr}^q[\bm{p}]$ fits into the following commutative diagram
\begin{equation*}
\begin{tikzcd}
 & \arrow[dl, "\bm{p}"'] \coprod_i \cM^{q,i} \times \cF_n^i \arrow[d, "\bm{b}_\mfr"] \arrow[r, "\subseteq"] & \cE_{n+1,\mfr}^q \arrow[d, "\bm{b}_\mfr"]\\
\coprod_i \cF_n^i & \arrow[l, "\bm{p}"] \coprod_j ( \coprod_i \cM(j,i) \times \cF_n^i ) \arrow[r, "\subseteq"] & \coprod_j \cF_{n+1}^j
\end{tikzcd}
\end{equation*}

\subsection{The path-lifting property for connecting maps}
\label{ss:path-lift}

We now establish a path-lifting property for $\bm{p} = \bm{p}_n$.
\bprop
\label{prop:path}
Suppose that $\xi_n: \: [0,1] \to G_n$ is a continuous path with the following properties: 
\setlength{\parindent}{0cm} \setlength{\parskip}{0cm}

\begin{enumerate}
\item[(P1)] There exist $0 = \mft_0 < \mft_1 < \dotso < \mft_D < \mft_{D+1} = 1$, $D \geq 0$, such that for all $0 \leq d \leq D$ and $I = [\mft_d,\mft_{d+1}]$, there exist $\gamma_{n,I} \in \cE_n$ and a continuous, monotonous function $\omega_{n,I}: \: I \to [0,1]$ with stop values at $\omega_{n,I}(I) \cap \Zz[\halb]$, i.e., such that, for all $t \in I$, $\xi_n(t) = [\omega_{n,I}(t),\gamma_{n,I}]$. 
\item[(P2)] There exist $d$ and $t \in I = [\mft_d,\mft_{d+1}]$ such that $\omega_{n,I}(t) \in \gekl{0,\halb,1}$ is a stop value of $\omega_{n,I}$. 
\end{enumerate}
Let $\xi_{n+1}^0, \xi_{n+1}^1 \in H_n$ satisfy $\bm{p}(\xi_{n+1}^\mfr) = \xi_n(\mfr)$ for $\mfr = 0,1$. Then there exists a continuous path $\xi_{n+1}: \: [0,1] \to H_n$ with properties (P1) and (P2) such that $\xi_{n+1}(\mfr) = \xi_{n+1}^\mfr$ for $\mfr = 0,1$ and $\bm{p} \circ \xi_{n+1} = \xi_n$.
\setlength{\parindent}{0cm} \setlength{\parskip}{0.5cm}

\textbf{Variation:} Suppose that $\xi_n$ has properties (P1) and (P2), with the following exception:
\setlength{\parindent}{0cm} \setlength{\parskip}{0cm}

\begin{enumerate}
\item[(P3a)] $\omega_{n,[\mft_0,\mft_1]}(0) \in \gekl{0,1}$ is not a stop value for $\omega_{n,[\mft_0,\mft_1]}$, 
\item[(P3b)] there exist $w_{n+1}^0 \in \gekl{0,1}, \mu_{n+1}^0, \gamma_n^0$ such that $\xi_{n+1}^0 = [w_{n+1}^0,(\mu_{n+1}^0,\gamma_n^0)]$ and $\omega_{n,[\mft_d,\mft_{d+1}]}(t) \in \img(\lambda_{\mu_{n+1}^0})$ for all $t \in [0,\mft] \cap [\mft_d,\mft_{d+1}]$ if $\mft_d < \mft$, where $\mft \defeq \min \menge{t>0}{t \in [\mft_d,\mft_{d+1}], \, \omega_{n,[\mft_d,\mft_{d+1}]}(t) \in \gekl{0,\halb,1}}$.
\end{enumerate}
Then we can arrange that $\xi_{n+1}$ has (P3a). We allow for a similar variation for $\omega_{n,[\mft_D,\mft_{D+1}]}(1)$ instead of $\omega_{n,[\mft_0,\mft_1]}(0)$.
\eprop
\setlength{\parindent}{0cm} \setlength{\parskip}{0cm}

Here $w$ is called a stop value of $\omega_{n,I}$ if $\omega_{n,I}$ takes the constant value $w$ on some closed subinterval of $I$ with positive length (see for instance \cite{FR}).
\bproof
By assumption, there are $0 = r_0 \leq t_0 < r_1 < t_1 < \dotso < r_c < t_c < r_{c+1} \leq t_{c+1} = 1$, $c \geq 0$, such that for every interval $I$ of the form $[t_b,r_{b+1}]$, we have $\xi_n(t) = [\omega_{n,I}(t),\gamma_{n,I}]$ for all $t \in I$ for some $\gamma_{n,I} \in G_n$ and $\omega_{n,I} \equiv 0, \halb$ or $1$, and every interval of the form $[r_b,t_b]$ of positive length splits into finitely many subintervals $I$ for which there are $\gamma_{n,I} \in \cE_n$ and continuous maps $\omega_{n,I}: \: I \to [0,1]$ as in (P1) and (P2) such that $\xi_n(t) = [\omega_{n,I}(t),\gamma_{n,I}]$ for all $t \in I$. Moreover, for $1 \leq b \leq c$ and $I \subseteq [r_b,t_b]$ as above, $\omega_{n,I}$ does not take the values $0$, $\halb$ or $1$ on $(r_b,t_b)$. Set $\xi_{n+1}[0] \defeq \xi_{n+1}^0$, $\xi_{n+1}[c+1] \defeq \xi_{n+1}^1$, and write $\xi_{n+1}[0] = [w[0],(\mu[0],\gamma[0])]$, $\xi_{n+1}[c+1] = [w[c+1],(\mu[c+1],\gamma[c+1])]$. If $w[0] \in (0,1)$, we arrange $t_0 > 0$ by replacing $t_0$ by $\half \cdot r_1$ if necessary. If $w[c+1] \in (0,1)$, we arrange $r_{c+1} < 1$ by replacing $r_{c+1}$ by $\half \cdot (t_c + 1)$ if necessary. As a result, if $w[0] \in (0,1)$, we must have $t_0 > 0$, and either $\omega_{n,I}$ does not take the values $0$, $\halb$ or $1$ on $I \cap [r_0,t_0)$ for each $I \subseteq [r_0,t_0]$ as above, or $\omega_{n,[r_0,t_0]}$ is constant with value $0$, $\halb$ or $1$. If $w[0] \in \gekl{0,1}$, then we must have $t_0 = 0$. For the variation, we must have $t_0 > 0$, and for each $I \subseteq [r_0,t_0]$ as above, $\omega_{n,I}$ does not take the values $0$, $\halb$ or $1$ on $I \cap (r_0,t_0)$, $\omega_{n,I}(0) \in \gekl{0,1}$ is not a stop value, and $w[0] \in \gekl{0,1}$. A similar statement holds for $w[c+1]$.
\setlength{\parindent}{0cm} \setlength{\parskip}{0.5cm}

For $1 \leq b \leq c$, take $s_b \in (r_b,t_b)$ and choose $\xi_{n+1}[b] = [w[b],(\mu[b],\gamma[b])] \in H_n$ such that $\bm{p}(\xi_{n+1}[b]) = \xi_n(s_b)$. Such $\xi_{n+1}[b]$ exist because of \eqref{e:m>1}. Define $s_0 \defeq 0$ and $s_{c+1} \defeq 1$. Now let $0 \leq b \leq c+1$. Suppose that $\xi_n(s_b)$ is of the form $[w,\gamma]$ with $w \notin \gekl{0,1}$, which is always the case if $1 \leq b \leq c$. Let $I \subseteq [r_b,t_b]$ be as above. Define $\gamma_{n+1,I} \defeq \gamma[b]$. If $\lambda_{\mu[b]} = \lambda^+$, define $\omega_{n+1,I} \defeq -1 + 2 \cdot \omega_{n,I}$, if $\lambda_{\mu[b]} = \lambda_+$, define $\omega_{n+1,I} \defeq 2 \cdot \omega_{n,I}$, if $\lambda_{\mu[b]} = \lambda^-$, define $\omega_{n+1,I} \defeq 2 - 2 \cdot \omega_{n,I}$, and if $\lambda_{\mu[b]} = \lambda_-$, define $\omega_{n+1,I} \defeq 1 - 2 \cdot \omega_{n,I}$. If $b=0$, $I = [r_0,t_0]$, $t_0 > 0$, i.e., $w[0] \in (0,1)$, and if $\omega_{n,I} \equiv 0, \halb$ or $1$, set $\gamma_{n+1,I} \defeq (\mu[0],\gamma[0])$ and let $\omega_{n+1,I}$ be a continuous path as in (P1) with $\omega_{n+1,I}(0) = w[0]$, $\omega_{n+1,I}(0) = 1$ ((P2) is then automatic). Such a path exists by \cite[Lemma~2.10]{FR}. Define $\gamma_{n+1,I}$ and $\omega_{n+1,I}$ similarly for $b=c+1$, $I = [r_{c+1},t_{c+1}]$, $r_{c+1} < 1$ and $\omega_{n,I} \equiv 0, \halb$ or $1$. For the variation, note that (P3b) implies that we can define $\gamma_{n+1,I}$ and $\omega_{n+1,I}$ for $I \subseteq [r_0,t_0]$ and $I \subseteq [r_{c+1},t_{c+1}]$ as above in the same way as for $I \subseteq [r_b,t_b]$, where $\xi_n(s_b)$ is of the form $[w,\gamma]$ with $w \notin \gekl{0,1}$. Now set $\xi_{n+1}(t) \defeq [\omega_{n+1,I}(t),\gamma_{n+1,I}]$ for all $t \in I$.

Next, consider $I = [t_b,r_{b+1}]$ for $0 \leq b \leq c$. First assume that $\omega_{n,I} \equiv \half$. Let $\gamma_{n,I} = \gamma$. Set $w \defeq \omega_{n+1,[r_b,t_b]}(t_b)$, $\bar{w} \defeq \omega_{n+1,[s_{b+1},t_{b+1}]}(s_{b+1})$ and let $\gamma_{n+1,[s_b,t_b]}(t_b) = (\mu,\gamma)$, $\gamma_{n+1,[r_{b+1},t_{b+1}]}(s_{b+1}) = (\bar{\mu},\gamma)$. Note that $w, \bar{w} \in \gekl{0,1}$. If $[w,(\mu,\gamma)] = [\bar{w},(\bar{\mu},\gamma)]$, then set $\omega_{n+1,I} \equiv w$ and $\gamma_{n+1,I} \defeq (\mu,\gamma)$. If $[w,(\mu,\gamma)] \neq [\bar{w},(\bar{\mu},\gamma)]$, then $\bm{b}_w(\mu,\gamma) = \gamma^{j_0^p} \in \cE_n^p \subseteq \cF_{n+1}^{j_0^p}$ and $\bm{b}_{\bar{w}}(\bar{\mu},\gamma) = \gamma^{j_1^p} \in \cE_n^p \subseteq \cF_{n+1}^{j_1^p}$ for $p \neq \grave{p}$ (or with $j_0^p$ and $j_1^p$ swapped), where $\gamma^{j_\mfr^p}$ denotes the copy of $\gamma$ in $\cF_{n+1}^{j_\mfr^p}$. Let $\omega_{n+1,I}$ be a continuous path as in (P1) such that $\omega_{n+1,I}(t_b) = 0$, $\omega_{n+1,I}(s_{b+1}) = 1$ (such a path exists by \cite[Lemma~2.10]{FR}, and (P2) is automatic), and define $\gamma_{n+1,I} \defeq \gamma^{j_0^p}$. Set $\xi_{n+1}(t) \defeq [\omega_{n+1,I}(t), \gamma_{n+1,I}]$ for all $t \in I$. Then by \eqref{e:bpgamma_FFinE}, we have $\bm{p}(\xi_{n+1}(t)) = [\halb,\gamma] = \xi_n(t)$ for all $t \in I$, as well as $\xi_{n+1}(t_b) = [0,\gamma^{j_0^p}] = [w,(\mu,\gamma)]$ since $\bm{b}_0(\gamma^{j_0^p}) = \gamma^{j_0^p}$ and $\xi_{n+1}(r_{b+1}) = [1,\gamma^{j_0^p}] = [0,\gamma^{j_1^p}] = [\bar{w},(\bar{\mu},\gamma)]$ as $\bm{b}_1(\gamma^{j_0^p}) = \gamma^{j_1^p} = \bm{b}_0(\gamma^{j_1^p})$.
\setlength{\parindent}{0.5cm} \setlength{\parskip}{0cm}

Now assume that $\omega_{n,I} \equiv 0$. Set $w \defeq \omega_{n+1,[r_b,t_b]}(t_b)$, $\bar{w} \defeq \omega_{n+1,[r_{b+1},t_{b+1}]}(r_{b+1})$ and let $\gamma_{n+1,[r_b,t_b]}(t_b) = (\mu,\gamma)$, $\gamma_{n+1,[r_{b+1},t_{b+1}]}(r_{b+1}) = (\bar{\mu},\bar{\gamma})$. We have $\xi_n(t) = [0,\gamma_{n,I}] = \bm{p}[w,(\mu,\gamma)] = \bm{p}[\bar{w},(\bar{\mu},\bar{\gamma})]$ for all $t \in I$. Note that $w, \bar{w} \in \gekl{0,1}$. Now $[w,(\mu,\gamma)] = [0,\bm{b}_w(\mu,\gamma)]$ and $[\bar{w},(\bar{\mu},\bar{\gamma})] = [0, \bm{b}_{\bar{w}}(\bar{\mu},\bar{\gamma})]$, where we view $\bm{b}_w(\mu,\gamma)$ and $\bm{b}_{\bar{w}}(\bar{\mu},\bar{\gamma})$ as elements of $\cE^{\mfq}_{\rm conn}$ (the analogue of $\cY^{\mfq}_{\rm conn}$ in \S~\ref{ss:BBCDiagPathConn}). We have $\bm{p}[0,\bm{b}_w(\mu,\gamma)] = \bm{p}[w,(\mu,\gamma)] = \bm{p}[\bar{w},(\bar{\mu},\bar{\gamma})] = \bm{p}[0, \bm{b}_{\bar{w}}(\bar{\mu},\bar{\gamma})]$, so that by the analogue of Remark~\ref{rem:Yconnconn} for $\cE^{\mfq}_{\rm conn}$ instead of $\cY^{\mfq}_{\rm conn}$, after possibly splitting $I$ into finitely many subintervals, we can find $\gamma_{n+1,I}$ and $\omega_{n+1,I}$ as in (P1) and (P2) such that, if we define $\xi_{n+1}(t) \defeq [\omega_{n+1,I}(t),\gamma_{n+1,I}]$ for all $t \in I$, then we have $\bm{p}(\xi_{n+1}(t)) = \bm{p}[0,\bm{b}_w(\mu,\gamma)] = \bm{p}[w,(\mu,\gamma)] = \xi_n(t)$ for all $t \in I$, and $\xi_{n+1}(t_b) = [0,\bm{b}_w(\mu,\gamma)] = [w,(\mu,\gamma)]$, $\xi_{n+1}(r_{b+1}) = [0, \bm{b}_{\bar{w}}(\bar{\mu},\bar{\gamma})] = [\bar{w},(\bar{\mu},\bar{\gamma})]$.

The case $\omega_{n,I} \equiv 1$ is similar.
\eproof
\setlength{\parindent}{0cm} \setlength{\parskip}{0.5cm}

\section{Constructing C*-diagonals with Menger manifold spectra}
\label{s:CDiagMenger}

Suppose that modification (path) produces the C*-algebra $A = \ilim_n \gekl{A_n,\varphi_n}$ with prescribed Elliott invariant $\cE$ as in \S~\ref{ss:(path)} and the C*-diagonal $B = \ilim_n \gekl{B_n,\varphi_n}$ of $A$ as in Lemma~\ref{lem:path:Ell}. In the following, we write $X_n \defeq \Spec B_n$, $X \defeq \Spec B$. Note that $X$ is metrizable, Hausdorff and compact (in the unital case) or locally compact (in the stably projectionless case), $X \cong \plim_n \gekl{X_n,\bm{p}_n}$ and $\dim X \leq 1$ (see \cite{Li18}). Our goal now is to determine $X$ further. Let $\bm{p}_{n,\infty}: \: X \to X_n$ be the map given by the inverse limit structure of $X$ and $\bm{p}_{n,N}: \: X_{N+1} \to X_n$ the composition $\bm{p}_{n,N} \defeq \bm{p}_n \circ \dotso \circ \bm{p}_N$. Moreover, the groupoid model $G_n$ for $A_n$ in \S~\ref{ss:GPDModels} yields descriptions $X_n \cong \big( ([0,1] \times_\bullet \cY_n) \amalg \cX_n \big) / {}_\sim$, where $\cY_n = \cE_n^{(0)}$, $\cX_n = \cF_n^{(0)}$, and with $\cY_{n,\mfr} \defeq \cY_n \cap \cE_{n,\mfr}$, $[0,1] \times_\bullet \cY_n \defeq \menge{(t,y) \in [0,1] \times \cY_n}{y \in \cY_{n,t} \ \text{if} \ t = 0,1}$, and $\sim$ is the equivalence relation on $([0,1] \times_\bullet \cY_n) \amalg \cX_n$ generated by $(\mfr,y) \sim \bm{b}_{n,\mfr}(y)$ for all $\mfr = 0,1$, $y \in \cY_{n,\mfr}$.

\bprop
\label{prop:pathconn}
The C*-diagonal $B$ has path-connected spectrum $X$.
\eprop
\setlength{\parindent}{0cm} \setlength{\parskip}{0cm}

\bproof
Let $\eta = (\eta_n)_n$, $\zeta = (\zeta_n)_n$ be two points in $X$. The induction start in the proof of Proposition~\ref{prop:conn} shows that there exists a continuous path $\xi_1: \: [0,1] \to X_1$ with $\xi_1(0) = \eta_1$, $\xi_1(1) = \zeta_1$. Using \cite[Lemma~2.10]{FR}, it is straightforward to see that $\xi_1$ can be chosen with property (P1) and (P2). Applying Proposition~\ref{prop:path} recursively, we obtain continuous paths $\xi_n: \: [0,1] \to X_n$ with $\xi_n(0) = \eta_n$, $\xi_n(1) = \zeta_n$ and $\bm{p}_n \circ \xi_{n+1} = \xi_n$. Hence $\xi(t) \defeq (\xi_n(t))_n$ defines a continuous path $[0,1] \to X \cong \plim_n \gekl{X_n, \bm{p}_n}$ with $\xi(0) = \eta$ and $\xi(1) = \zeta$.
\eproof
\setlength{\parindent}{0cm} \setlength{\parskip}{0.5cm}

\bprop
\label{prop:lpc}
The spectrum $X$ of $B$ is locally path-connected.
\eprop
\setlength{\parindent}{0cm} \setlength{\parskip}{0cm}

\bproof
Consider a point $\bm{c} = ([w_n,y_n])_n \in X$ with $w_n \in [0,1]$, $y_n \in \cY_n$, and an open set $V$ of $X$ with $\bm{c} \in V$. First suppose that there is $\underline{n}$ such that $w_n \notin \gekl{0,1}$ for all $n \geq \underline{n}$. Then there exists $n \geq \underline{n}$, an open interval $I_n \subseteq (0,1)$, $\alpha, e \in \Zz_{\geq 0}$ such that $\tfrac{\alpha}{2^e} < w_n < \tfrac{\alpha + 1}{2^e}$, $[\tfrac{\alpha}{2^e}, \tfrac{\alpha + 1}{2^e}] \subseteq I_n$ and $\bm{p}_{n,\infty}^{-1}[I_n \times \gekl{y_n}] \subseteq V$. It is straightforward to see that if there exists an open interval $I_{n+m} \subseteq (0,1)$ of length at least $\frac{1}{2^{e-m}}$ with $w_{n+m} \in I_{n+m}$ and $\halb \notin I_{n+m}$ such that $[I_{n+m} \times \gekl{y_{n+m}}] \subseteq \bm{p}_{n,n+m}^{-1}[I_n \times \gekl{y_n}]$, then there exists an open interval $I_{n+m+1} \subseteq (0,1)$ of length at least $\frac{1}{2^{e-m-1}}$ with $w_{n+m+1} \in I_{n+m+1}$ and such that $[I_{n+m+1} \times \gekl{y_{n+m+1}}] \subseteq \bm{p}_{n,n+m+1}^{-1}[I_n \times \gekl{y_n}]$. Thus there exists $m \leq e-1$ and an open interval $I_{n+m} \subseteq (0,1)$ with $w_{n+m}, \halb \in I_{n+m}$ such that $[I_{n+m} \times \gekl{y_{n+m}}] \subseteq \bm{p}_{n,n+m}^{-1}[I_n \times \gekl{y_n}]$. Hence $U \defeq \bm{p}_{n+m,\infty}^{-1}[I_{n+m} \times \gekl{y_{n+m}}]$ satisfies $\bm{c} \in U \subseteq V$. Set $\bar{n} \defeq n+m$, so that $U = \bm{p}_{\bar{n},\infty}^{-1}[I_{\bar{n}} \times \gekl{y_{\bar{n}}}]$. Now assume that for all $\underline{n}$ there is $n \geq \underline{n}$ such that $w_n \in \gekl{0,1}$. Then, since $V$ is open, there exists $\bar{n} \geq \underline{n}$ and, for every $(\mfr,y) \sim (w_{\bar{n}},y_{\bar{n}})$, half-open intervals $I(\mfr,y)$ containing $\mfr$ such that $U \defeq \bigcup_{(\mfr,y) \sim (w_{\bar{n}},y_{\bar{n}})} \bm{p}_{\bar{n},\infty}^{-1}[I(\mfr,y) \times \gekl{y}] \subseteq V$.
\setlength{\parindent}{0cm} \setlength{\parskip}{0.5cm}

We claim that in both cases above, $U$ is path-connected. Let $\eta = (\eta_n), \zeta = (\zeta_n) \in U$. We construct a path $\xi_{\bar{n}}: \: [0,1] \to X_{\bar{n}}$ with (P1) and (P2) such that $\xi_{\bar{n}}(0) = \eta_{\bar{n}}$ and $\xi_{\bar{n}}(1) = \zeta_{\bar{n}}$. Let us treat the first case ($w_{\bar{n}} \notin \gekl{0,1}$). We have $\eta_{\bar{n}} = [w_{\bar{n}}^0,y_{\bar{n}}]$, $\zeta_{\bar{n}} = [w_{\bar{n}}^1,y_{\bar{n}}]$. Define $\xi_{\bar{n}}$ as in (P1), with $D=1$, $\mft_1 = \halb$, for $I = [\mft_0,\mft_1] = [0,\halb]$, $\gamma_{\bar{n},I} \defeq y_{\bar{n}}$, $\omega_{\bar{n},I}: \: [0,\halb] \to [0,1]$ as in (P1) with $\omega_{\bar{n},I}(0) = w_{\bar{n}}^0$, $\omega_{\bar{n},I}(\halb) = \halb$, and for $I = [\mft_1,\mft_2] = [\halb,1]$, $\gamma_{\bar{n},I} \defeq y_{\bar{n}}$, $\omega_{\bar{n},I}: \: [\halb,1] \to [0,1]$ as in (P1) with $\omega_{\bar{n},I}(\halb) = \halb$, $\omega_{\bar{n},I}(1) = w_{\bar{n}}^1$ (such paths exist by \cite[Lemma~2.10]{FR} and have (P2)). In the second case ($w_{\bar{n}} \in \gekl{0,1}$), let $\eta_{\bar{n}} = [w_{\bar{n}}^0,y_{\bar{n}}^0]$, $\zeta_{\bar{n}} = [w_{\bar{n}}^1,y_{\bar{n}}^1]$. There must exist $\mfr^0, \mfr^1 \in \gekl{0,1}$ with $(\mfr^0,y_{\bar{n}}^0) \sim (w_{\bar{n}},y_{\bar{n}}) \sim (\mfr^1,y_{\bar{n}}^1)$. Define $\xi_{\bar{n}}$ as in (P1), with $D=1$, $\mft_1 = \halb$, for $I = [\mft_0,\mft_1] = [0,\halb]$, $\gamma_{\bar{n},I} \defeq y_{\bar{n}}^0$, $\omega_{\bar{n},I}: \: [0,\halb] \to [0,1]$ as in (P1) with $\omega_{\bar{n},I}(0) = w_{\bar{n}}^0$, $\omega_{\bar{n},I}(\halb) = \mfr^0$, and for $I = [\mft_1,\mft_2] = [\halb,1]$, $\gamma_{\bar{n},I} \defeq y_{\bar{n}}^1$, $\omega_{\bar{n},I}: \: [\halb,1] \to [0,1]$ as in (P1) with $\omega_{\bar{n},I}(\halb) = \mfr^1$, $\omega_{\bar{n},I}(1) = w_{\bar{n}}^1$ (such paths exist by \cite[Lemma~2.10]{FR} and have (P2)). Now apply Proposition~\ref{prop:path} to obtain paths $\xi_n$ and thus a path $\xi$ connecting $\eta$ and $\zeta$ as for Proposition~\ref{prop:pathconn}.
\eproof
\setlength{\parindent}{0cm} \setlength{\parskip}{0.5cm}

\bcor
\label{cor:Peano}
$X$ is a Peano continuum in the unital case and a generalized Peano continuum in the stably projectionless case (see for instance \cite[Chapter~I, \S~9]{BQ} for the definition of a generalized Peano continuum).
\ecor

Our next goal is to show that we can always arrange $X$ to have no local cut points. In the following, we keep the same notations as in \S~\ref{ss:GPDModels}. First, we observe that in modification (path), because multiplicities in the original C*-algebra models in \cite{Ell, EV, GLN} can be chosen bigger than a fixed constant, by conjugating by suitable permutation matrices, we can always arrange the following conditions for all $n$:
\setlength{\parindent}{0cm} \setlength{\parskip}{0cm}

\begin{enumerate}
\item[(nlc$_1$)] For all $p$ and $m = m^+, m_+, m^-$ or $m_-$, we either have $\sum_q m(q,p) = 0$ or $\sum_q m(q,p) \geq 2$, and $\sum_q (\underline{m}(q,p) + \overline{m}(q,p)) \geq 2$; and for all $i$, we have $\sum_q m^{q,i} \geq 2$;
\item[(nlc$_2$)] For all $p$, $m = m^+, m_+, m^-$ or $m_-$, $\lambda = \lambda^+, \lambda_+, \lambda^-$ or $\lambda_-$ correspondingly, $\mfr, \mfs \in \gekl{0,1}$ with $\lambda(\mfr) = \mfs$, rank-one projections $d \in DM_{m(p,i)}$ and $\delta$ the image of $d \otimes 1_{F_n^i}$ under $\tailarr$ in the description \eqref{e:betanpi} of $\beta_{n,\mfs}^p$, and rank-one projections $\mfd \in DM_{m(q,p)}$, there exists a rank-one projection $\mfd' \in DM_{m(q',p)}$ orthogonal to $\mfd$, and orthogonal projections $\mff, \mff' \in DF_{n+1}$ such that, if $\Delta$ is the image of $\mfd \otimes \delta$ under $\tailarr$ in the description \eqref{e:varphiC} of $\varphi_C$, then
$$
\mfd \otimes (\delta \cdot \beta_{n,\mfs}^p(a) \cdot \delta) = \mfd \otimes (\delta \cdot f^p(\mfs) \cdot \delta) \tailarr \Delta \cdot \varphi_C^q(f,a)(\mfr) \cdot \Delta = \Delta \cdot \beta_\mfr(\mff \cdot \varphi_F(f,a) \cdot \mff) \cdot \Delta
$$
for all $(f,a) \in A_n$ with respect to the description \eqref{e:varphiC} of $\varphi_C$, and similarly for $\mfd'$ and $\mff'$.
\end{enumerate} 
On the groupoid level, with the same notation as in \S~\ref{ss:GPDModels}, (nlc)$_2$ means that for all $\gamma \in \cM(p,i) \times \cF_n^i \into \cE_{n,\mfs}^p$, $\mu \in \cM(q,p)$, where $\cM = \cM^+, \cM_+, \cM^-$ or $\cM_-$, $\lambda_\mu(\mfr) = \mfs$, there exists $\nu \in \cM(q',p)$ such that $\bm{b}_\mfr(\mu,\gamma) \neq \bm{b}_\mfr(\nu,\gamma)$, i.e., $[\mfr,(\mu,\gamma)] \neq [\mfr,(\nu,\gamma)]$, and $\lambda_\mu = \lambda_\nu$.

\bprop
\label{prop:nlc}
If we arrange (nlc$_1$), (nlc$_2$) in modification (path), then we obtain a C*-diagonal $B$ whose spectrum $X$ has no local cut points (i.e., for all $\bm{c} \in X$ and open connected sets $V \subseteq X$ containing $\bm{c}$, $V \setminus \gekl{\bm{c}}$ is still connected).
\eprop
\bproof
Let $\bm{c}$, $V$ and $U$ be as in Proposition~\ref{prop:lpc}. It suffices to show that $U \setminus \gekl{\bm{c}}$ is path-connected. Let $\eta, \zeta \in U \setminus \gekl{\bm{c}}$ and $\xi$ a path in $U$ connecting $\eta$ and $\zeta$ as in the proof of Proposition~\ref{prop:lpc}. If $\xi$ hits $\bm{c}$, our goal is to modify $\xi$ to obtain a path in $U$ from $\eta$ to $\zeta$ which avoids $\bm{c}$. First of all, we may assume that $\xi$ hits $\bm{c}$ only once, i.e., there exists $\check{t} \in [0,1]$ such that $\xi(\check{t}) = \bm{c}$ and $\xi(t) \neq \bm{c}$ for all $t \in [0,1] \setminus \gekl{\check{t}}$. Otherwise we could define $t^{\min} \defeq \min \menge{t \in [0,1]}{\xi(t) = \bm{c}}$, $t^{\max} \defeq \max \menge{t \in [0,1]}{\xi(t) = \bm{c}}$, and concatenate $\xi \vert_{[0,t^{\min}]}$ with $\xi \vert_{[t^{\max},1]}$ (and re-parametrize to get a map defined on $[0,1]$). Let $\xi_n$ be as in the proof of Proposition~\ref{prop:lpc}, obtained from Proposition~\ref{prop:path}, and let $\gamma_{n,I}$ and $\omega_{n,I}$ be as in (P1) for $\xi_n$. Choose $n$ such that, with $\bm{c}_n = [w_n,y_n]$, we have $[w_n,y_n] \neq \eta_n = [w_n^0,y_n^0]$ and $[w_n,y_n] \neq \zeta_n = [w_n^1,y_n^1]$. 
\setlength{\parindent}{0cm} \setlength{\parskip}{0.5cm}

If $w_n^0 \notin \gekl{0,1}$, then either $y_n \neq y_n^0$, in which case $\omega_n([0,\check{t}]) \defeq \bigcup_I \omega_{n,I}(I \cap [0,\check{t}])$ must contain either $[0,w_n^0]$ or $[w_n^0,1]$, or $w_n \neq w_n^0$, in which case $\omega_n([0,\check{t}])$ must contain the interval between $w_n$ and $w_n^0$. If $w_n^0 \in \gekl{0,1}$ and $w_n \notin \gekl{0,1}$, then $\omega_n([0,\check{t}])$ must contain the interval between $w_n$ and $w_n^0$ (or $1 - w_n^0$). If $w_n^0, w_n \in \gekl{0,1}$, then since $[w_n,y_n] \neq [w_n^0,y_n^0]$, we must have $\omega_n([0,\check{t}]) = [0,1]$. We conclude that in any case, $\omega_n([0,\check{t}]) \cap \Zz[\halb] \neq \emptyset$ and $\omega_n([0,\check{t}]) \cap (\Zz[\halb])^c \neq \emptyset$. Similarly, with $\omega_n([\check{t},1]) \defeq \bigcup_I \omega_{n,I}(I \cap [\check{t},1])$, $\omega_n([\check{t},1]) \cap \Zz[\halb] \neq \emptyset$ and $\omega_n([\check{t},1]) \cap (\Zz[\halb])^c \neq \emptyset$. By increasing $n$, we can arrange that $0, \halb \ \text{or} \ 1 \in \omega_n([0,\check{t}])$ and $\omega_n([0,\check{t}]) \cap (\Zz[\halb])^c \neq \emptyset$, as well as $0, \halb \ \text{or} \ 1 \in \omega_n([\check{t},1])$ and $\omega_n([\check{t},1]) \cap (\Zz[\halb])^c \neq \emptyset$. If we now let $0 = r_0 \leq t_0 < r_1 < t_1 < \dotso < r_c < t_c < r_{c+1} \leq t_{c+1} = 1$ be as in the proof of Proposition~\ref{prop:path}, then we must have $\check{t} \in (r_1,t_c)$.

First assume that $w_n \neq 0, \halb, 1$ and $w_{n+1} \neq 0, 1$. Then $\check{t} \in I \subseteq [r_b,t_b]$, and $\check{t}$ must lie in the interior of $I$. Let $s_b$ and $\xi_{n+1}[b] = [w[b],(\mu[b],y[b])]$ be as in the proof of Proposition~\ref{prop:path}. By condition (nlc$_1$), we can find $\bar{\mu}[b] \neq \mu[b]$ with $\lambda_{\bar{\mu}[b]} = \lambda_{\mu[b]}$, so that we can replace $\xi_{n+1}[b]$ by $[w[b],(\bar{\mu}[b],y[b])]$ since we still have $\bm{p}_n[w[b],(\bar{\mu}[b],y[b])] = \xi_n(s_b)$. Now let $\gamma_{n+1,I} \defeq (\bar{\mu}[b],y[b])$ and follow the recipe in the proof of Proposition~\ref{prop:path} to get $\omega_{n+1,I}$. Recursive application of Proposition~\ref{prop:path} gives us the desired path, which will not hit $\bm{c}$ in $(r_b,t_b)$ by construction, on $[t_{b-1},r_b]$ and $[t_b,r_{b+1}]$, we have $\omega_{n,I} \equiv 0, \halb \ \text{or} \ 1$, so that we will not hit $\bm{c}$ there, either, and on the rest of $[0,1]$, we keep our path $\xi$ and hence will not hit $\bm{c}$ there, either.

Secondly, assume that $w_n = 0, \halb, \ \text{or} \ 1$ and $w_{n+1} \neq 0, 1$. Then $\check{t} \in I = [r_b,t_b]$, and $\check{t}$ must lie in the interior of $I$. Let $\gamma_{n+1,I} = (\mu,y)$. We must have $\lambda_{\mu} \equiv w_n$. By condition (nlc$_1$), there exists $\bar{\mu} \neq \mu$ with $\lambda_{\bar{\mu}} = \lambda_\mu$. Now replace $\gamma_{n+1,I}$ by $(\bar{\mu},y)$ and follow the recipe in the proof of Proposition~\ref{prop:path} to get $\omega_{n+1,I}$. Recursive application of Proposition~\ref{prop:path} gives us the desired path, which will not hit $\bm{c}$ in $I$ by construction, and on $[0,1] \setminus I$, we keep our path $\xi$ and hence will not hit $\bm{c}$ there, either. 

Thirdly, assume $w_n \in \gekl{0, \halb, 1}$ and $w_{n+1} \in \gekl{0,1}$. By increasing $n$ if necessary, we may assume $w_n \in \gekl{0,1}$. We have $\check{t} \in I \defeq [t_b,r_{b+1}]$. If $\check{t} \in (t_b,r_{b+1})$, then choose a different path between $[w,(\mu,y)]$ and $[\bar{w},(\bar{\mu},\bar{y})]$ (here we are using the same notation as in the proof of Proposition~\ref{prop:path}). There are always two such paths only overlapping at their end points (see Remark~\ref{rem:Yconnconn} and the proof of Proposition~\ref{prop:conn} it refers to). Complete the construction of $\xi$ on $[t_b,r_{b+1}]$ using Proposition~\ref{prop:path} repeatedly. Keep $\xi$ on $[0,1] \setminus (t_b,r_{b+1})$. This yields the desired path which does not hit $\bm{c}$ anymore. Now assume that $\check{t} = t_b < r_{b+1}$. By condition (nlc$_2$), there exists $\nu$ with $\lambda_\nu = \lambda_\mu$ such that $[w,(\nu,y)] \neq [w,(\mu,y)]$. Construct a path as in the proof of Proposition~\ref{prop:path} connecting $[w,(\nu,y)]$ and $[\bar{w},(\bar{\mu},\bar{y})]$ not hitting $\bm{c}_{n+1}=[w_{n+1},y_{n+1}]$. This is possible because there are always two paths connecting these points and only overlapping at their end points (see Remark~\ref{rem:Yconnconn} and the proof of Proposition~\ref{prop:conn} it refers to). On the interval $\dot{I}$ before $t_b$, re-define $\xi_{n+1}$ by setting $\gamma_{n+1,\dot{I}} \defeq (\nu,y)$ and following the recipe in the proof of Proposition~\ref{prop:path} for $\omega_{n+1,\dot{I}}$. On the interval $\ddot{I}$ before $\dot{I}$, either $\omega_{n,\ddot{I}} \equiv \halb$, in which case simply following the recipe in the proof of Proposition~\ref{prop:path} for $\omega_{n+1,\ddot{I}}$ will make sure that we do not hit $\bm{c}$, or $\omega_{n,\ddot{I}} \equiv w$, in which case we construct $\xi_{n+1}$ on $\ddot{I}$ avoiding $\bm{c}_{n+1}$ using as before that Remark~\ref{rem:Yconnconn} (and the proof of Proposition~\ref{prop:conn} it refers to) always provides two paths we can choose from to connect end points as required. Complete the construction of $\xi$ on $I$, $\dot{I}$ and $\ddot{I}$ using Proposition~\ref{prop:path} repeatedly, and keep $\xi$ on the remaining part of $[0,1]$. This yields the desired path not hitting $\bm{c}$. Finally, suppose that $\check{t} = t_b = r_{b+1}$ (this can happen due to our re-parametrisation). By condition (nlc$_2$), there exist $\nu$, $\bar{\nu}$ with $\lambda_\nu = \lambda_\mu$, $\lambda_{\bar{\nu}} = \lambda_{\bar{\mu}}$ and $[w,(\nu,y)] \neq [w,(\mu,y)]$, $[\bar{w},(\bar{\nu},\bar{y})] \neq [\bar{w},(\bar{\mu},\bar{y})]$. As in the previous cases, we can construct a path connecting $[w,(\nu,y)]$ and $[\bar{w},(\bar{\nu},\bar{y})]$ not hitting $[w_{n+1},y_{n+1}]$. Complete the construction of $\xi$ on $I$ using Proposition~\ref{prop:path} repeatedly. On the two intervals before and after $t_b = r_{b+1}$, construct $\xi$ as in the previous case so that we do not hit $\bm{c}$ there. On the remaining part of $[0,1]$, keep the path $\xi$, which does not hit $\bm{c}$ by assumption. This yields the desired path not hitting $\bm{c}$.
\eproof
\setlength{\parindent}{0cm} \setlength{\parskip}{0.5cm}

Next, we show that we can always arrange $X$ so that no non-empty open subset of $X$ is planar. For this purpose, we observe that in modification (path), with the same notations as in \S~\ref{ss:GPDModels}, for the same reasons why we can always arrange (nlc$_1$) and (nlc$_2$), we can always arrange the following conditions for all $n$:
\setlength{\parindent}{0cm} \setlength{\parskip}{0cm}

\begin{enumerate}
\item[(nop$_1$)] For all $p$, $\sum_q m^+(q,p) \geq 1$ and $\sum_q m_+(q,p) \geq 1$, and $\sum_q \underline{m}(q,p) \geq 9$ or $\sum_q \overline{m}(q,p) \geq 9$;
\item[(nop$_2$)] The analogue of (nlc$_2$), implying on the groupoid level that for all $\gamma \in \cE_{n,0}^p$, there exist $\nu_+^i \in \cM_+(q^i,p)$, $i = 1, 2, 3$, such that $\bm{b}_0(\nu_+^i,\gamma)$ are pairwise distinct, i.e., $[0,(\nu_+^i,\gamma)]$ are pairwise distinct, and for all $\gamma \in \cE_{n,1}^p$, there exist $\nu^+_i \in \cM^+(q_i,p)$, $i = 1, 2, 3$, with $\bm{b}_1(\nu^+_i,\gamma)$ pairwise distinct, i.e., $[1,(\nu^+_i,\gamma)]$ pairwise distinct.
\end{enumerate} 

\bprop
\label{prop:nop}
If we arrange (nop$_1$) and (nop$_2$) in modification (path), then we obtain a C*-diagonal $B$ whose spectrum $X$ has the property that no non-empty open subset of $X$ is planar.
\eprop
\bproof
Let $\emptyset \neq V \subseteq X$ be open. A similar argument as in the beginning of the proof of Proposition~\ref{prop:lpc} shows that there exists $n$ and an open subset $U_n$ of $X_n$ such that $[\halb,y] \in U_n$ for some $y \in \cE_n$ and $\bm{p}_{n,\infty}^{-1}(U_n) \subseteq V$. By condition (nop$_1$), there exist $\mu^{ij}$, $1 \leq i, j \leq 3$ with $\lambda_{\mu^{ij}} \equiv \halb$ and $[0,(\mu^{ij},y)] = [1,(\mu^{kl},y)]$ for all $i, j, k, l$. Let $\xi_{n+1}^{ij}(t) \defeq [\omega_{n+1}(t),(\mu^{ij},y)]$, where $\omega_{n+1}: \: [0,1] \to [0,1]$ is as in (P1), with $\omega_{n+1}(0) = 0$, $\omega_{n+1}(1) = 1$ and $\omega_{n+1}(0)$, $\omega_{n+1}(1)$ are not stop values ($\omega_{n+1}$ exists by \cite[Lemma~2.10]{FR} and automatically has (P2)). Set $y_{n+1} \defeq (\mu^{11},y)$. By condition (nop$_2$), we can find $\nu_+^i$, $i = 1, 2, 3$, such that $\lambda_{\nu_+^i} = \halb \cdot \id$ and, with $y_{n+2,i}^0 \defeq (\nu_+^i,y_{n+1})$, we have that $[0, y_{n+2,i}^0]$ are pairwise distinct for $i = 1, 2, 3$. Similarly, by condition (nop$_2$), we can find $\nu^+_j$, $j = 1, 2, 3$, such that $\lambda_{\nu^+_i} = \halb + \halb \cdot \id$ and, with $y_{n+2,j}^1 \defeq (\nu^+_j,y_{n+1})$, we have that $[1, y_{n+2,j}^1]$ are pairwise distinct for $j = 1, 2, 3$. By the variation of Proposition~\ref{prop:path}, we can find a path $\xi_{n+2}^{ij}$ satisfying (P1), (P2) and (P3a) such that $\xi_{n+2}^{ij}(0) = [0, y_{n+2,i}^0]$, $\xi_{n+2}^{ij}(1) = [1, y_{n+2,j}^1]$ and $\bm{p}_{n+1} \circ \xi_{n+2}^{ij} = \xi_{n+1}^{ij}$. Now define recursively $y_{N,i}^0$ and $y_{N,j}^1$ for all $N \geq n+2$ by setting $y_{N+1,i}^0 \defeq (\mu_+,y_{N,i}^0)$ for some $\mu_+$ with $\lambda_{\mu_+} = \halb \cdot \id$ and $y_{N+1,j}^1 \defeq (\mu^+,y_{N,j}^1)$ for some $\mu^+$ with $\lambda_{\mu^+} = \halb + \halb \cdot \id$. We can find such $\mu_+$ and $\mu^+$ by the first part of condition (nop$_1$). By the variation of Proposition~\ref{prop:path}, we can find paths $\xi_N^{ij}$ satisfying (P1), (P2) and (P3a) such that $\xi_N^{ij}(0) = [0, y_{N,i}^0]$, $\xi_N^{ij}(1) = [1, y_{N,j}^1]$ and $\bm{p}_{N-1} \circ \xi_N^{ij} = \xi_{N-1}^{ij}$ for all $N \geq n+2$. This gives rise to paths $\xi^{ij}$, $1 \leq i, j \leq 3$, with $\xi^{ij}(t) \defeq (\xi_N^{ij}(t))_N$. As $\bm{p}_{n,\infty}(\xi^{ij}(t)) = \bm{p}_n(\xi_{n+1}^{ij}(t)) \in U_n$ for all $t \in [0,1]$, we must have $\img(\xi^{ij}) \subseteq \bm{p}_{n,\infty}^{-1}(U_n) \subseteq V$ for all $i$, $j$. Now define $v_i^0 \defeq ([0,y_{N,i}^0])_N$ and $v_j^1 \defeq ([1,y_{N,j}^1])_N$. By construction, we have $\img(\xi^{i,j}) \cap \img(\xi^{k,l}) = \{ v^0_i \}$ if $i = k$ and $j \neq l$, $\img(\xi^{i,j}) \cap \img(\xi^{k,l}) = \{ v^1_j \}$ if $i \neq k$ and $j = l$, and $\img(\xi^{i,j}) \cap \img(\xi^{k,l}) = \emptyset$ if $i \neq k$ and $j \neq l$. As $\img(\xi^{i,j})$ is a compact, connected, locally connected metric space, it is arcwise connected (see for instance \cite[\S~31]{Wil}). Hence we can find arcs $\xi^{(i \to j)}$ such that $\xi^{(i \to j)}(0) = v^0_i$, $\xi^{(i \to j)}(1) = v^1_j$, and $\img(\xi^{(i \to j)}) \subseteq \img(\xi^{i,j})$ for $i, j \in \gekl{1,2,3}$. Then we still have that $\img(\xi^{(i \to j)}) \cap \img(\xi^{(k \to l)}) = \{ v^0_i \}$ if $i = k$ and $j \neq l$, $\img(\xi^{(i \to j)}) \cap \img(\xi^{(k \to l)}) = \{ v^1_j \}$ if $i \neq k$ and $j = l$, and $\img(\xi^{(i \to j)}) \cap \img(\xi^{(k \to l)}) = \emptyset$ if $i \neq k$ and $j \neq l$. Now let $K_{3,3}$ be the bipartite graph consisting of vertices $e^i(0)$, $e^j(1)$, $1 \leq i,j \leq 3$, and edges $e^{(i \to j)}$, $1 \leq i, j \leq 3$, connecting $e^i(0)$ to $e^j(1)$ for all $i,j$, such that we have $e^{(i \to j)} \cap e^{(k \to l)} = \gekl{e^i(0)}$ if $i = k$ and $j \neq l$, $e^{(i \to j)} \cap e^{(k \to l)} = \gekl{e^j(1)}$ if $i \neq k$ and $j = l$, and $e^{(i \to j)} \cap e^{(k \to l)} = \emptyset$ if $i \neq k$ and $j \neq l$. By construction, we obtain a continuous map $K_{3,3} \to V$ which is a homeomorphism onto its image by sending $e^i(0)$ to $v^0_i$, $e^j(1)$ to $v^1_j$ and $e^{(i \to j)}$ to $\xi^{(i \to j)}$. Since $K_{3,3}$ is not planar by \cite{Kur}, this shows that $V$ is not planar.
\eproof
\setlength{\parindent}{0cm} \setlength{\parskip}{0.5cm}

\bcor
\label{cor:Menger_unital}
Suppose that we are in the unital case and that we arrange (nlc$_1$), (nlc$_2$), (nop$_1$) and (nop$_2$) in modification (path). Then we obtain a C*-diagonal $B$ whose spectrum $X$ is homeomorphic to the Menger curve.
\ecor
\setlength{\parindent}{0cm} \setlength{\parskip}{0cm}

\bproof
Anderson characterized the Menger curve as the (up to homeomorphism) unique Peano continuum with no local cut points and for which no non-empty open subset is planar (see \cite{And58_1,And58_2}). Our result thus follows from Corollary~\ref{cor:Peano} combined with Propositions~\ref{prop:nlc} and \ref{prop:nop}.
\eproof
\setlength{\parindent}{0cm} \setlength{\parskip}{0.5cm}

Our next goal is to identify $X = \Spec B$ in the stably projectionless case. We show that $X \cong \bm{M} \setminus \iota(C)$, where $\iota$ is an embedding of the Cantor space $C$ into the Menger curve $\bm{M}$ such that $\iota(C)$ is a non-locally-separating subset of $\bm{M}$. By \cite{MOT}, the homeomorphism type of $\bm{M} \setminus \iota(C)$ does not depend on the choice of $\iota$, and hence we denote the space by $\bm{M}_{\setminus C} \defeq \bm{M} \setminus \iota(C)$. More precisely, we will show that the Freudenthal compactification $\overline{X}{}^{^F}$ of $X$ is homeomorphic to $\bm{M}$, that the space of Freudenthal ends ${\rm End}_F(X)$ is homeomorphic to $C$, and that ${\rm End}_F(X)$ is a non-locally-separating subset of $\overline{X}{}^{^F}$. It follows that $X$ is homeomorphic to $\bm{M}_{\setminus C}$. We refer the reader to \cite{Fre} and \cite[Chapter~I, \S~9]{BQ} for details about the Freudenthal compactification. We follow the exposition in \cite[Chapter~I, \S~9]{BQ}.

First of all, in the stably projectionless case, we define $\overline{X}_n \defeq \big( ([0,1] \times \cY_n) \amalg \cX_n \big) / {}_\sim$, where we extend the equivalence relation describing $X_n$ (see the beginning of \S~\ref{s:CDiagMenger}) trivially from $([0,1] \times_\bullet \cY_n) \amalg \cX_n$ to $([0,1] \times \cY_n) \amalg \cX_n$. By our arrangement, for all $n$, there exists exactly one index $\grave{p}$ such that $\beta_{n,0}^{\grave{p}}$ is unital and $\beta_{n,1}^{\grave{p}}$ is non-unital, while $\beta_{n,\bullet}^p$ is unital for all other $p \neq \grave{p}$. This means that $\cY_{n,0} = \cY_n$ and $\cY_n \setminus \cY_{n,1} = \cY_n^{\grave{p}} \setminus \cY_{n,1}^{\grave{p}}$. Hence $\overline{X}_n \setminus X_n = \big\{ [1,y_n]: \: y_n \in \cY_n^{\grave{p}} \setminus \cY_{n,1}^{\grave{p}} \big\}$. Let $\bar{\bm{p}}_n: \: \overline{X}_{n+1} \to \overline{X}_n$ be the unique continuous extension of $\bm{p}_n$. Every $y_{n+1} \in \cY_{n+1}^{\grave{q}} \setminus \cY_{n+1,1}^{\grave{q}}$ is of the form $y_{n+1} = (\mu,y_n)$ for some $y_n \in \cY_n^{\grave{p}} \setminus \cY_{n,1}^{\grave{p}}$, $\mu \in \cM^+(\grave{q},\grave{p})$, and we have $\bar{\bm{p}}_n[1,y_{n+1}] = [1,y_n]$. Define $\overline{X} \defeq \plim_n \big\{ \overline{X}_n,\bar{\bm{p}}_n \big\}$. 
\blemma
\label{lem:XXF}
$\id_X$ extends to a homeomorphism $\overline{X} \isom \overline{X}{}^{^F}$.
\elemma
\setlength{\parindent}{0cm} \setlength{\parskip}{0cm}

\bproof
For $y \notin \cY_{n,1}$, define $I_y \defeq [0,\tfrac{1}{3}]$, and for $y \in \cY_{n,1}$, set $I_y \defeq [0,1]$. Define $K_n \defeq \bm{p}_{n,\infty}^{-1}[\bigcup_{y \in \cY_n} I_y \times \gekl{y}]$. Then $K_n$ is compact because $K_n \cong \plim_N \gekl{\bm{p}_{n,N}^{-1}[\bigcup_{y \in \cY_n} I_y \times \gekl{y}], \bm{p}_N}$ and $\bm{p}_{n,N}$ is proper (see \cite[\S~7]{Li18}). Every $y_{n+1} \notin \cY_{n+1,1}$ is of the form $y_{n+1} = (\mu,y_n)$ for some $y_n \notin \cY_{n,1}$ with $\lambda_\mu = \halb + \halb \cdot \id$, so that $\bm{p}_n[t,y_{n+1}] = [\halb + \tfrac{t}{2},y_n] \notin [[0,\tfrac{1}{3}] \times \gekl{y_n}]$ for all $t \in [0,1]$. Hence $\bm{p}_n^{-1}[\bigcup_{y_n \in \cY_n} I_{y_n} \times \gekl{y_n}] \subseteq [\bigcup_{y \in \cY_{n+1,1}} [0,1] \times \gekl{y}] \subseteq {\rm int}([\bigcup_{y \in \cY_{n+1}} I_y \times \gekl{y}])$. Thus $K_n \subseteq {\rm int}(K_{n+1})$ for all $n$. Moreover, $X \setminus K_n = \bm{p}_{n,\infty}^{-1}[\bigcup_{y \notin \cY_{n,1}} (\tfrac{1}{3},1) \times \gekl{y}] = \bigcup_{y \notin \cY_{n,1}} \bm{p}_{n,\infty}^{-1}[(\tfrac{1}{3},1) \times \gekl{y}]$. Using Proposition~\ref{prop:path}, the same argument as for Proposition~\ref{prop:pathconn} shows that $\bm{p}_{n,\infty}^{-1}[(\tfrac{1}{3},1) \times \gekl{y}]$ is path-connected, and we obtain $\menge{[1,y]}{y \notin \cY_{n,1}} \isom \Pi_0(X \setminus K_n), \, [1,y] \ma \bm{p}_{n,\infty}^{-1}[(\tfrac{1}{3},1) \times \gekl{y}]$. This induces a homeomorphism $ \overline{X} \setminus X = \ilim_n \gekl{\menge{[1,y_n]}{y_n \notin \cY_{n,1}}, \bm{p}_n} \isom \plim_n \Pi_0(X \setminus K_n) = {\rm End}_F(X)$ and hence a (set-theoretic) bijection $\overline{X} \isom \overline{X}{}^{^F}$ extending $\id_X$. For this description of ${\rm End}_F(X)$, we are using that $X$ is a generalized Peano continuum (see Corollary~\ref{cor:Peano}). It is now straightforward to see that this bijection is a homeomorphism.
\eproof
\setlength{\parindent}{0cm} \setlength{\parskip}{0.5cm}

To study properties of $\overline{X}$, we need the following observation.
\bremark
\label{rem:path-lift_spl}
In $\overline{X}$, the analogue of Proposition~\ref{prop:path} holds for a path $\xi_n$ with $\xi_n(0) \in \overline{X}_n \setminus X_n$, $\xi_n(1) \in X_n$ and $\xi_{n+1}^0 \in \overline{X}_{n+1} \setminus X_{n+1}$, $\xi_{n+1}^1 \in X_{n+1}$ with $\bm{p}(\xi_{n+1}^\mfr) = \xi_n(\mfr)$ for $\mfr = 0,1$. We also have the analogue of the variation, but we only need (P3a) because (P3b) is automatic in the present situation since we must have $\lambda_{\mu_{n+1}^0} = \half + \half \cdot \id$, and we get the additional statement that if $\xi_n(t) \in X_n \ \forall \ t \in (0,1]$, then $\xi_{n+1}(t) \in X_{n+1} \ \forall \ t \in (0,1]$.
\eremark

\bprop
\label{prop:clX...}
$\overline{X}$ is compact, path-connected and locally path-connected. If we arrange (nlc$_1$) and (nlc$_2$) in modification (path), then $\overline{X}$ has no local cut points. If we arrange (nop$_1$) and (nop$_2$) in modification (path), then no non-empty subset of $\overline{X}$ is planar.
\eprop
\setlength{\parindent}{0cm} \setlength{\parskip}{0cm}

\bproof
Clearly, $\overline{X}$ is compact. To see that $\overline{X}$ is path-connected, consider $\eta, \zeta \in \overline{X}$. If both $\eta$ and $\zeta$ lie in $X$, then Proposition~\ref{prop:pathconn} provides a path connecting them. If $\eta \in \overline{X} \setminus X$ and $\zeta \in X$, we produce a path connecting them as in the proof of Proposition~\ref{prop:pathconn} using the analogue of Proposition~\ref{prop:path} from Remark~\ref{rem:path-lift_spl}. If both $\eta$ and $\zeta$ lie in $\overline{X} \setminus X$, just connect them to some point in $X$ and concatenate the two paths. To see that $\overline{X}$ is locally path-connected, we follow the same strategy as for Proposition~\ref{prop:lpc}. We only need to consider $\bm{c} = ([1,y_n])_n \in \overline{X} \setminus X$. Choose $U$ in the proof of Proposition~\ref{prop:lpc} of the form $U = \bm{p}_{n,\infty}^{-1}[I \times \gekl{y_n}]$, where $I$ is a half-open interval containing $\halb$ and $1$. Then the same proof as for Proposition~\ref{prop:lpc}, using the analogue of Proposition~\ref{prop:path} from Remark~\ref{rem:path-lift_spl}, shows that $U$ is path-connected. To show that $\overline{X}$ has no local cut points if (nlc$_1$) and (nlc$_2$) hold, we again only need to consider $\bm{c} = ([1,y_n])_n \in \overline{X} \setminus X$. Choose $U$ as before and take $\eta, \zeta \in U \setminus \gekl{\bm{c}}$. If both $\eta$ and $\zeta$ lie in $X$, then Proposition~\ref{prop:nlc} yields a path in $U \setminus \gekl{\bm{c}}$ connecting $\eta$ and $\zeta$ because $U \cap X$ is of the form as in Proposition~\ref{prop:nlc}. If $\eta \in \overline{X} \setminus X$ and $\zeta \in X$, then we can construct a path $\xi$ in $U$ with $\xi(0) = \eta$, $\xi(1) = \zeta$ and $\xi(t) \in X$ for all $t \in (0,1]$, using the analogue of the variation in Proposition~\ref{prop:path} from Remark~\ref{rem:path-lift_spl}. Then $\xi(t) \neq \bm{c}$ for all $t \in (0,1]$, and we also have $\xi(0) = \eta \neq \bm{c}$. If both $\eta$ and $\zeta$ lie in $\overline{X} \setminus X$, then pick a point $u \in U \cap X$, connect $\eta$ and $\zeta$ to $u$ in $U \setminus \gekl{\bm{c}}$ as in the previous case, and concatenate the two paths. Finally, to see that no non-empty open subset of $\overline{X}$ is planar if (nlc$_1$) and (nlc$_2$) hold, just observe that every non-empty open subset $V$ of $\overline{X}$ gives rise to a non-empty open subset $V \cap X$ of $X$, and apply Proposition~\ref{prop:nop} to $V \cap X$.
\eproof
\setlength{\parindent}{0cm} \setlength{\parskip}{0.5cm}

The same reasoning as for Corollary~\ref{cor:Menger_unital} yields
\bcor
If we arrange (nlc$_1$), (nlc$_2$), (nop$_1$) and (nop$_2$) in modification (path), then $\overline{X}$ is homeomorphic to the Menger curve.
\ecor

\blemma
\label{lem:Cantor_spl}
If (nlc$_1$) holds, then $\overline{X} \setminus X$ is homeomorphic to the Cantor space.
\elemma
\setlength{\parindent}{0cm} \setlength{\parskip}{0cm}

\bproof
(nlc$_1$) implies that we always have $m^+(\grave{q},\grave{p}) \geq 2$, so that for all $y_n \notin \cY_{n,1}$, $\# \bm{p}_n^{-1}[1,y_n] \geq 2$. Now it is straightforward to see that $ \overline{X} \setminus X = \ilim_n \gekl{\menge{[1,y_n]}{y_n \notin \cY_{n,1}}, \bm{p}_n}$ is homeomorphic to the Cantor space.
\eproof
\setlength{\parindent}{0cm} \setlength{\parskip}{0.5cm}

\bprop
$\overline{X} \setminus X$ is a non-locally-separating subset of $\overline{X}$, i.e., for every connected open subset $V \subseteq \overline{X}$, $V \setminus (\overline{X} \setminus X) = V \cap X$ is connected.
\eprop
\setlength{\parindent}{0cm} \setlength{\parskip}{0cm}

\bproof
$V$ is open and connected, hence locally path-connected by Proposition~\ref{prop:clX...} and thus path-connected. Take $\eta, \zeta \in V \cap X$ and a continuous path $\xi: \: [0,1] \to \overline{X}$ with $\xi(0) = \eta$ and $\xi(1) = \zeta$. It is straightforward to see that we can find $0 = t_0 < t_1 < \dotso < t_l < t_{l+1} = 1$ and for each $0 \leq k \leq l$ an open subset $\bm{U}_k \subseteq V$ as in the proof that $X$ and $\overline{X}$ are locally path-connected (see Propositions~\ref{prop:lpc} and \ref{prop:clX...}) such that $\xi([t_k,t_{k+1}]) \subseteq \bm{U}_k$ for all $0 \leq k \leq l$. Set $\xi[0] \defeq \xi(0)$, $\xi[1] \defeq \xi(1)$, and for $1 \leq k \leq l$, set $\xi[t_k] \defeq \xi(t_k)$ if $\xi(t_k) \in X$ and pick some $\xi[t_k] \in \bm{U}_{k-1} \cap \bm{U}_k \cap X$ otherwise. Since $\bm{U}_k \cap X$ is an open set of the form as in the proof of Proposition~\ref{prop:lpc}, it is path-connected, so that we can find paths connecting $\xi[t_k]$ and $\xi[t_{k+1}]$ in $\bm{U}_k \cap X$ for all $0 \leq k \leq l$. Now concatenate these paths to obtain a path in $V \cap X$ connecting $\eta$ and $\zeta$.
\eproof
\setlength{\parindent}{0cm} \setlength{\parskip}{0.5cm}

All in all, we obtain the following consequence.
\bcor
\label{cor:M-C_spl}
Suppose that we are in the stably projectionless case and that we arrange (nlc$_1$), (nlc$_2$), (nop$_1$) and (nop$_2$) in modification (path). Then we obtain a C*-diagonal $B$ whose spectrum $X$ is homeomorphic to $\bm{M}_{\setminus C}$.
\ecor

\bremark
\label{rem:KMenger}
The K-groups of $C(\bm{M})$ and $C_0(\bm{M}_{\setminus C})$ are given as follows: We have $K_0(C(\bm{M})) = \Zz[1]$, $K_1(C(\bm{M})) \cong \bigoplus_{i=1}^\infty \Zz$ (see for instance \cite[Equation~(32)]{Li18}), and it follows that $K_0(C_0(\bm{M}_{\setminus C})) \cong \gekl{0}$, $K_1(C_0(\bm{M}_{\setminus C})) \cong \bigoplus_{i=1}^\infty \Zz$.
\eremark

\section{Constructing continuum many non-conjugate C*-diagonals with Menger manifold spectra}
\label{s:ManyCDiagMenger}

Let us present two further modifications of our constructions of C*-diagonals in classifiable C*-algebras which will allow us to produce continuum many pairwise non-conjugate C*-diagonals in all our classifiable C*-algebras. First of all, we recall the construction of the groupoid model $\bar{G}$ for the pair $(A,B)$, where $A$ is our classifiable C*-algebra with prescribed Elliot invariant $\cE$ as in \S~\ref{ss:(path)} and $B$ the C*-diagonal of $A$ produced by modification (path). Let $G_n$, $H_n$ and $\bm{p}_n$ be as in \S~\ref{ss:GPDModels}. Following \cite[\S~5]{Li18}, we define $G_{n,0} \defeq G_n$, $G_{n,m+1} \defeq \bm{p}_{n+m}^{-1}(G_{n,m}) \subseteq H_{n+m}$ for all $n$ and $m = 0, 1, \dotsc$, $\bar{G}_n \defeq \plim_m \gekl{G_{n,m}, \bm{p}_{n+m}}$ for all $n$. Moreover, the inclusions $H_n \into G_{n+1}$ induce embeddings with open image $\bm{i}_n: \: \bar{G}_n \into \bar{G}_{n+1}$, allowing us to define $\bar{G} \defeq \ilim \gekl{\bar{G}_n, \bm{i}_n}$. We will identify $\bar{G}_n$ with its image in $\bar{G}$. As explained in \cite[\S~5]{Li18}, $\bar{G}$ is a groupoid model for $(A,B)$ in the sense that we have a canonical isomorphism $A \isom C^*_r(\bar{G})$ sending $B$ to $C_0(\bar{G}^{(0)})$. In the following, we let $\bm{p}_{n+m,\infty}: \: \bar{G}_n \onto G_{n,m}$ be the canonical projection from the inverse limit structure of $\bar{G}_n$, and $\bm{p}_{n+m,n+\bar{m}}: \: G_{n,\bar{m}+1} \onto G_{n,m}$ denotes the composition $\bm{p}_{n+m} \circ \dotso \circ \bm{p}_{n+\bar{m}}$.

\subsection{Constructing closed subgroupoids}

Recall the description of $\beta_{n,\mfr}^{p,i}$ in \eqref{e:betanpi}. We observe that in modification (path), with the same notations as in \S~\ref{ss:GPDModels}, we can always arrange the following condition for all $n$ by adding $\id_{F_{n+1}^j}$ to $\beta_{n+1,\bullet}^q$, enlarging $E_{n+1}^q$ accordingly, and conjugation $\beta_{n+1,\bullet}^q$ by suitable permutation matrices as in modification (conn) without changing the properties or Elliott invariant of the classifiable C*-algebra we construct:
\setlength{\parindent}{0cm} \setlength{\parskip}{0cm}

\begin{enumerate}
\item[(clsg)] For all $q, p$, $m = m^+, m_+, m^-$ or $m_-$ and $\lambda = \lambda^+, \lambda_+, \lambda^-$ or $\lambda_-$ correspondingly, $\mfr, \mfs \in \gekl{0,1}$ with $\lambda(\mfr) = \mfs$, $\mfd \in DM_{m(q,p)}$, if we denote by $\fD(p,i)$ the set of rank-one projections in $DM_{m_\mfs(p,i)}$ and by $\fD(q,j)$ the set of rank-one projections in $DM_{m_\mfr(q,j)}$, then there is an injective map $\coprod_i \fD(p,i) \into \coprod_j \fD(q,j), \, d(p,i) \ma d(q,j)$ and $d \in DM_{m(j,i)}$ attached to each pair $d(p,i)$ and $d(q,j)$ such that, if we denote by $\delta(p,i)$ the image of $d(p,i) \otimes 1_{F_n^i}$ under $\tailarr$ in the description of $\beta_{n,\mfs}^p$ in \eqref{e:betanpi}, by $\Delta$ the image of $\mfd \otimes \delta(p,i)$ under $\tailarr$ in the description \eqref{e:varphiC} of $\varphi_C$, and by $\delta$ the image of $d \otimes 1_{F_n^i}$ under $\tailarr$ in the description of $\varphi_F^j \vert_{F_n^i}$ as in \eqref{e:varphiFj_NEW}, we have that for all $(f,a) \in A_n$
$$
d(p,i) \otimes a \tailarr \delta(p,i) \cdot f^p(\mfs) \cdot \delta(p,i)
$$
under $\tailarr$ in the description of $\beta_{n,\mfs}^p$ in \eqref{e:betanpi},
$$
\mfd \otimes (\delta(p,i) \cdot f^p(\mfs) \cdot \delta(p,i)) \tailarr \Delta \cdot \varphi_C^q(f,a)(\mfr) \cdot \Delta = \Delta \cdot \beta_\mfr(\varphi_F(f,a)) \cdot \Delta
$$
under $\tailarr$ in the description \eqref{e:varphiC} of $\varphi_C$, and
$$
d(q,j) \otimes \delta \cdot \varphi_F^j(a) \cdot \delta \tailarr \Delta \cdot \beta_\mfr(\varphi_F(f,a)) \cdot \Delta
$$
under $\tailarr$ in the description of $\beta_{n+1,\mfr}^q$ in \eqref{e:betanpi}.
\end{enumerate} 
On the groupoid level, with the same notation as in \S~\ref{ss:GPDModels}, (clsg) means that in \eqref{e:BiggestCD}, if we start at $\coprod_i \cM_\mfs(p,i) \times \cF_n^i$, go up to $\cE_{n,\mfs}^p$, follow the horizontal arrows to $\cE_{n+1,\mfr}^q$ and go down to $\coprod_j \cM_\mfr(q,j) \times \cF_{n+1}^j$, we get a map of the form $(\mu(p,i),\gamma_n^i) \ma (\mu(q,j),\ti{\gamma}_n^i)$ such that the assignment $\mu(p,i) \ma \mu(q,j)$ is injective (and $\gamma_n^i \ma \ti{\gamma}_n^i$ corresponds to the composition $a \ma d(j,i) \otimes a \tailarr \delta(j,i) \cdot \varphi_F^j(a) \cdot \delta(j,i)$ as in the description of $\varphi_F^j \vert_{F_n^i}$ in \eqref{e:varphiFj_NEW}).
\setlength{\parindent}{0cm} \setlength{\parskip}{0.5cm}

\blemma
\label{lem:HnClosed}
If we arrange (clsg) in modification (path), then $H_n$ is a closed subgroupoid of $G_{n+1}$.
\elemma
\setlength{\parindent}{0cm} \setlength{\parskip}{0cm}

\bproof
We make use of the descriptions of $H_n$ and $G_{n+1}$ in \S~\ref{ss:GPDModels}. Suppose $[t_k,\mu_k,\gamma_k] \in H_n$ converges in $G_{n+1}$ to $[t,\gamma_{n+1}] \in G_{n+1}$. Our goal is to show that $[t,\gamma_{n+1}]$ lies in $H_n$. As there are only finitely many possibilities for $(\mu_k,\gamma_k)$, we may assume that $(\mu_k,\gamma_k) = (\mu,\gamma)$ is constant (independent of $k$), and thus $\gamma_{n+1} = (\mu,\gamma)$. If $t \in (0,1)$, then we have $[t,\gamma_{n+1}] \in H_n$. Now suppose that $t \in \gekl{0,1}$. If $\mu \in \overline{\cM}(q,p) \amalg \underline{\cM}(q,p)$, $\gamma \in \cE_n^p$, or $\mu \in \cM(q,i)$, $\gamma \in \cF_n^i$, then $(\mu,\gamma) \in \cE_{n+1,t}^q$ and hence $[t,\gamma_{n+1}] \in H_n$. Finally, assume that $\mu \in \cM^+(q,p) \amalg \cM_+(q,p) \amalg \cM^-(q,p) \amalg \cM_-(q,p)$ and $\gamma \in \cE_n^p$. Let $\lambda = \lambda^+, \lambda_+, \lambda^-, \ \text{or} \ \lambda_-$ accordingly. Since $[t,(\mu,\gamma)] \in G_{n+1}$, $s(\mu,\gamma)$ and $r(\mu,\gamma)$ must be mapped to elements in $\coprod_j \cM_t(q,j) \times \cX_{n+1}^j$ with the same $\cM_t(q,j)$-component in \eqref{e:BiggestCD}, which then implies by (clsg) that $s(\gamma)$ and $r(\gamma)$ are mapped to elements in $\coprod_i \cM_{\lambda(t)}(p,i) \times \cX_n^i$ with the same $\cM_{\lambda(t)}(p,i)$-component in \eqref{e:BiggestCD}, which in turn implies that $\gamma \in \cE_{n,\lambda(t)}^p$ and thus $[t,(\mu,\gamma)] \in H_n$.
\eproof
\setlength{\parindent}{0cm} \setlength{\parskip}{0.5cm}

\bcor
\label{cor:barGn-clopen}
If we arrange (clsg) in modification (path), then $\bar{G}_n$ is a clopen subset of $\bar{G}$ for all $n = 1, 2, \dotsc$.
\ecor
\setlength{\parindent}{0cm} \setlength{\parskip}{0cm}

\bproof
An easy induction on $m$ shows that $G_{n,m}$ is an open subset of $G_{n+m}$: $G_{n,1} = H_n$ is open in $G_{n+1}$ by construction (see \cite[\S~6.2]{Li18}), and for the induction step, use the recursive definition of $G_{n,m}$ together with continuity of $\bm{p}_{n+m}$ and the observation that $H_{n+m}$ is open in $G_{n+m+1}$. Hence, for all $n$, $\bar{G}_n = \bm{p}_{n+m,\infty}^{-1}(G_{n,m})$ is an open subset of $\bar{G}_{n+m}$ for all $m = 0, 1, \dotsc$. By definition of the inductive limit topology, this shows that $\bar{G}_n$ is open in $\bar{G}$.
\setlength{\parindent}{0.5cm} \setlength{\parskip}{0cm}

To see that $\bar{G}_n$ is closed in $G$, let $(\bm{g}_k)_k$ be a sequence in $\bar{G}_n$ converging to $\bm{g} \in \bar{G}$. Suppose $\bm{g} \notin \bar{G}_n$. Then let $m \geq 1$ be minimal with $\bm{g} \in \bar{G}_{n+m}$. We have $\bm{g}_k \in \bar{G}_n = \bm{p}_{n+m,\infty}^{-1}(G_{n,m}) \subseteq \bm{p}_{n+m,\infty}^{-1}(H_{n+m-1})$ for all $k$. Since $H_{n+m-1}$ is closed in $G_{n+m}$ by Lemma~\ref{lem:HnClosed}, we must have $\bm{g} \in \bm{p}_{n+m,\infty}^{-1}(H_{n+m-1}) = \bm{p}_{n+m,\infty}^{-1}(G_{n+m-1,1}) = \bar{G}_{n+m-1}$. But this contradicts minimality of $m$. Hence $\bm{g} \in \bar{G}_n$, and thus $\bar{G}_n$ is closed in $\bar{G}$.
\eproof
\setlength{\parindent}{0cm} \setlength{\parskip}{0.5cm}

\subsection{Modification (sccb)}

We now present a further modification (path), which works in a similar way as modification (conn) or the second step in modification (path) and for the same reasons will not change the properties or Elliott invariant of the classifiable C*-algebra we construct. Let us use the same notations as in \S~\ref{ss:GPDModels}. Suppose that the first step in modification (path) produces the first building block $A_1$. For all $d \in DF_1^i$ choose a permutation matrix $w_{i,d} \in F_1^i$ such that $w_{i,d} d w_{i,d}^* = d$ and $w_{i,d} \hat{d} w_{i,d}^* \neq \hat{d}$ for all $\hat{d} \in DF_1^i$ with $\hat{d} \neq d$. Let $w \defeq (w_{i,d})_{i,d}$. Choose an index $\ti{q}$ and replace $E_1^{\ti{q}}$ by $M_{ \gekl{1,\ti{q}} + \sum_{i,d} [1,i] }$, $\beta_{1,0}^{\ti{q}}$ by 
$
\rukl{
\begin{smallmatrix}
\beta_{1,0}^{\ti{q}} & 0 \\
0 & \id_{\big( \bigoplus_{i,d} F_1^i \big)}
\end{smallmatrix}
}
$
and $\beta_{1,1}^{\ti{q}}$ by 
$
\rukl{
\begin{smallmatrix}
\beta_{1,1}^{\ti{q}} & 0 \\
0 & \Ad(w) \circ \id_{\big( \bigoplus_{i,d} F_1^i \big)}
\end{smallmatrix}
}
$. Now suppose that we have produced 
$A_1 \overset{\varphi_1}{\lori} A_2 \overset{\varphi_2}{\lori} \dotso \overset{\varphi_{n-1}}{\lori} A_n$, and that the next step of modification (path) yields $\varphi_n: \: A_n \to A_{n+1}$. We use the description of $\varphi_F^j \vert_{F_n^i}$ in \eqref{e:varphiFj_NEW}. Let $\fD(j,i)$ be the set of one-dimensional projections in $DM_{m(j,i)}$. For each $d \in \fD(j,i)$, define a permutation $w_{j,i,d} \in F_{n+1}^j$ such that, identifying $d' \otimes \mff$ (for $d' \in \fD(j,i')$ and $\mff \in DF_n^{i'}$) with its image under $\tailarr$ in the description of $\varphi_F^j \vert_{F_n^{i'}}$ in \eqref{e:varphiFj_NEW}, we have $w_{j,i,d} (d \otimes \mff) w_{j,i,d}^* = d \otimes \mff$ and, for all $\hat{d} \in \fD(j,\hat{i})$ with $\hat{d} \neq d$, $w_{j,i,d} (\hat{d} \otimes \mff) w_{j,i,d}^* = \check{d} \otimes \mff$ for some $\check{d} \in \fD(j,\hat{i})$ with $\check{d} \neq \hat{d}$, for all $\mff \in DF_n^{\hat{i}}$. Set $w \defeq (w_{j,i,d})_{j,i,d}$. Choose an index $\ti{q}$ and replace $E_{n+1}^{\ti{q}}$ by $M_{ \gekl{n+1,\ti{q}} + \sum_{j, i, d} [n+1,j] }$, $\beta_{n+1,0}^{\ti{q}}$ by 
$
\rukl{
\begin{smallmatrix}
\beta_{n+1,0}^{\ti{q}} & 0 \\
0 & \id_{\big( \bigoplus_{j,i,d} F_{n+1}^j \big)}
\end{smallmatrix}
}
$
and $\beta_{n+1,1}^{\ti{q}}$ by 
$
\rukl{
\begin{smallmatrix}
\beta_{n+1,1}^{\ti{q}} & 0 \\
0 & \Ad(w) \circ \id_{\big( \bigoplus_{j,i,d} F_{n+1}^j \big)}
\end{smallmatrix}
}
$. Modify $A_{n+1}$ and $\varphi_n$ accordingly as in modification (conn) or the second step of modification (path). Recursive application of this procedure completes modification (sccb).

\blemma
\label{lem:eta-zeta}
After modification (path) combined with modification (sccb), we have the following: For all $\eta \in \cF_1^i$, there exists a continuous path $\xi: \: [0,1] \to G_1$ of the form $\xi(t) = [\omega(t),\gamma]$ with (P1) and (P2) such that $\omega(0) = 0$, $\omega(1) = 1$, $\xi(0) = \eta$, and $\zeta \defeq \xi(1)$ lies in $\cF_1^i$ and satisfies $s(\zeta) = s(\eta)$ but $r(\zeta) \neq r(\eta)$ or $r(\zeta) = r(\eta)$ but $s(\zeta) \neq s(\eta)$. For all $n \geq 1$, $j$ and $\eta \in \cF_{n+1}^j \setminus \cF_{n+1}^j[\bm{p}]$, there exists a continuous path $\xi: \: [0,1] \to G_{n+1}$ of the form $\xi(t) = [\omega(t),\gamma]$ with (P1) and (P2) such that $\omega(0) = 0$, $\omega(1) = 1$, $\xi(0) = \eta$ and $\zeta \defeq \xi(1)$ lies in $\cF_{n+1}^j$ and satisfies $s(\zeta) = s(\eta)$ but $r(\zeta) \neq r(\eta)$ or $r(\zeta) = r(\eta)$ but $s(\zeta) \neq s(\eta)$.
\elemma
\setlength{\parindent}{0cm} \setlength{\parskip}{0cm}

\bproof
Let us start with the first part ($n=1$). Suppose that $s(\eta) = x$ and $r(\eta) = y$ with $x, y \in \cX_1^i$. We think of $x$ corresponding to $d$, $y$ corresponding to $\hat{d}$ and take $\nu$ corresponding to $(i,d)$ in the notation of modification (sccb). Let $\gamma \defeq (\nu,\eta) \in \cE_1^{\ti{q}}$, let $\omega$ be as in (P1) and (P2) with $\omega(0) = 0$ and $\omega(1) = 1$, and define $\xi(t) \defeq [\omega(t),\gamma]$. Then we have $\xi(0) = [0,(\nu,\eta)] = \eta$ as $\bm{b}_{1,0}(\nu,\eta) = \eta$, $s(\xi(1)) = s[1,(\nu,\eta)] = [1,(\nu,x)] = x = s(\eta)$ as $\bm{b}_{1,1}(\nu,x) = x$, but $r(\xi(1)) = r[1,(\nu,\eta)] = [1,(\nu,y)] \neq y = r(\eta)$ as $\bm{b}_{1,1}(\nu,y) \neq y$ by construction.
\setlength{\parindent}{0cm} \setlength{\parskip}{0.5cm}

Now we treat the second part. First suppose that $s(\eta) = (\mu,x)$ and $r(\eta) = (\hat{\mu},y)$ with $\mu \in \cM(j,i)$, $\hat{\mu} \in \cM(j,\hat{i})$, $x \in \cX_n^i$, $y \in \cX_n^{\hat{i}}$. We think of $\mu$ corresponding to $d$, $\hat{\mu}$ corresponding to $\hat{d}$ and take $\nu$ corresponding to $(j,i,d)$ in the notation of modification (sccb). Let $\gamma \defeq (\nu,\eta) \in \cE_{n+1}^{\ti{q}}$, let $\omega$ be as in (P1) and (P2) with $\omega(0) = 0$ and $\omega(1) = 1$, and define $\xi(t) \defeq [\omega(t),\gamma]$. Then we have $\xi(0) = [0,(\nu,\eta)] = \eta$ as $\bm{b}_{n+1,0}(\nu,\eta) = \eta$, $s(\xi(1)) = s[1,(\nu,\eta)] = [1,(\nu,(\mu,x))] = (\mu,x) = s(\eta)$ as $\bm{b}_{n+1,1}(\nu,(\mu,x)) = (\mu,x)$, but $r(\xi(1)) = r[1,(\nu,\eta)] = [1,(\nu,(\hat{\mu},y))] = (\check{\mu},y) \neq (\hat{\mu},y) = r(\eta)$ as $\bm{b}_{n+1,1}(\nu,(\hat{\mu},y)) = (\check{\mu},y)$ for some $\check{\mu} \in \cM(j,\hat{i})$ with $\check{\mu} \neq \hat{\mu}$ by construction. Now suppose that $s(\eta) = (\hat{\mu},x)$ and $r(\eta) = y$ with $\hat{\mu} \in \cM(j,i)$, $x \in \cX_n^i$ and $y \in \cE_n^p \subseteq \cF_{n+1}^j$ (or the other way round, with $s$ and $r$ swapped). We think of $\hat{\mu}$ corresponding to $\hat{d}$ and take $\nu$ corresponding to $(j,i,d)$ for some $d \neq \hat{d}$ in the notation of modification (sccb). Let $\gamma \defeq (\nu,\eta) \in \cE_{n+1}^{\ti{q}}$, let $\omega$ be as in (P1) and (P2) with $\omega(0) = 0$ and $\omega(1) = 1$, and define $\xi(t) \defeq [\omega(t),\gamma]$. Then we have $\xi(0) = [0,(\nu,\eta)] = \eta$ as $\bm{b}_{n+1,0}(\nu,\eta) = \eta$, $r(\xi(1)) = r[1,(\nu,\eta)] = [1,(\nu,y)] = y = r(\eta)$ as $\bm{b}_{n+1,1}(\nu,y) = y$, but $s(\xi(1)) = s[1,(\nu,\eta)] = [1,(\nu,(\hat{\mu},x))] = (\check{\mu},x) \neq (\hat{\mu},x) = s(\eta)$ as $\bm{b}_{n+1,1}(\nu,(\hat{\mu},x)) = (\check{\mu},x)$ for some $\check{\mu} \in \cM(j,\hat{i})$ with $\check{\mu} \neq \hat{\mu}$ by construction.
\eproof
\setlength{\parindent}{0cm} \setlength{\parskip}{0.5cm}

\bprop
\label{prop:eta-zeta_bG}
If we combine modification (path) with modification (sccb) and arrange condition (clsg), then the following holds: Let $\check{\bm{\eta}} \in \bar{G}$, and let $n \geq 0$ be such that $\check{\bm{\eta}} \in \bar{G}_{n+1} \setminus \bar{G}_n$. Suppose that $\bm{p}_{n+1,\infty}(\check{\bm{\eta}})$ is not of the form $[t,\gamma]$ for some $t \in (0,1)$ and $\gamma \in \cE_{n+1}$ with $\gamma \notin \cE_{n+1,0}$ and $\gamma \notin \cE_{n+1,1}$. Then there exist $\bm{\eta}, \bm{\zeta} \in \bar{G}$ such that $\check{\bm{\eta}} \sim_{\rm conn} \bm{\eta} \sim_{\rm conn} \bm{\zeta}$ in $\bar{G}$, and $s(\bm{\zeta}) = s(\bm{\eta})$ but $r(\bm{\zeta}) \neq r(\bm{\eta})$ or $r(\bm{\zeta}) = r(\bm{\eta})$ but $s(\bm{\zeta}) \neq s(\bm{\eta})$.
\eprop
\setlength{\parindent}{0cm} \setlength{\parskip}{0cm}

\bproof
Write $\check{\bm{\eta}} = (\check{\eta}_n)$. We have $\check{\eta}_{n+1} = [t,\gamma]$ for some $\gamma$ in $\cE_{n+1,\mfr}$ for $\mfr = 0 \ \text{or} \ 1$. Construct a path in $G_{n+1}$ connecting $[t,\gamma]$ with $[\mfr,\gamma]$, and, using Proposition~\ref{prop:path}, lift it to a path in $\bar{G}_{n+1}$ connecting $\check{\bm{\eta}}$ to an element $\bm{\eta} = (\eta_n) \in \bar{G}_{n+1}$ with $\eta_{n+1} = [\mfr,\gamma]$. Corollary~\ref{cor:barGn-clopen} implies that $\bar{G}_{n+1} \setminus \bar{G}_n$ is clopen. Since $\check{\bm{\eta}}$ lies in $\bar{G}_{n+1} \setminus \bar{G}_n$ and $\check{\bm{\eta}} \sim_{\rm conn} \bm{\eta}$, it follows that $\bm{\eta}$ lies in $\bar{G}_{n+1} \setminus \bar{G}_n$, too. Therefore, if $n=0$, we must have $\eta_1 = \eta_{n+1} \in \cF_1$, and if $n \geq 1$, we must have $\eta_{n+1} \in \cF_{n+1} \setminus \cF_{n+1}[\bm{p}]$. In both cases, Lemma~\ref{lem:eta-zeta} provides an element $\zeta_{n+1}$ such that $\zeta_{n+1} \in \cF_{n+1}^j$ if $\eta_{n+1} \in \cF_{n+1}^j$, and $s(\zeta_{n+1}) = s(\eta_{n+1})$ but $r(\zeta_{n+1}) \neq r(\eta_{n+1})$ or $r(\zeta_{n+1}) = r(\eta_{n+1})$ but $s(\zeta_{n+1}) \neq s(\eta_{n+1})$, together with a path $\xi_{n+1}$ in $G_{n+1}$ with (P1) and (P2) connecting $\eta_{n+1}$ and $\zeta_{n+1}$. Since $\bm{\eta}$ lies in $\bar{G}_{n+1}$, we must have $\eta_{N+1} = (\mu_N, \eta_N)$ for all $N \geq n+1$. Define $\bm{\zeta} \in \bar{G}$ by setting $\zeta_{N+1} \defeq (\mu_N, \zeta_N)$ for all $N \geq n+1$ and $\bm{\zeta} \defeq (\zeta_N)_{N \geq n+1}$. Then $\bm{\zeta}$ inherits the property from $\zeta_{n+1}$ that $s(\bm{\zeta}) = s(\bm{\eta})$ but $r(\bm{\zeta}) \neq r(\bm{\eta})$ or $r(\bm{\zeta}) = r(\bm{\eta})$ but $s(\bm{\zeta}) \neq s(\bm{\eta})$. Now apply Proposition~\ref{prop:path} recursively to construct paths $\xi_N$, $N \geq n+1$, which connect $\eta_N$ and $\zeta_N$ and satisfy $\bm{p}_N \circ \xi_{N+1} = \xi_N$. It follows that $\xi(t) \defeq (\xi_N(t))_{N \geq n+1}$ defines the desired path connecting $\bm{\eta}$ and $\bm{\zeta}$.
\eproof
\setlength{\parindent}{0cm} \setlength{\parskip}{0.5cm}

We now start to study connected components of $\bar{G}$.
\blemma
\label{lem:FiniteCC}
After modification (path), $\bar{G}_n$ has only finitely many connected components for all $n$.
\elemma
\setlength{\parindent}{0cm} \setlength{\parskip}{0cm}

\bproof
Given $\gamma \in \cE_n$, let $I_{\gamma} \defeq [0,1]$ if $\gamma \in \cE_{n,0}$ and $\gamma \in \cE_{n,1}$, $I_{\gamma} \defeq [0,1)$ if $\gamma \in \cE_{n,0}$ and $\gamma \notin \cE_{n,1}$, $I_{\gamma} \defeq (0,1]$ if $\gamma \notin \cE_{n,0}$ and $\gamma \in \cE_{n,1}$, and $I_{\gamma} \defeq (0,1)$ if $\gamma \notin \cE_{n,0}$ and $\gamma \notin \cE_{n,1}$. Using Proposition~\ref{prop:path} as in the proof of Proposition~\ref{prop:pathconn}, it is straightforward to see that $\bm{p}_{n,\infty}^{-1}[I_\gamma \times \gekl{\gamma}]$ is path-connected in $\bar{G}_n$. Moreover, it is clear that $\bar{G}_n = \bigcup_{\gamma \in \cE_n} \bm{p}_{n,\infty}^{-1}[I_\gamma \times \gekl{\gamma}]$.
\eproof
\setlength{\parindent}{0cm} \setlength{\parskip}{0.5cm}

The following is an immediate consequence of Lemma~\ref{lem:FiniteCC} and Corollary~\ref{cor:barGn-clopen}.
\bcor
If we arrange condition (clsg) in modification (path), then the connected components in $\bar{G}$ are open.
\ecor

\bprop
\label{prop:sccb}
If we combine modification (path) with modification (sccb) and arrange condition (clsg), then the only connected components of $\bar{G}$ which are also bisections (i.e., source and range maps restrict to bijections) are precisely of the form $\bm{p}_{n,\infty}^{-1}[(0,1) \times \gekl{\gamma}] \subseteq \bar{G}_n$ for some $n$ and $\gamma \in \cE_n$ with $\gamma \notin \cE_{n,0}$ and $\gamma \notin \cE_{n,1}$.
\eprop
\setlength{\parindent}{0cm} \setlength{\parskip}{0cm}

\bproof
Let us first show that sets of the form $\bm{p}_{n,\infty}^{-1}[(0,1) \times \gekl{\gamma}]$ for $\gamma$ as in the proposition are indeed connected components and bisections. First of all, an application of Proposition~\ref{prop:path} as in the proof of Proposition~\ref{prop:pathconn}shows that $\bm{p}_{n,\infty}^{-1}[(0,1) \times \gekl{\gamma}]$ is path-connected, hence connected. If $C$ is a connected subset of $\bar{G}$ containing $\bm{p}_{n,\infty}^{-1}[(0,1) \times \gekl{\gamma}]$, then we must have $C \subseteq \bar{G}_n$ because $\bar{G}_n$ is clopen by Corollary~\ref{cor:barGn-clopen}. Moreover, $[(0,1) \times \gekl{\gamma}] \subseteq \bm{p}_{n,\infty}(C)$. As $[(0,1) \times \gekl{\gamma}]$ is a connected component in $G_n$, it follows that $\bm{p}_{n,\infty}(C) \subseteq [(0,1) \times \gekl{\gamma}]$ and thus $C \subseteq \bm{p}_{n,\infty}^{-1}(\bm{p}_{n,\infty}(C)) \subseteq \bm{p}_{n,\infty}^{-1}[(0,1) \times \gekl{\gamma}]$. This shows that $\bm{p}_{n,\infty}^{-1}[(0,1) \times \gekl{\gamma}]$ is a connected component. It is also a bisection: Let $\bm{\eta}, \bm{\zeta} \in \bm{p}_{n,\infty}^{-1}[(0,1) \times \gekl{\gamma}]$ with $s(\bm{\eta}) = s(\bm{\zeta})$. (The case of equal range is analogous.) Write $\bm{\eta} = (\eta_N)_N$, $\bm{\zeta} = (\zeta_N)_N$. It follows that $\eta_n = \zeta_n$. So we have for all $N \geq n$ that $s(\eta_N) = s(\zeta_N)$ and $\bm{p}_{n,N}(\eta_N) = \bm{p}_{n,N}(\zeta_N)$. Since $\bm{p}_{n,N}$ is a fibrewise bijection (see \cite[\S~6.2 and \S~7]{Li18}), it follows that $\eta_N = \zeta_N$ for all $N \geq n$ and hence $\bm{\eta} = \bm{\zeta}$.
\setlength{\parindent}{0cm} \setlength{\parskip}{0.5cm}

Now let $C$ be a connected component in $\bar{G}$. As $\bar{G}_n$ is clopen for each $n$ by Corollary~\ref{cor:barGn-clopen}, there exists $n$ such that $C \subseteq \bar{G}_{n+1}$ and $C \not\subseteq \bar{G}_n$, so that we must have $C \subseteq \bar{G}_{n+1} \setminus \bar{G}_n$. Suppose that $C$ is not of the form $\bm{p}_{n+1,\infty}^{-1}[(0,1) \times \gekl{\gamma}]$ for some $\gamma$ as in the proposition. It follows that $C$ must contain some element $\check{\bm{\eta}}$ as in Proposition~\ref{prop:eta-zeta_bG}, and hence it follows that there exist $\bm{\eta}, \bm{\zeta} \in \bar{G}$ such that $\check{\bm{\eta}} \sim_{\rm conn} \bm{\eta} \sim_{\rm conn} \bm{\zeta}$ in $\bar{G}$, and $s(\bm{\zeta}) = s(\bm{\eta})$ but $r(\bm{\zeta}) \neq r(\bm{\eta})$ or $r(\bm{\zeta}) = r(\bm{\eta})$ but $s(\bm{\zeta}) \neq s(\bm{\eta})$. As $C$ is a connected component, we must have $\bm{\eta}, \bm{\zeta} \in C$. Thus $C$ cannot be a bisection.
\eproof
\setlength{\parindent}{0cm} \setlength{\parskip}{0.5cm}

Let $\fC_{\rm bi}$ be the set of connected components of $\bar{G}$ which are bisections. For two bisections $U$ and $V$ in $\bar{G}$, we define the product $U \cdot V$ only if $s(U) = r(V)$, and in this case $U \cdot V \defeq \menge{\bm{u} \bm{v}}{\bm{u} \in U, \, \bm{v} \in V}$ is another bisection. Let $\spkl{\fC_{\rm bi}}$ be the smallest collection of bisections in $\bar{G}$ closed under products and containing $\fC_{\rm bi}$, i.e., the set of all finite products of elements in $\fC_{\rm bi}$.
\blemma
\label{lem:<Cbi>}
If we combine modification (path) with modification (sccb) and arrange condition (clsg) as well as
\begin{equation}
\label{e:np>4ni}
\gekl{n,p} > 4 [n,i] \quad \forall \ n,p,i,
\end{equation}
then $\spkl{\fC_{\rm bi}} = \menge{\bm{p}_{n,\infty}^{-1}[(0,1) \times \gekl{\gamma}]}{n \in \Zz_{\geq 1}, \, \gamma \in \cE_n}$.
\elemma
\setlength{\parindent}{0cm} \setlength{\parskip}{0cm}

\bproof
Clearly, $\menge{\bm{p}_{n,\infty}^{-1}[(0,1) \times \gekl{\gamma}]}{n \in \Zz_{\geq 1}, \, \gamma \in \cE_n}$ is a collection of bisections in $\bar{G}$ closed under products and containing $\fC_{\rm bi}$. This shows \an{$\subseteq$}. To prove \an{$\supseteq$}, we show that for all $n, p$ and $\gamma \in \cE_n^p$ with $s(\gamma) = y^s$ and $r(\gamma) = y^r$ that $\bm{p}_{n,\infty}^{-1}[(0,1) \times \gekl{\gamma}] \in \spkl{\fC_{\rm bi}}$. Recall that $\bm{b}_{n,\bullet}^p$ is given by a composition of the form $\cE_{n,\bullet}^p \isom \coprod_i (\cM_\bullet(p,i) \times \cF_n^i) \onto \coprod_i \cF_n^i = \cF_n$ (see \eqref{e:betanpi}). Let us denote the induced bijections $\cY_{n,\bullet}^p \isom \coprod_i (\cM_\bullet(p,i) \times \cX_n^i)$ by $y \ma (y)_\bullet$. Suppose that $(y^*)_\bullet = (\mu_\bullet^*,x_\bullet^*)$, with $x_\bullet^* \in \cX_n^{i_\bullet^*}$, for $* = r,s$ and $\bullet = 0,1$. We claim that there exists $y \in \cY_n^p$ such that $(y)_\bullet \notin \gekl{\mu_\bullet^*} \times \cX_n^{i_\bullet^*}$ for all $* = r,s$ and $\bullet = 0,1$. Indeed, since $\# \big\{ y \in \cY_n^p: \: (y)_\bullet \in \gekl{\mu_\bullet^*} \times \cX_n^{i_\bullet^*} \big\} = \# \cX_n^{i_\bullet^*} = [n,i_\bullet^*]$ for all $* = r,s$ and $\bullet = 0,1$, we have that $\# \big\{ y \in \cY_n^p: \: (y)_\bullet \in \gekl{\mu_\bullet^*} \times \cX_n^{i_\bullet^*} \ \text{for some} \ * = r,s, \bullet = 0,1 \big\} = [n,i_0^s] + [n,i_0^r] + [n,i_1^s] + [n,i_1^r] \leq 4 \max \menge{[n,i_\bullet^*]}{* = r,s, \bullet = 0,1} < \gekl{n,p} = \# \cY_n^p$ by condition \eqref{e:np>4ni}. Now take $y$ with these properties, and let $\gamma_1, \gamma_2 \in \cE_n^p$ be such that $s(\gamma_1) = y^s$, $r(\gamma_1) = y$, $s(\gamma_2) = y$, $r(\gamma_2) = y^r$. Then $\gamma_1, \gamma_2 \notin \cE_{n,0}$ and $\gamma_1, \gamma_2 \notin \cE_{n,1}$, so that $\bm{p}_{n,\infty}^{-1}[(0,1) \times \gekl{\gamma_i}] \in \fC_{\rm bi}$ for $i=1,2$. It follows that $\bm{p}_{n,\infty}^{-1}[(0,1) \times \gekl{\gamma}] = \bm{p}_{n,\infty}^{-1}[(0,1) \times \gekl{\gamma_2}] \cdot \bm{p}_{n,\infty}^{-1}[(0,1) \times \gekl{\gamma_1}] \in \spkl{\fC_{\rm bi}}$.
\eproof
\setlength{\parindent}{0cm} \setlength{\parskip}{0.5cm}

\bdefin[compare {\cite[Definition~3.1]{Nek}}]
Let $\cY$ be a finite set. We call a multisection in $\spkl{\fC_{\rm bi}}$ the image of an injective map $\cY \times \cY \to \spkl{\fC_{\rm bi}}, \, (x,y) \ma U_{x,y}$ such that $U_{x,y} \cdot U_{y',z}$ is only defined if $y = y'$, and in that case $U_{x,y} \cdot U_{y,z} = U_{x,z}$. We call $\# \cY$ the degree of $\gekl{U_{x,y}}$.
\edefin
\setlength{\parindent}{0cm} \setlength{\parskip}{0cm}

\bcor
In the situation of Lemma~\ref{lem:<Cbi>}, multisections in $\spkl{\fC_{\rm bi}}$ are precisely of the form $\bm{p}_{n,\infty}^{-1}[(0,1) \times \gekl{\gamma}]$ for some $n, p$ and $\gamma \in \cE_n^p$.
\ecor
\setlength{\parindent}{0cm} \setlength{\parskip}{0.5cm}

As it is clear that the degree can be read off from a multisection, Lemma~\ref{lem:<Cbi>} implies the following.
\bcor
\label{cor:GG'}
Suppose that we combine modification (path) with modification (sccb) and arrange (clsg), \eqref{e:np>4ni} to obtain classifiable C*-algebras $A$ and $A'$ with the same prescribed Elliott invariant together with C*-diagonals $B$ and $B'$. Let $\bar{G}$ and $\bar{G}'$ be the groupoid models for $(A,B)$ and $(A',B')$ and $\cY_n^p(\bar{G})$, $\cY_n^p(\bar{G}')$ the analogues of $\cY_n^p$ above for $\bar{G}$, $\bar{G}'$. If $\bar{G} \cong \bar{G}'$ as topological groupoids (i.e., if $(A,B) \cong (A',B')$), then we must have $\gekl{\# \cY_n^p(\bar{G})}_{n,p} = \gekl{\# \cY_n^p(\bar{G}')}_{n,p}$.
\ecor
Let us now construct, for every sequence $\mfm = (\mfm_n)$ of non-negative integers a groupoid model $\bar{G}(\mfm)$ for our classifiable C*-algebra such that for any two sequences $\mfm$ and $\mfn$, we have $\bar{G}(\mfm) \not\cong \bar{G}(\mfn)$ if $\mfm \neq \mfn$. First combine modification (path) with modification (sccb) and arrange (clsg), \eqref{e:np>4ni} to obtain a classifiable C*-algebra $A$ with prescribed Elliott invariant $\cE$ as in \S~\ref{ss:(path)} and C*-diagonal $B$. Now we modify the construction. For all $n \geq 1$, choose a direct summand $F_n^j$ of $F_n$ such that, for all $n \geq 1$, we have $F_{n+1}^j \neq F_{n+1}^{j_\mfr^p}$ for all $p$ and $\mfr = 0,1$. Given a sequence $\mfm = (\mfm_n)$, we modify $(A,B)$ by adding $\id_{(F_n^{j_n})^{\oplus \mfm_n}}$ to $\beta_{n,\bullet}^p$ for all $p$ and enlarging $E_n^p$ correspondingly. In this way, we obtain for each $\mfm$ a classifiable C*-algebra $A(\mfm)$ with the same prescribed Elliott invariant $\cE$ and the same properties as $A$, together with a C*-diagonal $B(\mfm)$ of $A(\mfm)$. Let $\bar{G}(\mfm)$ be the groupoid model of $(A(\mfm),B(\mfm))$.
\bprop
If $\mfm \neq \mfn$, then $\bar{G}(\mfm) \not\cong \bar{G}(\mfn)$, i.e., $(A(\mfm),B(\mfm)) \not\cong (A(\mfn),B(\mfn))$.
\eprop
\setlength{\parindent}{0cm} \setlength{\parskip}{0cm}

\bproof
Let $\bar{G} \defeq \bar{G}(\mfm)$ and $\bar{G}' \defeq \bar{G}(\mfn)$. Suppose that $\mfm_n = \mfn_n$ for all $n \leq N-1$ and that $\mfm_N \neq \mfn_N$, say $\mfm_N < \mfn_N$. As the first $N-1$ steps of the construction coincide, we have $\menge{\# \cY_n^p(\bar{G})}{n \leq N-1,p} = \menge{\# \cY_n^p(\bar{G}')}{n \leq N-1,p}$. Now we have $\# \cY_N^p(\bar{G}') = \# \cY_N^p(\bar{G}) + (\mfn_N - \mfm_N) \cdot \# \cX_N^{j_N}$ for all $p$. Hence $\# \cY_N^{p'}(\bar{G}') > \min_p \# \cY_N^p(\bar{G})$ for all $p'$. As $\# \cY_{\bar{n}+1}^{q'}(\bar{G}') > \# \cY_{\bar{n}}^{p'}(\bar{G}')$ for all $\bar{n}$, $q'$ and $p'$ , it follows that $\min_p \# \cY_N^p(\bar{G})$ does not appear in $\gekl{\# \cY_n^p(\bar{G}')}_{n,p}$, while it appears in $\gekl{\# \cY_n^p(\bar{G})}_{n,p}$. Hence Corollary~\ref{cor:GG'} implies that $\bar{G} \not\cong \bar{G}'$, i.e., $(A(\mfm),B(\mfm)) \not\cong (A(\mfn),B(\mfn))$.
\eproof
\setlength{\parindent}{0cm} \setlength{\parskip}{0.5cm}

All in all, in combination with Corollaries~\ref{cor:Menger_unital} and \ref{cor:M-C_spl}, we obtain
\btheo
\label{thm:ManyMenger_GPD_Ell}
For every sequence $\mfm$ in $\Zz_{\geq 0}$ and every prescribed Elliott invariant $(G_0, G_0^+, u, T, r, G_1)$ as in \cite[Theorem~1.2]{Li18} with torsion-free $G_0$ and trivial $G_1$, our construction produces topological groupoids $\bar{G}(\mfm)$ with the same properties as in \cite[Theorem~1.2]{Li18} (in particular, $C^*_r(\bar{G}(\mfm))$ is a classifiable unital C*-algebra satisfying ${\rm Ell}(C^*_r(\bar{G}(\mfm))) \cong (G_0, G_0^+, u, T, r, G_1)$), such that $\bar{G}(\mfm)^{(0)} \cong \bm{M}$, and $\bar{G}(\mfm) \not\cong \bar{G}(\mfn)$ if $\mfm \neq \mfn$.

For every sequence $\mfm$ in $\Zz_{\geq 0}$ and every prescribed Elliott invariant $(G_0, T, \rho, G_1)$ as in \cite[Theorem~1.3]{Li18} with torsion-free $G_0$ and trivial $G_1$, our construction produces topological groupoids $\bar{G}(\mfm)$ with the same properties as in \cite[Theorem~1.3]{Li18} (in particular, $C^*_r(\bar{G}(\mfm))$ is a classifiable stably projectionless C*-algebra with continuous scale satisfying ${\rm Ell}(C^*_r(\bar{G}(\mfm))) \cong (G_0, \gekl{0}, T, \rho, G_1)$), such that $\bar{G}(\mfm)^{(0)} \cong \bm{M}_{\setminus C}$, and $\bar{G}(\mfm) \not\cong \bar{G}(\mfn)$ if $\mfm \neq \mfn$.
\etheo

In combination with the classification result in \cite{Rob}, this yields the following
\btheo
\label{thm:ManyMenger_Diag_Ell}
For every prescribed Elliott invariant $(G_0, G_0^+, u, T, r, G_1)$ as in \cite[Theorem~1.2]{Li18} with torsion-free $G_0$ and trivial $G_1$, our construction produces a classifiable unital C*-algebra $A$ with ${\rm Ell}(A) \cong (G_0, G_0^+, u, T, r, G_1)$ and continuum many pairwise non-conjugate C*-diagonals of $A$ whose spectra are all homeomorphic to $\bm{M}$.

For every Elliott invariant $(G_0, T, \rho, G_1)$ as in \cite[Theorem~1.3]{Li18} with torsion-free $G_0$ and $G_1 = \gekl{0}$, our construction produces a classifiable stably projectionless C*-algebra $A$ having continuous scale with ${\rm Ell}(A) \cong (G_0, \gekl{0}, T, \rho, G_1)$ and continuum many pairwise non-conjugate C*-diagonals of $A$ whose spectra are all homeomorphic to $\bm{M}_{\setminus C}$.
\etheo
\setlength{\parindent}{0cm} \setlength{\parskip}{0cm}

This theorem, combined with classification results for all classifiable C*-algebras, implies Theorems~\ref{thm:main2_unital} and \ref{thm:main2_spl}.
\setlength{\parindent}{0cm} \setlength{\parskip}{0.5cm}

\end{document}